\documentclass[11pt]{article}
\usepackage{amsmath}
\usepackage{amssymb,amsmath}
\usepackage{mathrsfs}
\usepackage{graphicx}
\usepackage{color}
\usepackage{esint}%
\usepackage[linktoc=page]{hyperref}
\usepackage[numbers,sort&compress]{natbib}
\usepackage{tocloft}
\usepackage[title]{appendix}

\def\dive{\mathop{\rm div}\nolimits}

\def\C{\mathop{\bf C\kern 0pt}\nolimits}
\def\DD{\mathop{\bf D\kern 0pt}\nolimits}
\def\K{\mathop{\bf K\kern 0pt}\nolimits}
\def\N{\mathop{\bf N\kern 0pt}\nolimits}
\def\Q{\mathop{\bf Q\kern 0pt}\nolimits}
\def\R{\mathop{\bf R\kern 0pt}\nolimits}

\newcommand{\beq}{\begin{equation}}
\newcommand{\eeq}{\end{equation}}
\newcommand{\ben}{\begin{eqnarray}}
\newcommand{\een}{\end{eqnarray}}
\newcommand{\beno}{\begin{eqnarray*}}
\newcommand{\eeno}{\end{eqnarray*}}

\newtheorem{Theorem}{Theorem}[section]
\newtheorem{Definition}[Theorem]{Definition}
\newtheorem{Proposition}[Theorem]{Proposition}
\newtheorem{Lemma}[Theorem]{Lemma}
\newtheorem{Corollary}[Theorem]{Corollary}
\newtheorem{Remark}[Theorem]{Remark}
\newtheorem{Question}{Question}

\numberwithin{equation}{section}

\allowdisplaybreaks \numberwithin{equation} {section}

\begin{document}
\title{Sharp non-uniqueness for the 2D hyper-dissipative Navier-Stokes equations  }
      \author{
       \\[2mm]
{\small  Lili Du  \ \ and \  \ Xinliang Li}}
         \date{}
         \maketitle
\noindent{\bf Abstract.} In this article, we study the non-uniqueness of weak solutions for the two-dimensional hyper-dissipative Navier-Stokes equations in the super-critical spaces $L_{t}^{\gamma}W_{x}^{s,p}$ when $\alpha\in[1,\frac{3}{2})$, and obtain the conclusion that the non-uniqueness of the weak solutions at the two endpoints is sharp in view of the generalized Lady\v{z}enskaya-Prodi-Serrin condition with the triplet $(s,\gamma,p)=(s,\infty, \frac{2}{2\alpha-1+s})$ and $(s, \frac{2\alpha}{2\alpha-1+s}, \infty)$. As a good observation, we use the intermittency of the temporal concentrated function in an almost optimal way. The research results extend the recent elegant works on 2D Navier-Stokes equations in [Cheskidov and Luo, Invent.  Math.,  229 (2022), pp. 987--1054; Cheskidov and Luo, Ann. PDE, 9:13 (2023)] to the hyper-dissipative case $\alpha \in(1,\frac{3}{2})$, and are also applicable in Lebesgue and Besov spaces. It is proved that even in the case of high viscosity, the behavior of the solution remains unpredictable and stochastic due to the lack of integrability and regularity.\ \ \

\vskip   0.2cm \footnotetext[1]{
College of Mathematical Sciences, Shenzhen University, Shenzhen 518061, P. R. China; Department of Mathematics, Sichuan University, Chengdu 610064, P. R. China. Email address, Lili Du: dulili@szu.edu.cn.}

\vskip   0.2cm \footnotetext[2]{
College of Mathematical Sciences, Shenzhen University, Shenzhen 518061, P. R. China; College of Physics and
Optoelectronic Engineering, Shenzhen University, Shenzhen 518061, P. R. China. Email address, Xinliang Li: lixinliangmaths@163.com.}

	\setlength{\baselineskip}{20pt}

\begin{center}
\tableofcontents
\end{center}

\begin{center}
\section{Introduction and main results}\label{A}
\end{center}
\subsection{Background.}\ \ \ \
In this paper, we consider the Cauchy problem of the following 2D generalized Navier-Stokes equations (gNSE) on the torus $\mathbb{T}^{2}=\mathbb{R}^{2}/\mathbb{Z}^{2}$:
\[
\left\{\begin{array}{l}
\partial_{t} u+\nu(-\Delta)^{\alpha} u+(u \cdot \nabla) u+\nabla\mathsf{P}=0 , \label{1.1}\tag{1.1}\\
\operatorname{div} u=0,\\
u|_{t=0}=u_{0},
\end{array}\right.
\]
where $u:\mathbb{T}^{2}\times[0,T] \rightarrow\mathbb{R}^{2}$ is the velocity field and $\mathsf{P}: \mathbb{T}^{2}\times[0,T]\rightarrow \mathbb{R}$ is the scalar pressure of the fluid, $\nu>0$ is the constant viscous coefficient, $\alpha$ is called viscosity exponent in the literature, and $(-\Delta)^{\alpha}$ is the fractional Laplacian defined by the Fourier transform on the flat torus as
$$
\mathcal{F}\left((-\Delta)^{\alpha} u\right)(\xi)=|\xi|^{2 \alpha} \mathcal{F}(u)(\xi), \quad \xi \in \mathbb{Z}^{2}.
$$
Sometimes, we write $\Lambda=(-\Delta)^{1 / 2}$ for notational convenience. When $\alpha < 1$, \eqref{1.1} is known as hypo-dissipative Navier-Stokes equations, and  $\alpha > 1$ is hyper-dissipative Navier-Stokes equations. When $\alpha = 1$, \eqref{1.1} turns to be the 2D classical Navier-Stokes equations (NSE), for which the classical conclusions \cite{ OA58, JL33,JL34,RT77} indicate that the 2D Leray-Hopf solution is the unique smooth solution. Furthermore, J.-L. Lions \cite{JL69} firstly proved that for any $d\geq2$ dimensional gNSE \eqref{1.1}, the Leray-Hopf solution is unique when $\alpha \geq \frac{d+2}{4}$, which also guarantees global regularity with smooth initial data presented in Wu \cite{JW03}.

For the 2D hyper-dissipative NSE \eqref{1.1}, we mainly consider the mixed Sobolev space, and the space $\mathbb{X}=L_{t}^{\gamma} \dot{W}_{x}^{s, p}$ is called critical if its norm $\|\cdot\|_{\mathbb{X}}$ is invariant under the  natural scaling:
\begin{equation*}
u(t, x) \mapsto \lambda^{2 \alpha-1} u\left(\lambda^{2 \alpha} t, \lambda x\right), \quad P(t, x) \mapsto \lambda^{4 \alpha-2} P\left(\lambda^{2 \alpha} t, \lambda x\right), \label{1.4}\tag{1.2}
\end{equation*}
where the exponents $(s, \gamma, p)$ satisfy
\begin{equation*}
\frac{2 \alpha}{\gamma}+\frac{2}{p}=2 \alpha-1+s \label{1.5}\tag{1.3}
\end{equation*}
for some $\gamma, p\in[1,\infty]$, $s \geq 0$, known as generalized Lady\v{z}enskaya-Prodi-Serrin (gLPS) condition.

In particular, the mixed Lebesgue space $L_{t}^{\gamma} L_{x}^{p}$ is called critical for the classical NSE ($\alpha=1$), when the exponents $(\gamma, p)$ satisfy the well-known Lady\v{z}enskaya-Prodi-Serrin (LPS) condition $\frac{2}{\gamma}+\frac{2}{p}=1$(see \cite{PRO59,SER62,LAD67}, for the case $d=2$). When $\frac{2}{\gamma}+\frac{2}{p}<1$,  the space $L_{t}^{\gamma} L_{x}^{p}$ is called sub-critical and super-critical when $\frac{2}{\gamma}+\frac{2}{p}>1$.
As far as we know, the current well-posedness conclusions are all discussed in (sub)critical space, as outlined below.

\subsection{Summary of existence and uniqueness results.}\ \ \ \
The research on the well-posedness of 2D NSE in certain critical spaces originated from the pioneering work of Fujita-Kato \cite{HF64}, in which they first introduced scaling analysis to determine the workspace of the NSE, and used Banach's fixed point theorem to consider the existence and uniqueness of the mild solutions. Some similar well-posedness results can be found in \cite{YG85,TK84} for the critical space $L_{x}^{2}$, in \cite{IG02} for the scaling invariant Besov space $B_{p, q}^{2 / p-1}$ with $p, q < \infty$, and in \cite{HK01} for the well-known critical space $\mathrm{BMO}^{-1}$. Moreover, Lei and Lin \cite{ZL11} innovatively established the global well-posedness with large initial data in $\chi^{-1}$ which is a subspace  of $\mathrm{BMO}^{-1}$. For a comprehensive overview of related works, we refer the readers to the monographs \cite{MC04,PG16,YM97}. According to the classical embedding relationship of critical spaces for $p, q < \infty$,
\begin{equation*}
L^{2}_{x} \hookrightarrow \dot{B}_{p, q}^{\frac{2}{p}-1} \hookrightarrow B M O^{-1} \hookrightarrow \dot{B}_{\infty, \infty}^{-1}, \label{1.7}\tag{1.4}
\end{equation*}
it can be seen that the existence and uniqueness of weak solutions to the 2D NSE has been relatively well established, excluding the endpoint $\dot{B}_{\infty, \infty}^{-1}$ which represents the largest critical space in scaling the NSE, regardless of spatial dimensions. Current research in this space necessitates additional assumptions, as evident in references \cite{JY06, JY11, MP14} and those cited within. Furthermore, the uniqueness of weak solutions within $C_t \dot{B}_{\infty, \infty}^{-1}$ remains a long standing problem, despite the availability of results on other forms of ill-posedness in dimensions $d \geq 3$, such as \cite{ JB08,PG08,TY10,BW15}, and for the hyperdissipative NSE in $B_{\infty, \infty}^{1-2\alpha}$ with $\alpha\in[1,\frac{5}{4})$ by Cheskidov and Schvydkoy in \cite{AC10}, as well as norm-infation instability in $\dot{B}_{\infty, q}^{-s}$ with $s\geq\alpha\geq1,\ q\in(2,\infty]$ explored by Cheskidov and Dai \cite{AC14}.

The known uniqueness results can be concisely summarized as stated in Theorem 1.3 of \cite{AC22}. Specifically, the significant findings reported in \cite{EB72,GF00,PL01} imply that any weak solution to the 2D NSE in the distributional sense, as defined in Definition \ref{A.1} below, within the (sub)critical spaces $L_{t}^{\gamma} L_{x}^{p}$ with the condition $\frac{2}{\gamma}+ \frac{2}{p} \leq 1$ for some $\gamma, p \in [1, \infty]$, is automatically the unique regular Leray-Hopf solution on $\mathbb{T}^{2}$. However, when considering the 2D hyper-dissipative NSE $(\alpha>1)$, there are few literatures discussing the uniqueness of weak solutions, such as \cite{ST98} for the one endpoint space $C_{t} L_{x}^{\frac{2}{2\alpha-1}}$ when $s=0$ in view of the gLPS condition \eqref{1.5}, and \cite{ST03} for the Besov space $\dot{B}_{p, q}^{\frac{2}{p}+1-2 \alpha}$. In fact, even when the viscosity index greater than ``1'' which is called 2D Lions exponent, only enough integrability and regularity can guarantee the uniqueness of the solution. Specifically, we summarize in Theorem \ref{I.1} of Appendix \ref{G2} as follows: For some $\gamma, p\in[1,\infty]$, $s \geq 0$ and $0 < \frac{1}{p}-\frac{s}{2} <\frac{1}{2}$, when \begin{equation*}
\frac{2\alpha}{\gamma}+\frac{2}{p} \leq 2\alpha-1+s,\label{1.41}\tag{1.5}
\end{equation*}
all weak solutions (in the sense of Definition \ref{A.1}) belonging to the class $X^{s,\gamma,p}$ are automatically Leray-Hopf and hence, by the celebrated conclusions \cite{JL69,JW03} when $\alpha\geq1$, unique and regular. In other words, within the scale of $X^{s,\gamma, p}$ spaces, sub-critical or critical weak solutions are classical solutions.
\subsection{Statement of the problem.}\ \ \ \
In the super-critical space, there remains many non-uniqueness questions unsolved, offering a stark contrast to the positive advancements. A recent remarkable contribution by Cheskidov and Luo \cite{AC22} has proved that the sharp non-uniqueness of the weak solutions to the 2D NSE at the one endpoint $(s, \gamma, p) = (0,2, \infty)$ in view of the LPS condition by firstly introducing the temporal concentrated function into the convex integration framework. Additional references \cite{AC21,ACX22} show the validity of the temporal intermittency for the transport equations, while \cite{YLZ22,YLZZ22} explores its implications for 3D Magneto-Hydrodynamic (MHD) equations. Recently, Li, Qu, Zeng and Zhang \cite{YL22} have proved that $L_{t}^{2}C_{x}$ weak solution is non-unique for the 2D hypo-dissipative NSE $(\alpha<1)$ which implies that $\alpha = 1$ is the sharp viscosity threshold for the well-posedness in the space $L_{t}^{2}C_{x}$. Regarding the alternative endpoint $(s, \gamma, p) = (0, \infty, 2)$, Luo and Qu \cite{LTP20} have presented that $C_{t}L_{x}^{2}$ weak solutions lack uniqueness when $\alpha<1$. Very recently, Cheskidov and Luo \cite{ACL22} proved that it is not feasible to construct non-unique $C_{t}L_{x}^{2}$ weak solutions for the 2D NSE, indicating that $\alpha = 1$ serves as the sharp viscosity threshold for the well-posedness within the space $C_{t}L_{x}^{2}$.

As concluded above, the current researches only show the non-unique results for \eqref{1.1} with $\alpha\leq1$ in the super-critical space, whether there exist non-unique results for $\alpha>1$ in 2D case remains unknown. Therefore, we deliberate upon the following two questions  pertaining to non-uniqueness of weak solutions for the 2D hyper-dissipative NSE:
\begin{Question}
In the hyper-dissipative case $\alpha >1$,  is it feasible to find non-unique and non-Leray-Hopf weak solutions, despite having the same initial data as Leray-Hopf solution?\label{Q1}
\end{Question}
\begin{Question}
In view of the gLPS condition \eqref{1.5} for some $\gamma, p\in[1,\infty]$, $s \geq 0$, do non-unique weak solutions exist for the 2D hyper-dissipative NSE \eqref{A.1} in super-critical spaces $L_{t}^{\gamma} W_{x}^{s, p}$, where the condition $2 \alpha / \gamma + 2 / p > 2 \alpha - 1 + s$ holds?\label{Q2}
\end{Question}

In this study, we provide positive answers to the Question \ref{Q1} and Question \ref{Q2} at the two endpoints $\gamma=\infty$ and $p=\infty$. These results are included in Theorem \ref{A.2} below, which addresses the non-uniqueness of weak solutions in the super-critical spaces  $L_{t}^{\gamma} W_{x}^{s, p}$ with $(s, \gamma, p)\in\mathcal{A}_{1}\cup\mathcal{A}_{2}$, when $\alpha\in[1,\frac{3}{2})$. Consequently, we recover the recent elegant works \cite{AC22,ACL22} when $\alpha=1,\ s=0$ and extend to the hyper-dissipative case $\alpha \in(1,\frac{3}{2})$ in Theorem \ref{A.2} by using the temporal intermittency in an almost optimal way.

\subsection{Convex integration technique.}\ \ \ \
Convex integration technique was originally introduced by Nash \cite{JO54} and Kuiper \cite{NI55} in geometry for isometric embeddings, and was introduced into fluid dynamics by De Lellis and  Sz\'{e}kelyhidi Jr. in their groundbreaking work \cite{CD09}. It is worth noting that the convex integration method has shown to be extremely beneficial to the fluids community and its evolution across a series of works \cite{ TB15, TB16,CD09, CD13, CD14} ultimately leading to the resolution of the flexible part of Onsager's conjecture for the 3D Euler equations by Isett \cite{PI18} and Buckmaster-De Lellis-Sz\'{e}kelyhidi-Vicol \cite{TB19} for dissipative solutions. See also the work \cite{Nov23} by Novack and Vicol on an intermittent Onsager theorem, building upon their previous collaboration \cite{BMNV23} with Buckmaster and Masmoudi. Additionally, Giri, Kwon and Novack \cite{VG23} constructed the $L^3$-based strong Onsager conjecture. Recently, Giri and Radu \cite{VGR23} settled the flexible part of the Onsager conjecture in 2D by a Newton-Nash convex integration method. For the related issue on anomalous dissipation, please refer to the recent study \cite{DR24} by De Rosa and Park.  Very recently, Bulut, Huynh and Palasek \cite{BU23} established that the regularity in 3D can be elevated to $C^{1 /2-\varepsilon}$ in the presence of force. For the applications in other models, please refer to \cite{Shv11, IV15, BSV19, MS20, Nov20, BBV20, FLS21} and surveys \cite{DLS19, BV21} for additional insightful works on convex integration framework.

When turning to consider the viscous fluids, Buckmaster and Vicol \cite{BV19} creatively constructed non-unique $C_{t}L_{x}^{2}$ weak solutions of the 3D NSE by using the spatial intermittent convex integration. Following this significant advancement, there have been many subsequent results including extensions of \cite{BV19} to the hyper-dissipative case \cite{LT20, YL22}, partial regularity in time \cite{BCV22}, high dimensions $d \geq 4$ \cite{Luo19}, MHD equations \cite{YLZ22,YLZZ22} and Hall-MHD equations \cite{Dai21} and so on. For the non-uniqueness of Leray-Hopf solutions to 3D hypo-dissipative NSE please refer to \cite{CDLDR18, DR19}.

In 2D, how to set up different building blocks in the convex integration framework without interfering or intersecting with each other is one of the key considerations. Additionally, in order to achieve the hyper-dissipative case $(\alpha>1)$, this study also requires the use of spatial-temporal intermittent convex integration method, and will focus on developing the crucial role of temporal concentration functions $g_{(k)}$ so that even if the different spatial building blocks $\mathbf{W}_{k}$ intersects with each other, the newly generated ``building blocks'' $g_{(k)}\mathbf{W}_{k}$ will not intersect under the effect of non-intersecting temporal building blocks by temporal shifts $t_{k}$, which will be particularly reflected in the discussion of the super-critical regime $\mathcal{A}_{2}$ where we are not necessary to consider the errors caused by support intersection in \cite{AC21}. And we will utilize temporal intermittency in an almost optimal way to achieve the sharp non-uniqueness at the two endpoint cases of gLPS condition, namely $(s,\infty, \frac{2}{2\alpha-1+s})$ and $(s, \frac{2\alpha}{2\alpha-1+s}, \infty)$. As will be demonstrated subsequently, in the case of extreme dissipativity where $\alpha$ approximates $\frac{3}{2}$, the appropriate temporal intermittency is approximately equivalent to 2D spatial intermittency in the super-critical regimes $\mathcal{A}_{1}$ defined in \eqref{1.8}, and correspondingly, 3D spatial intermittency in the super-critical regimes $\mathcal{A}_{2}$ defined in \eqref{1.9}.

\subsection{Notations.}\ \ \ \
The notation $a \lesssim b$ means that $a \leq C b$ for some non-negative constant $C$. The mean of $u \in L^1\left(\mathbb{T}^2\right)$ is given by $\fint_{\mathbb{T}^2} u d x=\left|\mathbb{T}^2\right|^{-1} \int_{\mathbb{T}^2} u d x$, where $|\cdot|$ denotes the Lebesgue measure. For $(p, \gamma) \in[1,+\infty]$ and $s \in \mathbb{R}$, we use the following shorthand notations,
$$
L_{t}^{p}:=L^{p}(0, T), \quad L_{x}^{p}:=L^{p}\left(\mathbb{T}^{2}\right), \quad H_{x}^{s}:=H^{s}\left(\mathbb{T}^{2}\right), \quad W_{x}^{s, p}:=W^{s, p}\left(\mathbb{T}^{2}\right)
$$
where $W_{x}^{s, p}$ is the usual Sobolev space and $H_{x}^{s}=W_{x}^{s, 2}$, $L_{t}^{\gamma} L_{x}^{p}$ denotes the usual Banach space $L^{\gamma}\left(0, T ; L^{p}\left(\mathbb{T}^{2}\right)\right)$. In particular, we write $L_{t, x}^{p}:=L_{t}^{p} L_{x}^{p}$ for brevity. Let
$$
\|u\|_{W_{t, x}^{N, p}}:=\sum_{0 \leq m+|\xi| \leq N}\left\|\partial_{t}^{m} \nabla^{\xi} u\right\|_{L_{t, x}^{p}}, \quad\|u\|_{C_{t, x}^{N}}:=\sum_{0 \leq m+|\xi| \leq N}\left\|\partial_{t}^{m} \nabla^{\xi} u\right\|_{C_{t, x}}
$$
where $\xi=\left(\xi_{1}, \xi_{2} \right)$ is the multi-index and $\nabla^{\xi}:=\partial_{x_{1}}^{\xi_{1}} \partial_{x_{2}}^{\xi_{2}}$. For any Banach space $X$, we use $C([0, T]; X)$ to denote the space of continuous functions from $[0, T]$ to $X$ with norm $\|u\|_{C_{t} X}:=\sup _{t \in[0, T]}\|u(t)\|_{X}$,
and $B_{p, q}^{s}\left(\mathbb{T}^{2}\right)$ to denote the Besov space with the norm
$$
\|u\|_{B_{p, q}^{s}\left(\mathbb{T}^{2}\right)}=\left(\sum_{j \geq -1}\left|2^{j s}\left\|\Delta_{j} u\right\|_{L^{p}\left(\mathbb{T}^{2}\right)}\right|^{q}\right)^{\frac{1}{q}}
$$
where $\left\{\Delta_{j}\right\}_{j \in \mathbb{Z}}$ is the inhomogeneous frequency localization operators. We denote the Triebel-Lizorkin space $F_{p, q}^{s}\left(\mathbb{T}^{2}\right)$ with the norm
$$
\|u\|_{F_{p, q}^{s}\left(\mathbb{T}^{2}\right)}=\left\|\left(\sum_{j \geq -1}\left|2^{j s} \Delta_{j} u\right|^{q}\right)^{\frac{1}{q}}\right\|_{L^{p}\left(\mathbb{T}^{2}\right)}.
$$
The homogeneous space Besov $\dot{B}_{p, q}^{s}\left(\mathbb{T}^{2}\right)$ and Triebel-Lizorkin space $\dot{F}_{p, q}^{s}\left(\mathbb{T}^{2}\right)$ are defined similarly please refer to \cite{HS87} for more details.

For any $\mathscr{A} \subseteq[0, T]$, we introduce the $\varepsilon_{*}$-neighborhood of $\mathscr{A}$ as
$$
\mathscr{N}_{\varepsilon_{*}}(\mathscr{A}):=\left\{t \in[0, T]: \exists\ s \in \mathscr{A} \text {, s.t. }|t-s| \leq \varepsilon_{*}\right\},
$$
where $\varepsilon_{*}>0$. Additionally, we set the viscosity coefficient $\nu=1$ for brevity.

\subsection{Main results.}\ \ \ \
Before presenting the main results, we introduce the concept of weak solutions in the distributional sense for the system \eqref{1.1}.

\begin{Definition} \label{A.1}
(Weak solution) For $\alpha\geq1$, given any weakly divergence-free datum $u_{0} \in L^{2}\left(\mathbb{T}^{2}\right)$, $u \in$ $L^{2}\left([0, T] \times \mathbb{T}^{2}\right)$ is a weak solution for the hyper-dissipative NSE \eqref{1.1}, if $u$ is divergence-free for all $t \in[0, T]$ and satisfies
$$
\int_{\mathbb{T}^{2}} u_{0} \mathscr{T}(0, \cdot) d x=-\int_{0}^{T} \int_{\mathbb{T}^{2}} u\left(\partial_{t} \mathscr{T}-(-\Delta)^{\alpha} \mathscr{T}+(u \cdot \nabla)\mathscr{T}\right) d x d t,
$$
where $\mathscr{T} \in C_{0}^{\infty}\left([0, T) \times \mathbb{T}^{2}\right)$ is the divergence-free test function.
\end{Definition}

In this paper, we focus on the following two super-critical regimes, whose borderlines include two endpoints of the gLPS condition \eqref{1.5}. More precisely, when $\alpha \in[1,\frac{3}{2})$ we consider the super-critical regime $\mathcal{A}_{1}$ given by
\begin{equation*}
\mathcal{A}_{1}:=\left\{(s, \gamma, p) \in[0,2) \times[1, \infty] \times[1, \infty]: 0 \leq s<\frac{4 \alpha-4}{\gamma}+\frac{2}{p}+1-2 \alpha\right\} \label{1.8}\tag{1.6}
\end{equation*}
and the super-critical regime $\mathcal{A}_{2}$ given by
\begin{equation*}
\mathcal{A}_{2}:=\left\{(s, \gamma, p) \in[0,2) \times[1, \infty] \times[1, \infty]: 0 \leq s<\frac{2 \alpha}{\gamma}+\frac{2 \alpha-2}{p}+1-2 \alpha\right\} . \label{1.9}\tag{1.7}
\end{equation*}
\vspace{-8mm}
\begin{figure}[h]
	\centering
	\includegraphics[width=0.67\linewidth]{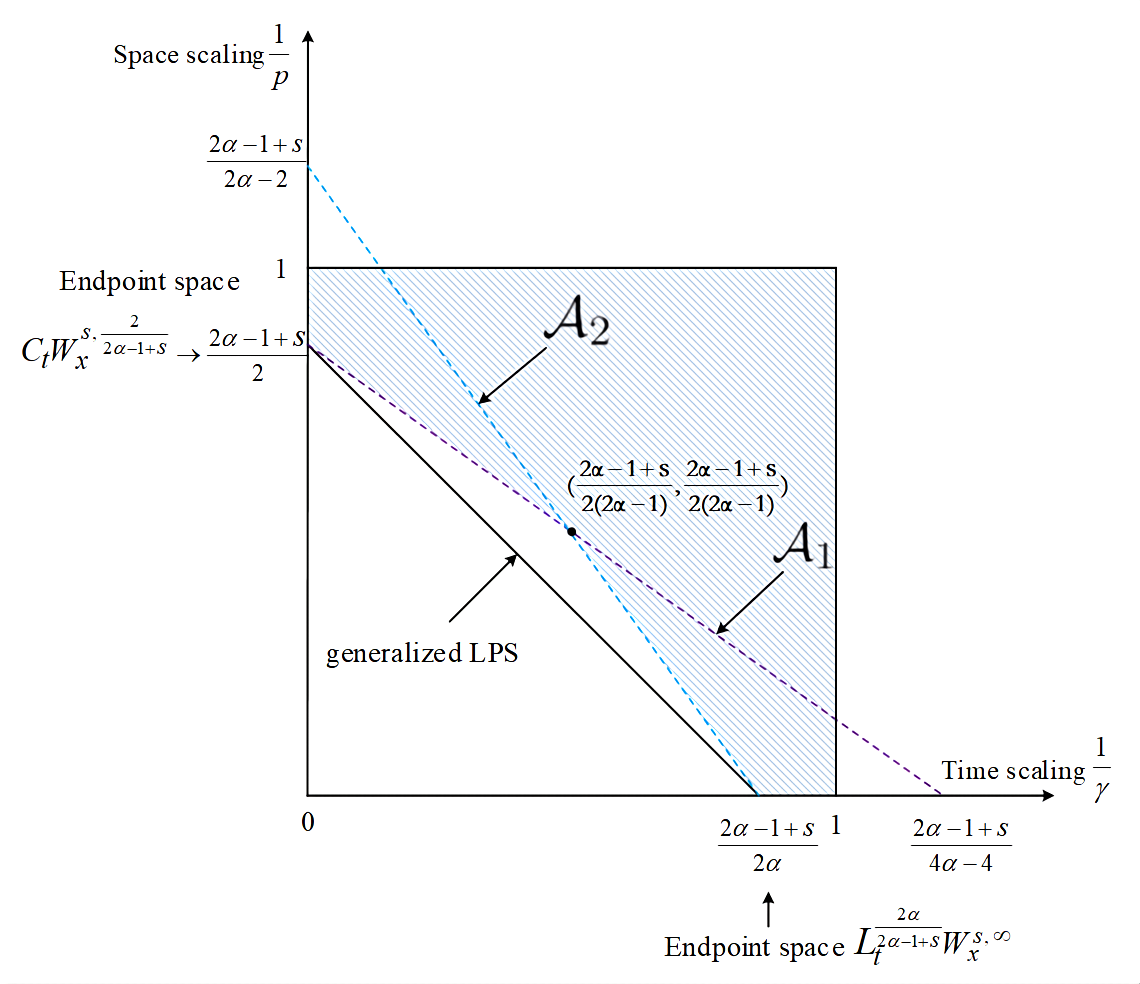}
	\caption{The case $\alpha \in\left[1, \frac{3}{2}\right), \ s\in[0,2).$}	
\label{fig:1}
\end{figure}

For a more intuitive expression, we present the super-critical regimes $\mathcal{A}_{1}$ and $\mathcal{A}_{2}$ in Figure \ref{fig:1} above, where the solution in the region below the generalized LPS condition \eqref{1.5} is shown to be unique in Theorem \ref{I.1} of Appendix \ref{G2}, and the shaded area represents the non-uniqueness results stated in Theorem \ref{A.2} below. In particular, the borderline of the super-critical regime $\mathcal{A}_{1}\cup\mathcal{A}_{2}$ contains two endpoint cases $(s,\infty, \frac{2}{2\alpha-1+s})$ and $(s, \frac{2\alpha}{2\alpha-1+s}, \infty)$, where $\alpha=1, s=0$ corresponds to the recent eye-catching works \cite{AC22,ACL22}. However, due to the $L_{t, x}^{2}$-criticality of the spatial-temporal convex integration method, it is still unknown whether the blank region between the generalized LPS condition and the borderline of $\mathcal{A}_{1}\cup\mathcal{A}_{2}$ is non-unique.

\begin{Theorem} \label{A.2}(Main theorem)
For $\alpha \in[1,\frac{3}{2})$, $(s, p, \gamma) \in \mathcal{A}_{1}\cup\mathcal{A}_{2}$. Given any smooth, divergence-free and mean-free vector field $\tilde{u}$ on $[0, T] \times \mathbb{T}^{2}$, the system \eqref{1.1} admits a velocity field $u$ and a set
$$
\mathscr{G}=\bigcup_{i=1}^{\infty}\left(a_{i}, b_{i}\right) \in[0, T]
$$
such that the following holds.

$(i)$ Weak solution: $u$ is a mean-free and divergence-free weak solution to the system \eqref{1.1} with the initial datum $\widetilde{u}(0)$.

$(ii)$ Regularity: $u \in H_{t, x}^{\beta^{\prime}} \cap L_{t}^{\gamma} W_{x}^{s, p}$ for some $\beta^{\prime}\in(0,1)$, and $u$ is a smooth solution on $\left(a_{i}, b_{i}\right)$ for every $i$, namely
$$
\left.u\right|_{\mathscr{G} \times \mathbb{T}^{2}} \in C^{\infty}\left(\mathscr{G} \times \mathbb{T}^{2}\right).
$$
In addition, $u$ is consistent with the unique smooth solution with the initial data
$\tilde{u}_{0}$ near $t=0$ and is also regular near $t= T$.

$(iii)$ Zero Hausdorff measure: the Hausdorff dimension of the singular set $\mathscr{B}=[0, T] / \mathscr{G}$ satisfies
$$
d_{\mathcal{H}}(\mathscr{B})<\eta_{*}, \quad \text{for some} \ \eta_{*}>0.
$$
Especially, the singular set $\mathscr{B}$ has a zero Hausdorff $\mathcal{H}^{\eta_{*}}$ measure, that is, $\mathcal{H}^{\eta_{*}}(\mathscr{B})=0$.

$(iv)$ $\varepsilon_{*}$-neighborhood of temporal support: there exists $\varepsilon_{*}>0$ such that
$$
\operatorname{supp}_{t} u \subseteq \mathscr{N}_{\varepsilon_{*}}\left(\operatorname{supp}_{t} \tilde{u}\right) .
$$

$(v)$ $\varepsilon_{*}$-close between the solution $u$ and the given vector field $\tilde{u}$ in $L_{t}^{1} L_{x}^{2}\cap L_{t}^{\gamma} W_{x}^{s, p}$:
$$
\|u-\tilde{u}\|_{L_{t}^{1} L_{x}^{2}\cap L_{t}^{\gamma} W_{x}^{s, p}} \leq \varepsilon_{*} .
$$
\end{Theorem}

The next result is our headline corollary, where we prove the strong non-uniqueness of weak solutions to the system \eqref{1.1} in the hyper-dissipative case where $\alpha \in[1,\frac{3}{2})$. This corollary is a direct consequence of our main theorem.  In particular, we settled the Question \ref{Q2} raised in this paper regarding the sharpness of gLPS condition at the two endpoint case.
\begin{Corollary}\label{A.3}
(Strong non-uniqueness) $\alpha \in[1,\frac{3}{2})$. For any weak solution $\widetilde{u}$ to the Cauchy problem \eqref{1.1}, there exists a different weak solution $u \in L_{t}^{\gamma} W_{x}^{s, p}$ to the system \eqref{1.1} with the same initial data, where $(s, \gamma, p) \in \mathcal{A}_{1} \cup \mathcal{A}_{2}$.

Furthermore, for any divergence-free initial data belonging to $L_{x}^{2}$, there exist infinitely many weak solutions in $L_{t}^{\gamma} W_{x}^{s, p}$ to the system \eqref{1.1} which are smooth outside a null set in time as described in Theorem \ref{A.2} $(iii)$.
\end{Corollary}

The sharp non-uniqueness results obtained in this article is also applicable to Lebesgue, Besov and Triebel-Lizorkin spaces, which has been verified in \cite{YL22}, and the following conclusion is also a direct consequence of Theorem \ref{A.2} in 2D case.
\begin{Corollary}\label{A.4}
For $\alpha \in[1,\frac{3}{2})$. There exist non-unique weak solutions to the Cauchy problem \eqref{1.1} in the super-critical spaces $C_{t} \mathbb{X}$, where $\mathbb{X}$ are the following three types of spaces:

$(i)$ $L^{p}$, for any $1 \leq p<\frac{2}{2 \alpha-1}$;

$(ii)$ $B_{p, q}^{s}$, for any $p,q\in[1,\infty]$, $-\infty<s<\frac{2}{p}+1-2 \alpha$;

$(iii)$ $F_{p, q}^{s}$, for any $p,q\in[1,\infty]$, $-\infty<s<\frac{2}{p}+1-2 \alpha$.
\end{Corollary}

Next, we list some remarks on our main results.
\begin{Remark}\label{A.5}(``2D intermittent jets'' $+$ ``temporal concentrated function'')
The fundamental spatial building blocks used in $\mathcal{A}_{1}$ can be replaced by 2D intermittent jets, i.e.
$$
\mathbf{W}_{k}:=-\varphi_{r_{\|}}\big(\lambda r_{\perp} N_{\Lambda}(k_{1} \cdot x +\mu t)\big) \psi_{r_{\perp}}^{\prime}\left(\lambda r_{\perp} N_{\Lambda} k \cdot x \right) k_{1}
$$
instead of accelerating jets \eqref{4.4} used in Section \ref{C}. Due to the disjoint nature of the temporal concentrated function $g_{(k)}$, the space-time ``building blocks'' $g_{(k)}\mathbf{W}_{k}$ do not intersect with each other, which is particularly evident in the proof of $\mathcal{A}_{2}$. Meanwhile, applying the time concentrated function $g_{(k)}(t)$ to the 2D intermittent jet selected in \cite{LTP20} seems to be a new spatial-temporal ``building block'' in 2D case, and the verification process has been carried out by the second author et al. in recent work \cite{LT24} , which once again indicates that sufficient intermittency has a significant impact on the viscosity exponent in dissipative dynamical systems.
\end{Remark}

After discussing the two-dimensional conclusion, it is naturally to consider what the result will be in the $d\geq2$ dimensional case. Therefore, we provide the following remark:
\begin{Remark}\label{A.6}
(Nonuniqueness of weak solutions for $d\geq2$)
We consider the following two super-critical regimes, whose borderlines contain two endpoints of the d-dimensional generalized Lady\v{z}enskaya-Prodi-Serrin condition \eqref{1.5}, i.e,
\begin{equation*}
\frac{2 \alpha}{\gamma}+\frac{d}{p}=2 \alpha-1+s.
\end{equation*}
More precisely, in the case $\alpha \in[\frac{d+2}{4},\frac{d+1}{2})$ we consider the super-critical regime $\mathcal{A}_{3}$ given by
\begin{equation*}
\mathcal{A}_{3}:=\left\{(s, \gamma, p) \in[0,d) \times[1, \infty] \times[1, \infty]: 0 \leq s<\frac{4 \alpha-(2+d)}{\gamma}+\frac{d}{p}+1-2 \alpha\right\} \label{1.71}\tag{1.8}
\end{equation*}
and in the case $\alpha \in[1,\frac{d+1}{2})$ we consider super-critical regime $\mathcal{A}_{4}$ given by
\begin{equation*}
\mathcal{A}_{4}:=\left\{(s, \gamma, p) \in[0,d) \times[1, \infty] \times[1, \infty]: 0 \leq s<\frac{2 \alpha}{\gamma}+\frac{2 \alpha-2}{p}+1-2 \alpha\right\} . \label{1.81}\tag{1.9}
\end{equation*}
\end{Remark}
\vspace{-3mm}
\begin{figure}[h]
	\centering
	\includegraphics[width=0.7\linewidth]{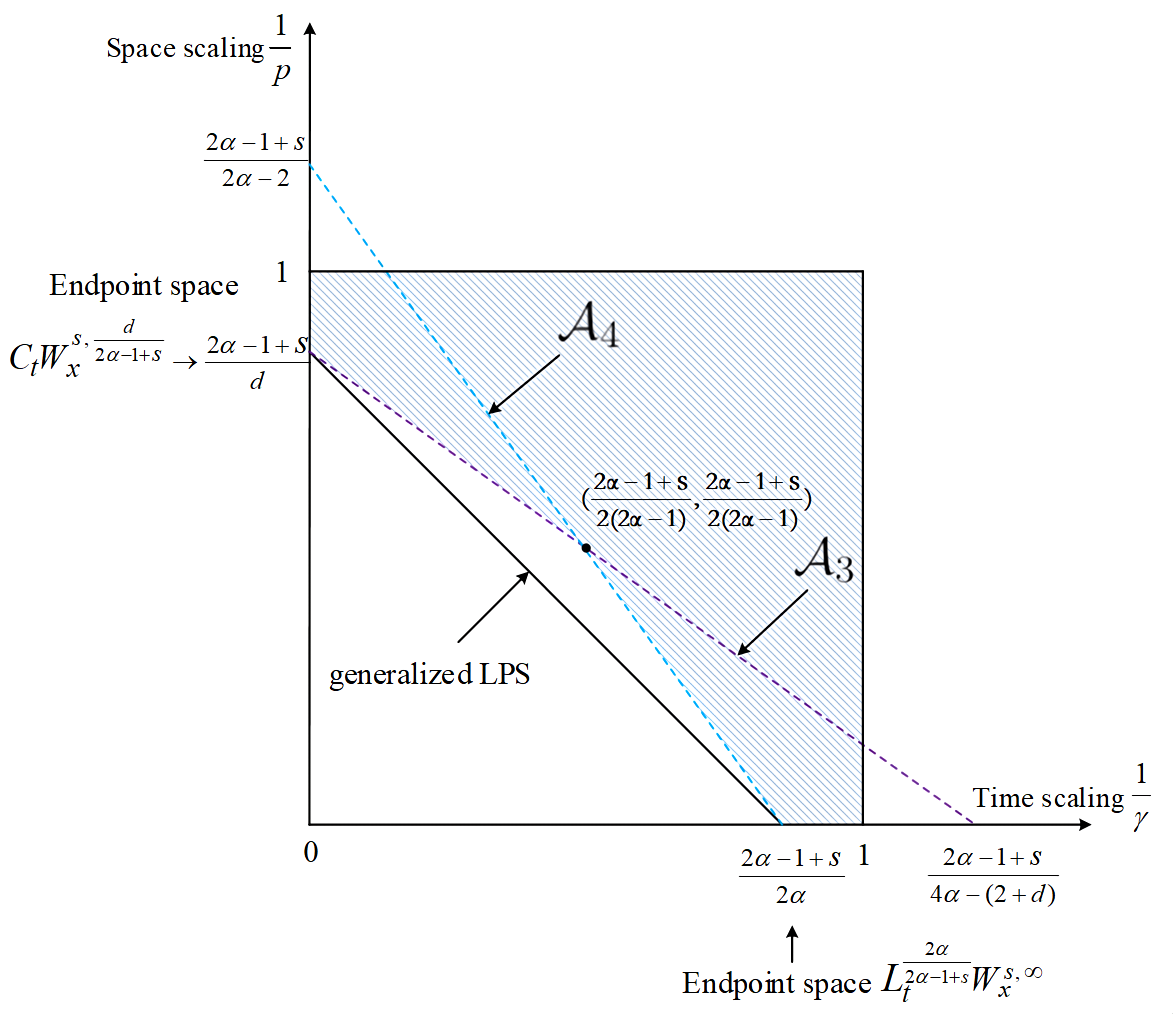}
	\caption{The case $\alpha \in[\frac{d+2}{4},\frac{d+1}{2}).$}
\label{fig:2}
\end{figure}

The super-critical regimes $\mathcal{A}_{3}$ and $\mathcal{A}_{4}$ can be seen in Figure \ref{fig:2} above. When the space dimension $d=2$,  $\mathcal{A}_{3}$ and  $\mathcal{A}_{4}$ correspond to  $\mathcal{A}_{1}$ and $\mathcal{A}_{2}$ discussed in Figure \ref{fig:1}, respectively. When $d=3$, the corresponding conclusions \eqref{1.71} and \eqref{1.81} have been obtained in \cite{YL22}. For the verification of the above Remark \ref{A.6}, it is only necessary to consider that the dimension $d>3$. Therefore, it is necessary to extend the 3D intermittent jets and 3D concentrated Mikado flow used in \cite{YL22} to the $d$-dimensional case. Additionally, in order to achieve the non-uniqueness near the two endpoint cases $(s,\infty, \frac{d}{2\alpha-1+s})$ and $(s, \frac{2\alpha}{2\alpha-1+s}, \infty)$ of $d$-dimensional gLPS condition, it is necessary to fully exploit temporal intermittency, that is, the appropriate temporal intermittency is approximately equivalent to $d$-dimensional spatial intermittency in the super-critical regimes $\mathcal{A}_{3}$ defined in \eqref{1.71}, and correspondingly, $(d+1)$-dimensional spatial intermittency in the super-critical regimes $\mathcal{A}_{4}$ defined in \eqref{1.81}. The detailed proof will be presented in our future research work \cite{DL24}.

\subsection{Comments on main results.}\ \ \ \
In this subsection, we will share some comments on the key observations drawn from our results.

{\bf (1) Strong non-uniqueness for the high viscosity $\alpha\geq1$.} In order to explore the strong non-uniqueness of weak solutions in the super-critical region $\mathcal{A}_{1}\cup\mathcal{A}_{2}$ when $\alpha\geq1$, it is necessary to fully exploit the temporal intermittency in an extremely way. To this end, we establish  the relationship between the intermittency of the temporal concentrated function $g_{(k)}$, the viscosity exponent $\alpha$ and the regularity of the weak solutions, see \eqref{4.1}, \eqref{6.1} and gLPS condition \eqref{1.5} for details. Our results recover the two-dimensional non-uniqueness results in \cite{AC22,ACL22} obtained by Cheskidov and Luo when $\alpha=1,\ s=0$ in Figure \ref{fig:1}, and are the first non-uniqueness results when $\alpha>1$ to the Cauchy problem \eqref{1.1}.

{\bf (2) Sharp non-uniqueness at two endpoints of gLPS condition.} In the remarkable papers \cite{AC22,ACL22}, Cheskidov and Luo firstly proved the sharp non-uniqueness at the two endpoint cases $(s, \gamma, p)=(0,2, \infty)$, $(s, \gamma, p)=(0, \infty, 2)$ for the 2D NSE.

In view of the uniqueness results in Theorem \ref{I.1} of Appendix \ref{G2}, the weak solutions at the two endpoints $(s,\gamma,p)=(s, \frac{2\alpha}{2\alpha-1+s}, \infty)$ and $(s,\gamma,p)=(s, \infty, \frac{2}{2\alpha-1+s})$ are unique, so Theorem \ref{A.2} provides the sharp non-uniqueness for the 2D hyper-dissipative NSE \eqref{1.1} at the two endpoints when $\alpha\geq1$ and $s\geq0$. Therefore, we extend the sharp non-uniqueness results in \cite{AC22,ACL22} to the hyper-dissipative case with $\alpha \in(1,\frac{3}{2})$.
In particularly, the viscosity exponent $\alpha=\frac{3}{2}$ seems the upper limit for the one endpoint case $(s,\infty, \frac{2}{2\alpha-1+s})$ due to $\frac{2}{2\alpha-1+s}>1$ with $s\in[0,2)$.

Combining the results in \cite{YL22} shows that the endpoint on the space coordinate axis in Figure \ref{fig:1} is related to the spatial dimension, and it remains unknown whether there exists non-unique solutions for 3D NSE in proximity to the endpoint $C_{t}L^{3}_{x}$. At present, due to the $L_{t, x}^{2}$-criticality of space-time convex integration method, the non-uniqueness of the remaining super-critical region in Figure \ref{fig:1} is unknown,
we are looking forward to better technical means for processing.

{\bf(3) Sharp non-uniqueness in the super-critical Lebesgue and Besov spaces.}
As far as we know, $(ii)$ and $(iii)$ of Corollary \ref{A.4} is currently the first non-unique result in Besov space, and in Triebel-Lizorkin space in 2D. As for $(i)$ of Corollary \ref{A.4}, we have extended the non-uniqueness results on 2D NSE $(\alpha=1)$ in \cite{ACL22} to the hyper-dissipative case $\alpha>1$.

Due to the well-posedness result in critical space $L_{x}^{2 /(2 \alpha-1)}$ by \cite{ST98} when $\alpha\geq1$, we proved that the non-uniqueness in Corollary \ref{A.4} $(i)$ is sharp in the super-critical Lebesgue spaces to 2D hyper-dissipative NSE \eqref{1.1}.

For Corollary \ref{A.4} $(ii)$, it has been shown that in \cite{ST03}, for $\alpha\geq1,$  \eqref{1.1} possess a unique global solution in $\dot{B}_{p, q}^{\frac{2}{p}+1-2 \alpha}$, which implies that the non-uniqueness in Corollary \ref{A.4} $(ii)$ is sharp in the super-critical Besov spaces. Due to the similarity with the definition of Besov space, verifying that the non-uniqueness in Corollary \ref{A.4} $(iii)$ is sharp in super-critical Triebel-Lizorkin space is trivial, and we omit the details here.

{\bf(4) Nonuniqueness of weak solutions for the hyper-dissipative NSE equations in all dimensions $d\geq2$.} In Remark \ref{A.6}, we expect to employ the spatial-temporal convex integration method to obtain the non-uniqueness of weak solutions to the hyper-dissipative NSE in all dimensions $d\geq2$ in the super-critical regime $\mathcal{A}_{3}$ and $\mathcal{A}_{4}$, which is entirely based on the intermittent development of the temporal concentrated function $g_{(k)}(t)$. The research results are highly generalizable and include most of the non-uniqueness results of weak solutions of the hyper-dissipative NSE equations for $\alpha\geq1$, such as \cite{AC22,ACL22,YL22}. Specifically, when $\alpha\geq\frac{d}{4}+\frac{1}{2}$, the known result in \cite{JL69} tells us that $C_{t}L^{2}$ weak solution is unique, however surprisingly, Remark \ref{A.6} implies that even the viscosity exponent $\alpha\geq\frac{d}{4}+\frac{1}{2}$, there also exist non-unique solutions in super-critical regimes $\mathcal{A}_{3}$ and $\mathcal{A}_{4}$ as shown in Figure \ref{fig:2}, where the dimensions $d=2,3$ have been discussed in this article and \cite{YL22}, respectively. Meanwhile, it is worthy noting that $L_{t}^{2}C_{x}$ is another endpoint of the LPS condition, regardless of the space dimension, Remark \ref{A.6} also implies that for dimensions $d\geq2$, there exist non-unique weak solutions in super-critical space when $\alpha\geq1$, where $\alpha=1$ has been discussed in \cite{AC22}. Moreover, we would expect to establish the relationship between the space dimension, temporal intermittency, viscosity exponent $\alpha$, and the regularity of weak solutions. Due to the difference between the two-dimensional spatial building blocks setting and high dimensions $d\geq3$, we first present the corresponding results in two dimensions in this article.

The remainder of this paper is organized as follows. we give the main iteration Proposition \ref{B.3} and concentrate the Reynolds error in Section \ref{B}. Sections \ref{C} and \ref{D} are primarily concerned with the endpoint case of $(s,\infty, \frac{2}{2\alpha-1+s})$. Specifically, Section \ref{C} involves the construction of velocity perturbations and the preparation of crucial algebraic identities and analytical estimates. Then, in Section \ref{D}, we address the Reynolds stress. The another endpoint case  $(s, \frac{2\alpha}{2\alpha-1+s}, \infty)$ is primarily addressed in Section \ref{E}. Finally, the proofs of the main results are presented in Section \ref{F}, while Appendix \ref{G1}, \ref{G2} and \ref{G3} contain some preliminary results used in the proof.

\begin{center}
\section{Main iteration and concentrating Reynolds error}\label{B}
\end{center}\ \ \ \
In this section, we first present the main iteration of the velocity and Reynolds stress, which is the heart of the proof for our main theorem.\\
\subsection{Main iteration}\ \ \ \ We consider the approximate solutions to the following 2D generalized Navier-Stokes-Reynolds system for each  integer $q \geq 0$,
\[
\left\{\begin{array}{l}\label{2.1}
\partial_{t} u_{q}+(-\Delta)^{\alpha} u_{q}+\operatorname{div}\left(u_{q} \otimes u_{q}\right)+\nabla\mathsf{P}_{q}=\operatorname{div} \mathring{\mathsf{R}}_{q},  \tag{2.1}\\
\operatorname{div} u_{q}=0,
\end{array}\right.
\]
where $\mathring{\mathsf{R}}_{q}: \mathbb{T}^{2}\times[0,T]\rightarrow \mathcal{S}_{0}^{2\times2}$ is a $2 \times 2$ symmetric traceless  matrix, known as Reynolds stress in the literature.

In addition the pressure term can be recoverd from the elliptic equation by taking the divergence of equations \eqref{B.1}. Specifically, the equation is given by
$$
\Delta \mathsf{P}_{q}=\operatorname{div} \operatorname{div}(\mathring{\mathsf{R}}_{q}-u_{q} \otimes u_{q}),
$$
which, in conjunction with the standard zero spatial mean condition $\fint_{\mathbb{T}^2} \mathsf{P} dx= 0$, uniquely determines the pressure.

To harness the delicate temporal singular set of approximate solutions, we embrace the concept of well-prepared solutions introduced in \cite{AC22}.

\begin{Definition}\label{B.1}(Well-preparedness) For any fixed $\eta^{*}\in(0,1)$, and $\eta\in\left(0, \eta_{*}\right)$. The smooth solution $\left(u_{q}, \mathring{\mathsf{R}}_{q}\right)$ to \eqref{2.1} on the interval $[0, T]$ is deemed to be well-prepared if there exist a set $\mathcal{I}$ and a length scale $\vartheta>0$, such that $\mathcal{I}$ is the union of at most $\vartheta^{-\eta}$ many closed intervals of length scale $5 \vartheta$ and
$$
\mathring{\mathsf{R}}_{q}(t, x)=0, \quad \text { if } \quad \operatorname{dist}\left(t, \mathcal{I}^{c}\right) \leq \vartheta.
$$
\end{Definition}

To measure the size of the relaxation solutions $\left(u_{q}, \mathring{\mathsf{R}}_{q}\right)$ for $q \in \mathbb{N}$, two crucial parameters are employed: the frequency parameter $\lambda_{q}$ and the amplitude parameter $\delta_{q+2}$:
\begin{equation*}
\lambda_{q}=a^{\left(b^{q}\right)}, \quad \delta_{q+2}=\lambda_{q+2}^{-2 \beta}. \label{2.2}\tag{2.2}
\end{equation*}
where $a\in\mathbb{N}$ is a large integer, the parameter $\beta>0$ represents regularity and $b\in 2 \mathbb{N}$ satisfies
\begin{equation*}
b>\frac{1000}{\varepsilon \eta_{*}}, \quad 0<\beta<\frac{1}{100 b^{2}}, \label{2.3}\tag{2.3}
\end{equation*}
where $\varepsilon \in \mathbb{Q}_{+}$ is sufficiently small such that, for any $(s, p, \gamma) \in \mathcal{A}_{1}$,
\begin{equation*}
\varepsilon \leq \frac{1}{20} \min \left\{3-2\alpha, \frac{4 \alpha-4}{\gamma}+\frac{2}{p}-(2 \alpha-1)-s\right\} \quad \text { and } \quad b \varepsilon \in \mathbb{N} , \label{2.4}\tag{2.4}
\end{equation*}
and for $(s, p, \gamma) \in \mathcal{A}_{2}$,
\begin{equation*}
\varepsilon \leq \frac{1}{20} \min \left\{3-2\alpha, \frac{2 \alpha}{\gamma}+\frac{2 \alpha-2}{p}-(2 \alpha-1)-s\right\} \quad \text { and } \quad b(3-2\alpha-8 \varepsilon) \in \mathbb{N} \text {. } \label{2.5}\tag{2.5}
\end{equation*}

The purpose of the iteration is to prove that when $q$ approaches infinity, the Reynolds stress disappears in an appropriate space and the limit value of $u_{q}$ is the solution of the initial equations \eqref{1.1}.  This approach is formally outlined in the subsequent iterative steps:
\begin{align*}
& \left\|u_{q}\right\|_{L_{t}^{\infty} H_{x}^{3}} \lesssim \lambda_{q}^{5},  \label{2.6}\tag{2.6}\\
& \left\|\partial_{t} u_{q}\right\|_{L_{t}^{\infty} H_{x}^{2}} \lesssim \lambda_{q}^{7}, \label{2.7}\tag{2.7}\\
& \left\|\mathring{\mathsf{R}}_{q}\right\|_{L_{t}^{\infty} H_{x}^{3}} \lesssim \lambda_{q}^{9},  \label{2.8}\tag{2.8}\\
& \left\|\mathring{\mathsf{R}}_{q}\right\|_{L_{t}^{\infty} H_{x}^{4}} \lesssim \lambda_{q}^{10},  \label{2.9}\tag{2.9}
\end{align*}
and
\begin{equation*}
\left\|\mathring{\mathsf{R}}_{q}\right\|_{L_{t, x}^{1}} \leq \lambda_{q}^{-2\varepsilon_{R}} \delta_{q+1},\label{2.10} \tag{2.10}
\end{equation*}
where the implicit constants $C$ are independent of $q$ and the small parameter $\varepsilon_{R}$ satisfies that
$$
0<\varepsilon_{R}<\frac{\varepsilon}{20}.
$$

Our main theorem relies on the following
crucial iteration results formulated below.
\begin{Proposition}\label{B.3}(Main iteration). Let $(s, p, \gamma) \in \mathcal{A}_{1}\cup \mathcal{A}_{2}$ for $\alpha \in[1,\frac{3}{2})$. Subsequently, there exist constants $\beta \in(0,1), M^{*}>0$ sufficiently large and $a_{0}=a_{0}\left(\beta, M^{*}\right)$, such that for any integer $a \geq a_{0}$, the following holds:

Suppose that $\left(u_{q}, \mathring{\mathsf{R}}_{q}\right)$ is a well-prepared solution to \eqref{2.1} within the set $\mathcal{I}_{q}$ and at the length scale $\vartheta_{q}$, satisfying conditions \eqref{2.6}-\eqref{2.10}. Then, there exists another well-prepared solution $\left(u_{q+1}, \mathring{\mathsf{R}}_{q+1}\right)$ to \eqref{2.1} fulfills conditions \eqref{2.6}-\eqref{2.10} with $q+1$  replacing $q$ within subset $\mathcal{I}_{q+1} \subseteq \mathcal{I}_{q}$, $0,T\notin \mathcal{I}_{q+1}$, and the length scale $\vartheta_{q+1}<\vartheta_{q} / 2$. Furthermore, we have
\begin{align*}
& \left\|u_{q+1}-u_{q}\right\|_{L_{t, x}^{2}} \leq M^{*} \delta_{q+1}^{\frac{1}{2}},  \label{2.11}\tag{2.11}\\
& \left\|u_{q+1}-u_{q}\right\|_{L_{t}^{1} L_{x}^{2}} \leq \delta_{q+2}^{\frac{1}{2}},  \label{2.12}\tag{2.12}\\
& \left\|u_{q+1}-u_{q}\right\|_{L_{t}^{\gamma} W_{x}^{s, p}} \leq \delta_{q+2}^{\frac{1}{2}}, \label{2.13}\tag{2.13}
\end{align*}
and
\begin{equation*}
\operatorname{supp}_{t}\left(u_{q+1}, \mathring{\mathsf{R}}_{q+1}\right) \subseteq N_{\delta_{q+2}^{\frac{1}{2}}}\left(\operatorname{supp}_{t}\left(u_{q}, \mathring{\mathsf{R}}_{q}\right)\right). \label{2.14}\tag{2.14}
\end{equation*}
\end{Proposition}

The proof of the main iteration theorem is firmly based on the concentrating the Reynolds error and the innovative approach of space-time intermittent convex integration. On the one hand, the procedure of concentrating the Reynolds stress error by gluing technique, which was first developed in \cite{TBU15,TB16,TB19,PI18} to solve the famous Onsager conjecture for 3D Euler equations. The gluing technique to obtain partial regularity in time of weak solutions for 3D NSE was first introduced in \cite{BCV22}, where the Hausdorff dimension of the constructed solution is strictly less than one.  Subsequently, Cheskidov and Luo \cite{AC22} employed this technique to establish the sharp non-uniqueness of weak solutions at the endpoint case $(s, \gamma, p)=(0,2, \infty)$, where the singular sets exhibited a small Hausdorff dimension over time. Furthermore, the gluing technique has also been utilized in \cite{CDLDR18,DR19} to demonstrate the non-uniqueness of Leray-Hopf solutions for the 3D hypo-dissipative NSE. For more results on the partial regularity of weak solutions in time, please refer to the corresponding summary description in \cite{YL22}.

In this article, we will also use the gluing technique to concentrate the Reynolds stress error. The brief description is as follows, and the detailed process can be found in \cite{YL22,AC22,BCV22}.

\subsection{Concentrating the Reynolds error}\ \ \ \
In order to focus the singular times into a null set characterized by a low Hausdorff dimension, we  introduce a properly prepared solution $\left(u_{q}, \mathring{\mathsf{R}}_{q}\right)$ to equations \eqref{2.1} at level $q$, and divide the time interval $[0, T]$ into $\ell_{q+1}$ equal-length sub-intervals $\left[t_{i}, t_{i+1}\right]$ with the length of $T / \ell_{q+1}$, where $t_{i}=i T / \ell_{q+1}$ and $\ell_{q+1}$ may depend on the given solution $\left(u_{q}, \mathring{\mathsf{R}}_{q}\right)$. Subsequently, we solve the following hyper-dissipative NSE on each of small intervals $\left[t_{i}, t_{i+1}+\vartheta_{q+1}\right]$:
\[
\left\{\begin{array}{l}
\partial_{t} v_{i}+(-\Delta)^{\alpha} v_{i}+\left(v_{i} \cdot \nabla\right) v_{i}+\nabla\mathsf{P}_{i}=0, \label{2.15} \tag{2.15}\\
\operatorname{div} v_{i}=0, \\
\left.v_{i}\right|_{t=t_{i}}=u_{q}\left(t_{i}\right),
\end{array}\right.
\]
where $\ell_{q+1}$ and $\vartheta_{q+1}$, is detailed as
\begin{equation*}
T / \ell_{q+1}=\lambda_{q}^{-12}, \quad \vartheta_{q+1}:=\left(T / \ell_{q+1}\right)^{1 / \eta} \simeq \lambda_{q}^{-\frac{12}{\eta}}, \label{3.1}\tag{2.16}
\end{equation*}
and $\eta$ is a small constant such that
$
0<\frac{\eta_{*}}{2}<\eta<\eta_{*}<1
$ with $\eta_{*}$ as stated in Theorem \ref{A.2}.

According to the definition of $\varepsilon$ defined in \eqref{2.4}-\eqref{2.5}, it is easy to deduce \begin{equation*}
\vartheta_{q+1}^{-20} \simeq \ell_{q+1}^{\frac{20}{\eta}} \ll \lambda_{q+1}^{\frac{\varepsilon}{2}} \label{2.16}\tag{2.17}
\end{equation*}
which suggests that the amplitudes of velocity perturbations in Lemmas \ref{D.5} and \ref{F.2} oscillate at a significantly lower frequency compared to those of temporal and spatial building blocks.

Moreover, in accordance with the established local well-posedness theory, Proposition \ref{C.2} can be applied with the substitution of $t_{*}$ and $t_{0}$ by $t_{i+1}+\vartheta_{q+1}$ and $t_{i}$ respectively, there exists a unique local smooth solution $v_{i}$ to equations \eqref{2.15} when $\vartheta_{q+1}$ is sufficiently small . Subsequently,  we meticulously glue together these local solutions in a manner that ensures the glued solution $\sum_{i} \mathscr{X} _{i} v_{i}$ accurately solves \eqref{1.1} for the majority of the time interval $[0, T]$. It is worthy mentioning that the new Reynolds stress error would exhibit smaller and disconnected temporal supports. More precisely, we introduce a smooth partition of unity on $[0, T]$ denoted by $\left\{\mathscr{X} _{i}\right\}_{i=0}^{\ell_{q+1}-1}$ such that $0 \leq \mathscr{X} _{i}(t) \leq 1, $
for $0<i<\ell_{q+1}-1$,
\[
\mathscr{X} _{i}= \begin{cases}1 & \text { if } t_{i}+\vartheta_{q+1} \leq t \leq t_{i+1},  \label{3.18}\tag{2.18}\\ 0 & \text { if } t \leq t_{i}, \text { or } t \geq t_{i+1}+\vartheta_{q+1},\end{cases}
\]
and for $i=0$,
\[
\mathscr{X} _{i}= \begin{cases}1 & \text { if } 0 \leq t \leq t_{i+1},  \label{3.19}\tag{2.19}\\ 0 & \text { if } t \geq t_{i+1}+\vartheta_{q+1},\end{cases}
\]
and for $i=\ell_{q+1}-1$,
\[
\mathscr{X} _{i}= \begin{cases}1 & \text { if } t_{i}+\vartheta_{q+1} \leq t \leq T,  \label{3.20}\tag{2.20}\\ 0 & \text { if } t \leq t_{i}.\end{cases}
\]
We assume the sequence $\mathscr{X} _{i}$ satisfy $
\left\|\partial_{t}^{M} \mathscr{X} _{i}\right\|_{L_{t}^{\infty}} \lesssim \vartheta_{q+1}^{-M}$, where the implicit constant is independent of both $\vartheta_{q+1}$ and the indices $i, M \geq 0$.

Then we define the gluing solution $\widetilde{u}_{q}$ as follows:
\begin{equation*}
\widetilde{u}_{q}:=\sum_{i=0}^{\ell_{q+1}-1} \mathscr{X} _{i} v_{i}, \label{3.21}\tag{2.21}
\end{equation*}
where $\widetilde{u}_{q}:[0, T] \times \mathbb{T}^{2} \rightarrow \mathbb{R}^{2}$ is divergence-free and mean-free. Additionally, for  $t \in\left[t_{i}, t_{i+1}\right]$, we have
$\tilde{u}_{q}=\left(1-\mathscr{X} _{i}\right) v_{i-1}+\mathscr{X} _{i} v_{i}$
and satisfies the following equation:
\begin{equation*}
\partial_{t} \widetilde{u}_{q}+(-\Delta)^{\alpha} \widetilde{u}_{q}+\operatorname{div}\left(\widetilde{u}_{q} \otimes \widetilde{u}_{q}\right)+\nabla \widetilde{p}=\operatorname{div} \mathring{\widetilde{\mathsf{R}}}_q, \label{3.22}\tag{2.22}
\end{equation*}
where the new Reynolds stress tensor given by
\begin{equation*}
\mathring{\widetilde{\mathsf{R}}}_q=\partial_{t} \mathscr{X} _{i} \mathcal{R}\left(v_{i}-v_{i-1}\right)-\mathscr{X} _{i}\left(1-\mathscr{X} _{i}\right)\left(\left(v_{i}-v_{i-1}\right) \mathring{\otimes}\left(v_{i}-v_{i-1}\right)\right), \label{3.23}\tag{2.23}
\end{equation*}
for $0\leq i\leq \ell_{q+1}-1$, with the convention that $v_{-1}\equiv0$ and the pressure $\widetilde{p}:[0,1] \times \mathbb{T}^{2}\rightarrow\mathbb{R}$ is defined by:
\begin{equation*}
\widetilde{p}=\mathscr{X} _{i} \mathsf{P}_{i}+(1-\mathscr{X} _{i}) \mathsf{P}_{i-1}-\mathscr{X} _{i}\left(1-\mathscr{X} _{i}\right)\left(\left|v_{i}-v_{i-1}\right|^{2}-\int_{\mathbb{T}^{2}}\left|v_{i}-v_{i-1}\right|^{2} \mathrm{~d} x\right), \label{3.24}\tag{2.24}
\end{equation*}
where the inverse-divergence operator $\mathcal{R}: C^{\infty}\left(\mathbb{T}^2, \mathbb{R}^2\right) \rightarrow C^{\infty}\left(\mathbb{T}^2, \mathcal{S}_0^{2 \times 2}\right)$ is defined as in \cite{CD13}, specifically
$$ (\mathcal{R} v)_{i j}=-\Delta^{-1} \partial_k \delta_{i j}v_k+\Delta^{-1} \partial_i \delta_{j k}v_k+\Delta^{-1} \partial_j \delta_{i k} v_k.
$$
By a direct computation, one can also show that
$\dive(\mathcal{R}v)=v$, where $v$ is a  mean-free smooth function satisfying $\int_{\mathbb{T}^{2}} v dx=0$. The antidivergence operator $\mathcal{R}$ is bounded on $L^p\left(\mathbb{T}^2\right)$ for any $1 < p < \infty$ (see also \cite{MS18}).

According to the definition of the set $\left\{\mathscr{X} _{i}\right\}_{i}$, the new Reynolds stress $\mathring{\widetilde{\mathsf{R}}}_{q}$  is localized within a $\vartheta_{q+1}$-neighborhood of each $t_{i}$ for all indices $i$ satisfying $0 \leq i<\ell_{q+1}$. Consequently, this implies that the singular set of $\widetilde{u}_{q}$ can be entirely covered by $\ell_{q+1}$ small intervals of length scale of $5\vartheta_{q+1}$, which leads to a small Hausdorff dimension.
Furthermore, as stated in reference \cite{BCV22}, given that $u_{q}$ is already a smooth solution to the system \eqref{1.1} on the majority of the interval $[0, T]$, specifically within the $\vartheta_{q}$-neighborhood of the complement of a small set $\mathcal{I}_{q}$, if both $t_{i-1}$ and $t_{i}$ reside in this region, then $\widetilde{u}_{q}=v_{i-1}=v_{i}=u_{q}$ holds true on the support of $\mathscr{X} _{i} \mathscr{X} _{i-1}$. Consequently, we introduce an index set $\mathscr{C}$
\begin{equation*}
\mathscr{C}:=\left\{i \in \mathbb{Z}: 1 \leq i \leq \ell_{q+1}-1 \text { and } \mathring{\mathsf{R}}_{q} \not \equiv 0 \text { on }\left[t_{i-1}, t_{i}+\vartheta_{q+1}\right] \cap[0, T]\right\} \label{2.18}\tag{2.25}
\end{equation*}
to identify regions where $\widetilde{u}_{q}$ may not be an exact solution to the system \eqref{1.1},
and choose $\mathcal{I}_{q+1}$ in the way
$$
\mathcal{I}_{q+1}:=\bigcup_{i \in \mathscr{C}}\left[t_i-2 \vartheta_{q+1}, t_i+3 \vartheta_{q+1}\right].
$$
Based on the setting of \eqref{2.2} and \eqref{3.1}, it is easy to verify that the bad set $\mathcal{I}_{q+1}$ associated with $\left(\widetilde{u}_{q}, \mathring{\widetilde{\mathsf{R}}}_q\right)$ is contained within $\mathcal{I}_{q}$. Then we have the following conclusion:
\begin{Proposition}\label{C.1} For $\alpha \in[1,\frac{3}{2})$, $(s, p, \gamma) \in \mathcal{A}_{1}\cup\mathcal{A}_{2}$. Given a smoothly prepared solution $\left(u_{q}, \mathring{\mathsf{R}}_{q}\right)$ to \eqref{2.1} on a set $\mathcal{I}_{q}$ and a length scale $\vartheta_{q}$. Subsequently, there exists another smoothly prepared solution $\left(\widetilde{u}_q, \mathring{\widetilde{\mathsf{R}}}_q\right)$ to the same equations \eqref{2.1} on a subset $\mathcal{I}_{q+1} \subseteq I_q$ and a smaller length scale $\vartheta_{q+1}\left(<\vartheta_q / 2\right)$ such that
\begin{align*}
& \mathring{\widetilde{\mathsf{R}}}_q(t, x)=0 \quad \text { if } \ \operatorname{dist}\left(t, \mathcal{I}_{q+1}^{c}\right) \leq 3 \vartheta_{q+1} / 2,  \label{3.3}\tag{2.26}\\
& \left\|\widetilde{u}_{q}\right\|_{L_{t}^{\infty} H_{x}^{3}} \lesssim \lambda_{q}^{5},  \label{3.4}\tag{2.27}\\
& \left\|\widetilde{u}_{q}-u_{q}\right\|_{L_{t}^{\infty} L_{x}^{2}} \lesssim \lambda_{q}^{-3},  \label{3.5}\tag{2.28}\\
& \left\|\mathring{\widetilde{\mathsf{R}}}_q\right\|_{L_{t, x}^{1}} \leq \lambda_{q}^{-\frac{2}{3}\varepsilon_{R}} \delta_{q+1},  \label{3.6}\tag{2.29}
\end{align*}
and
\begin{equation*}
\left\|\partial_{t}^{M} \nabla^{N} \mathring{\widetilde{\mathsf{R}}}_q\right\|_{L_{t}^{\infty} H_{x}^{3}} \lesssim \vartheta_{q+1}^{-M-1} \ell_{q+1}^{-\frac{N}{2 \alpha}} \lambda_{q}^{5} \lesssim\vartheta_{q+1}^{-M-N-1} \lambda_{q}^{5}, \label{3.7}\tag{2.30}
\end{equation*}
where the implicit constants are independent of the variable $q$, and  $0, T \notin \mathcal{I}_{q+1}$.
\end{Proposition}
{\bf Proof.} The proof of Proposition \ref{C.1} can be directly obtained from Proposition \ref{C.2} and Lemma \ref{C.3} in the Appendix \ref{G3}. Specifically, the concentration process of time singular sets is independent of the spatial dimension, for the proof of its related details, please refer to \cite{BCV22,YL22}.

It is worthy noting that the newly introduced Reynolds stress exhibits a minor penalty of $\lambda_{q}^{-4 \varepsilon_{R} / 3}$ in terms of decay rate as demonstrated in \eqref{3.6}, which is deemed tolerable for the subsequent convex integration stage. For further reference, please consult \cite{AC22}.

{\centering
\section{Velocity perturbations in the super-critical regime $\mathcal{A}_{1}$}\label{C}}

In this section,  we mainly focus on constructing the velocity perturbations and selecting appropriate intermittent or oscillatory parameters, to ensure the corresponding inductive estimates in Proposition \ref{B.3} hold in convex integration schemes.

In order to achieve the sharp non-uniqueness result at the one endpoint case $(s, \gamma, p)=(s,\infty, \frac{2}{2\alpha-1+s})$ with $\alpha \in[1,\frac{3}{2})$, we select the acceleration jets $\mathbf{W}_{k}$, introduced in \cite{ACL22}, as our primary spatial building blocks (see \eqref{4.5} below), a key characteristic of which is its nearly 2D intermittency, i.e.,
\begin{equation*}
\left\|\mathbf{W}_{k}\right\|_{L_{t}^{\infty} L_{x}^{1}} \lesssim \lambda^{-1+}. \label{2.19}
\end{equation*}
Additionally, the appropriate temporal intermittency closely corresponds to $(4 \alpha-4)$-dimensional spatial intermittency. Specially, when the viscosity exponent $\alpha$ approaches $\frac{3}{2}$, the temporal intermittency nearly attains 2D spatial intermittency. Furthermore, we provide a precise selection of six admissible parameters as follows:
\begin{equation*}
r_{\perp}:=\lambda_{q+1}^{-1+2 \varepsilon}, \ r_{\|}:=\lambda_{q+1}^{-1+10 \varepsilon},\ \lambda:=\lambda_{q+1},\ \mu:=\lambda_{q+1}^{1 +2 \varepsilon},\ \tau:=\lambda_{q+1}^{4 \alpha-4+16 \varepsilon}, \ \sigma:=\lambda_{q+1}^{2 \varepsilon} \label{4.1}\tag{3.1}
\end{equation*}
where $\varepsilon$ denotes a sufficiently small constant that fulfills the condition stated in \eqref{2.4}. \\
\subsection{Temporal building blocks.}\ \ \ \
Firstly, we introduce the temporal building blocks that possess two crucial parameters, $\tau$ and $\sigma$, which govern concentration and oscillation over time, respectively. It is crucial to consider suitable temporal shifts $t_{k}$ within these building blocks to ensure that the supports of different temporal building blocks are disjoint. This construction is motivated by recent advancements in temporal intermittency, as detailed in references \cite{AC21,AC22,ACL22,ACX22}, which permits to obtain sharp non-uniqueness results.

More precisely, let $\Lambda$ be a subset of $\mathbb{S} \cap \mathbb{Q}^{2}$ representing the wavevector set, as specified in Geometric Lemma \ref{G.1}. Define $\left\{g_k\right\}_{k \in \Lambda} \subset C_c^{\infty}([0, T])$ as cut-off functions with mean zero,  such that for $k \neq k^{\prime}$, the temporal supports of $g_k$ and $g_{k^{\prime}}$ are disjoint.  Furthermore, for all $k \in \Lambda$, we have
\begin{align*}
\fint_0^T g_k^2(t) \mathrm{d} t=1.
\label{4.13}\tag{3.2}
\end{align*}
Since there are finitely many wavevectors in $\Lambda$, the existence of such $\left\{g_k\right\}_{k \in \Lambda}$ can be guaranteed by constructing them as, for example, $g_k=g\left(t-t_k\right)$, where $g \in C_c^{\infty}([0, T])$ has a very small support, and $\left\{t_k\right\}_{k \in \Lambda}$ are temporal shifts chosen such that the supports of $\left\{g_k\right\}_{k \in \Lambda}$ are disjoint. Subsequently, for each $k \in \Lambda$, we rescale $g_k$ by
$g_{k, \tau}(t)=\tau^{\frac{1}{2}} g_k(\tau t),$ where the concentration parameter $\tau$ is defined by \eqref{4.1}. Note that the size of those time intervals $\sim \tau^{-1}$ then determines the level of temporal concentration. By an abuse of notation, we consider $g_{k, \tau}$ as a periodic function on $[0, T]$ by extension.
Additionally, we define
\begin{align*}
h_{k, \tau}(t):=\int_0^t\left(g_{k, \tau}^2(s)-1\right) d s, \quad t \in[0, T],
\label{4.14}\tag{3.3}
\end{align*}
and
\begin{align*}
g_{(k)}(t):=g_{k, \tau}(\sigma t), \quad h_{(k)}(t):=h_{k, \tau}(\sigma t) .
\label{4.15}\tag{3.4}
\end{align*}

Subsequently, we obtain
\begin{align*}
\partial_t\left(\sigma^{-1} h_{(k)}\right)=g_{(k)}^2-1=g_{(k)}^2-\fint_0^T g_{(k)}^2(t) \mathrm{d} t,
\label{4.16}\tag{3.5}
\end{align*}
where $\sigma$ is given by \eqref{4.1}.

We draw upon the crucial estimates of $g_{(k)}$ and $h_{(k)}$ established in previous works \cite{AC21,YLZ22} and present them in the following Lemma.

\begin{Lemma}\label{D.3} (Estimates of temporal intermittency). For $\gamma \in[1, \infty], M \in \mathbb{N}$, we have
\begin{equation*}
\left\|\partial_{t}^{M} g_{(k)}\right\|_{L_{t}^{\gamma}} \lesssim \sigma^{M} \tau^{M+\frac{1}{2}-\frac{1}{\gamma}}, \label{4.17}\tag{3.6}
\end{equation*}
where the implicit constants do not rely on both $\sigma$ and $\tau$. Moreover, we have the following estimate for $h_{(k)}$:
\begin{equation*}
\left\|h_{(k)}\right\|_{C_{t}} \leq 1. \label{4.18}\tag{3.7}
\end{equation*}
\end{Lemma}

\subsection{Spatial building blocks.}\label{D.3.1}\ \ \ \
Define $\varphi, \psi:\mathbb{R}\rightarrow \mathbb{R}$ as smooth, mean-free functions supported on a ball of radius 1, satisfying the following conditions
\begin{equation*}
\frac{1}{2 \pi} \int_{\mathbb{R}} \varphi^{2}(x) \mathrm{d} x=1\quad \text{and}\quad \frac{1}{2 \pi} \int_{\mathbb{R}} \psi^{2}(x) \mathrm{d} x=1, \quad \operatorname{supp} \varphi,\ \psi \subseteq[-1,1]. \label{4.2}\tag{3.8}
\end{equation*}

The corresponding rescaled cut-off functions are defined as follows:
\begin{equation*}
\varphi_{r_{\|}}(x):=r_{\|}^{-\frac{1}{2}}\varphi\left(\frac{x}{r_{\|}}\right), \quad \psi_{r_{\perp}}(x):=r_{\perp}^{-\frac{1}{2}} \psi\left(\frac{x}{r_{\perp}}\right) .\label{4.3}\tag{3.9}
\end{equation*}
In the scaling process, the functions $\varphi_{r_{\|}}$ and $\psi_{r_{\perp}}$ are supported within the balls of radius $r_{\|}$ and $r_{\perp}$ in $\mathbb{R}$, respectively. For notational convenience, we also extend the definition of $\varphi_{r_{\|}}$ and $\psi_{r_{\perp}}$ to periodic functions on the torus $\mathbb{T}$.

For each $k \in \Lambda$, let $\left(k, k_{1}\right)$ denotes an orthonormal basis. Then the accelerating jets are defined by
\begin{align*}
\mathbf{W}_{k}:=-\varphi_{r_{\|}}\left(\lambda r_{\perp} N_{\Lambda}k_{1} \cdot \left(x-\alpha_{(k)}\right)+N_{\Lambda}\phi_{(k)}(t)\right) \psi_{r_{\perp}}^{\prime}\left(\lambda r_{\perp} N_{\Lambda} k \cdot\left(x-\alpha_{(k)}\right)\right) k_{1}, \label{4.4}\tag{3.10}
\end{align*}
where the quantity $N_{\Lambda}$ is defined according to \eqref{7.2}, while the function $\phi_{(k)}(t)$ is determined by the relationship $\phi_{(k)}^{\prime}(t)=\mu g_{(k)}$, and the term $g_{(k)}$ has been specified in \eqref{4.15}. By design, $\mathbf{W}_{k}$ is $(\mathbb{T}/\lambda r_{\perp})^2$-periodic in space and $\sigma^{-1}$-periodic in time.

 The parameters $r_{\|}$ and $r_{\perp}$ serve as indicators of the concentration effect exerted by accelerating jets. Additionally, $\mu$ represents the acceleration of the flow within the building block, while $\phi_{(k)}$ establishes a relationship between the velocity of the moving support sets and the intermittent oscillator $g_{(k)}$. By employing an appropriate discretization of temporal velocity, it is possible to avoid collisions among different $\mathbf{W}_{k}$ in 2D. The shifts $\alpha_{(k)} \in \mathbb{R}^{2}$ are carefully selected to ensure that $\mathbf{W}_{k}$ and $\mathbf{W}_{k^{\prime}}$ have disjoint supports whenever $k \neq k^{\prime}$. The existence of such $\alpha_{(k)}$ can be assured by maintaining sufficiently small values for $r_{\perp}$ and $r_{\|}$ (as demonstrated in \cite{TBV19}).

For brevity, we set
\begin{align*}
& \varphi_{\left(k_{1}\right)}(x):=\varphi_{r_{\|}}\left(\lambda r_{\perp} N_{\Lambda}k_{1} \cdot \left(x-\alpha_{(k)}\right)+N_{\Lambda}\phi_{(k)}(t)\right)=\varphi_{r_{\|}}\left(\lambda r_{\perp} N_{\Lambda}x_{k}+N_{\Lambda}\phi_{(k)}(t)\right), \\
& \psi_{(k)}(x):=\psi_{r_{\perp}}\left(\lambda r_{\perp} N_{\Lambda} k \cdot\left(x-\alpha_{(k)}\right)\right)=\psi_{r_{\perp}}\left(\lambda r_{\perp} N_{\Lambda} y_{k}\right),
\end{align*}
where the pair $(x_{k}, y_{k})$ is employed to represent the spatial coordinates in $\mathbb{R}^{2}$. Therefore, the accelerating jets is simplified as
\begin{equation*}
\mathbf{W}_{k}=-\varphi_{\left(k_{1}\right)} \psi_{(k)}^{\prime} k_{1}, \quad k \in \Lambda . \label{4.5}\tag{3.11}
\end{equation*}
Since $\mathbf{W}_{k}$ is not divergence-free, we introduce the corrector $\mathbf{W}_{k}^{c}$  as
\begin{equation*}
\mathbf{W}_{k}^{c}:=\frac{r_{\perp}}{r_{\|}}\varphi_{\left(k_{1}\right)}^{\prime} \psi_{(k)} k. \label{4.6}\tag{3.12}
\end{equation*}
Furthermore, we define periodic potentials $\boldsymbol{\Psi}_{k} \in C_{c}^{\infty}\left(\mathbb{R}^{2}\right)$ as
\begin{equation*}
\boldsymbol{\Psi}_{k} =r_{\perp}\varphi_{\left(k_{1}\right)} \psi_{(k)}. \label{4.7}\tag{3.13}
\end{equation*}
By straightforward computations,
$$
\begin{aligned}
(\lambda r_{\perp}N_{\Lambda})^{-1} \nabla^{\perp} \boldsymbol{\Psi}_{k} & =(\lambda r_{\perp}N_{\Lambda})^{-1}(-\partial_{y_{k}} \boldsymbol{\Psi}_{k} k_{1}+\partial_{x_{k}} \boldsymbol{\Psi}_{k} k) \\
& =\mathbf{W}_{k}+\mathbf{W}_{k}^{(c)}.
\end{aligned}
$$
Also, we have the important identities
\begin{equation*}
\partial_{t}\left|\mathbf{W}_{k}\right|^{2} k_{1}=(\lambda r_{\perp})^{-1} \mu g_{(k)} \operatorname{div}\left(\mathbf{W}_{k} \otimes \mathbf{W}_{k}\right), \label{4.8}\tag{3.14}
\end{equation*}
and
\begin{equation*}
\partial_{t} \boldsymbol{\Psi}_{k}=(\lambda r_{\perp})^{-1} \mu g_{(k)}\left(k_{1} \cdot \nabla\right) \boldsymbol{\Psi}_{k}. \label{4.9}\tag{3.15}
\end{equation*}

The following lemma is the crucial estimates pertaining to the accelerating jets.

\begin{Lemma}\label{D.1}(Estimates of accelerating jets). For any $p \in[1, \infty]$, $N, M \in \mathbb{N}$, the following bounds hold,
\begin{align*}
& \left\|\nabla^{N} \partial_{t}^{M} \varphi_{\left(k_{1}\right)}\right\|_{C_{t} L_{x}^{p}}+\left\|\nabla^{N} \partial_{t}^{M} \varphi_{\left(k_{1}\right)}^{\prime}\right\|_{C_{t} L_{x}^{p}}\lesssim r_{\|}^{\frac{1}{p}-\frac{1}{2}}\left(\frac{r_{\perp} \lambda}{r_{\|}}\right)^{N}\left(\frac{\mu \tau^{\frac{1}{2}}}{r_{\|}}\right)^{M},  \label{4.10}\tag{3.16}\\
& \left\|\nabla^{N} \psi_{(k)}\right\|_{L_{x}^{p}}+\left\|\nabla^{N} \psi_{(k)}^{\prime}\right\|_{L_{x}^{p}} \lesssim r_{\perp}^{\frac{1}{p}-\frac{1}{2}} \lambda^{N}. \label{4.11}\tag{3.17}
\end{align*}
Furthermore, it holds that
\begin{align*}
& \left\|\nabla^{N} \partial_{t}^{M} \mathbf{W}_{k}\right\|_{C_{t} L_{x}^{p}}+\frac{r_{\|}}{r_{\perp}}\left\|\nabla^{N} \partial_{t}^{M} \mathbf{W}_{k}^{c}\right\|_{C_{t} L_{x}^{p}}+r_{\perp}^{-1}\left\|\nabla^{N} \partial_{t}^{M} \boldsymbol{\Psi}_{k}\right\|_{C_{t} L_{x}^{p}} \\
& \quad \lesssim (r_{\perp}r_{\|})^{\frac{1}{p}-\frac{1}{2}}\lambda^{N}\left(\frac{\mu\tau^{\frac{1}{2}}}{r_{\|}}\right)^{M}, \label{4.12} \tag{3.18}
\end{align*}
where the implicit constants $C$ are independent of $r_{\perp}, r_{\|}, \lambda, \mu$, and $k \in \Lambda$.
\end{Lemma}
{\bf Proof.} Combining the relationship $\phi_{(k)}^{\prime}(t)=\mu g_{(k)}$ with \eqref{4.17}, due to the fact  $\frac{\mu\tau^{\frac{1}{2}}}{r_{\|}}\gg\sigma\tau$, we deduce that
$$\|\nabla^{N} \partial_{t}^{M}\varphi_{(k_{1})}\|_{C_{t}L_{x}^{p}}\lesssim
r_{\|}^{\frac{1}{p}-\frac{1}{2}}\left(\frac{r_{\perp} \lambda}{r_{\|}}\right)^{N}\sum\limits_{i=1}^{M}(\frac{\mu\tau^{\frac{1}{2}}}{r_{\|}})^{i}
(\sigma\tau)^{M-i}\lesssim r_{\|}^{\frac{1}{p}-\frac{1}{2}}\left(\frac{r_{\perp} \lambda}{r_{\|}}\right)^{N}\left(\frac{\mu \tau^{\frac{1}{2}}}{r_{\|}}\right)^{M},$$
\eqref{4.10} is verified. Similarly, \eqref{4.11} also holds true. Then by applying the Fubini theorem to combine the estimates obeyed by $\varphi_{(k_{1})}$ and $\psi_{(k)}$ (which are both 1D function), to obtain estimates for the 2D functions $  \mathbf{W}_{k}, \mathbf{W}_{k}^{c}$ and $\mathbf{\Psi}_{k}$. Therefore the Lemma \ref{D.1} is verified.

\begin{Remark}\label{D.2} After defining the temporal concentrated function $g_{(k)}$ and 2D accelerating jets $\mathbf{W}_{k}$, we will utilize $g_{(k)} \mathbf{W}_{k}$ as our ``building block" in spatial-temporal convex integration scheme. Owing to the temporal concentration of $g_{(k)}$, there will be no overlap or interference among the distinct elements $\mathbf{W}_k$. Consequently, the supports of all $g_{(k)} \mathbf{W}_{k}$ are mutually disjoint if $k\neq k'$ on $\mathbb{T}^2 \times[0,T]$.
\end{Remark}

\subsection{Velocity perturbations.}\label{C.3.3}\ \ \ \
In the following sections, we will divide the construction of velocity perturbations into four components, encompassing the principal perturbation, incompressible corrector, and two temporal correctors. Firstly, we will provide the specific form of the amplitudes of these perturbations, which has been presented in previous research works \cite{TBV19,YL22,LTP20}. Here, we state that the core purpose of its establishment is achieving a cancellation between the low frequency part of the nonlinearity and the old Reynolds stress $\mathring{\widetilde{\mathsf{R}}}_q$.

\noindent {\bf  Amplitudes of perturbations.}
Choosing $\mathscr{X} :[0, \infty) \rightarrow \mathbb{R}$ to be a smooth cut-off function that satisfies
\[
\mathscr{X} (z)= \begin{cases}1, & 0 \leq z \leq 1,  \label{4.19}\tag{3.19}\\ z, & z \geq 2,\end{cases}
\]
and
\begin{equation*}
\frac{1}{2} z \leq \mathscr{X} (z) \leq 2 z, \quad \text { for } \quad z \in(1,2). \label{4.20}\tag{3.20}
\end{equation*}
Define
\begin{equation*}
\rho_{u}(t, x):=2 C_{\mathsf{R}}^{-1} \lambda_q^{-\frac{2}{3}\varepsilon_R} \delta_{q+1} \mathscr{X} \left(\frac{\left|\mathring{\widetilde{\mathsf{R}}}_q(t, x)\right|}{\lambda_q^{-\frac{2}{3}\varepsilon_R} \delta_{q+1}}\right), \label{4.21}\tag{3.21}
\end{equation*}
where $C_{\mathsf{R}}$ is the positive constant defined in the Geometric Lemma \ref{G.1}. By \eqref{4.19}, \eqref{4.20} and \eqref{4.21},
\begin{equation*}
\left|\frac{\mathring{\widetilde{\mathsf{R}}}_q}{\rho_{u}}\right|=\left|\frac{\mathring{\widetilde{\mathsf{R}}}_q}{2 C_{\mathsf{R}}^{-1} \lambda_{q}^{-\frac{2}{3}\varepsilon_R} \delta_{q+1} \mathscr{X} \left(\lambda_{q}^{\frac{2}{3}\varepsilon_R} \delta_{q+1}^{-1}\left|\mathring{\widetilde{\mathsf{R}}}_q\right|\right)}\right| \leq C_{\mathsf{R}}, \label{4.22}\tag{3.22}
\end{equation*}
and for any $p \in[1, \infty]$,
\begin{align*}
& \rho_{u} \geq C_{\mathsf{R}}^{-1} \lambda_{q}^{-\frac{2}{3}\varepsilon_R} \delta_{q+1},\label{4.23}\tag{3.23}
\end{align*}
\begin{align*}
\|\rho_{u}\|_{L_{t, x}^p} \lesssim C_{\mathsf{R}}^{-1}\left(\lambda_q^{-\frac{2}{3}\varepsilon_R} \delta_{q+1}+\left\|\mathring{\widetilde{\mathsf{R}}}_q\right\|_{L_{t, x}^p}\right).\label{4.24}\tag{3.24}
\end{align*}
Furthermore, combining \eqref{3.7}, \eqref{4.23} with the standard H\"{o}lder's estimates (see \cite{TB15}), for $1 \leq N \leq 9$,
\begin{gather*}
\|\rho_{u}\|_{C_{t, x}} \lesssim \vartheta_{q+1}^{-2}, \quad\|\rho_{u}\|_{C_{t, x}^{N}} \lesssim \vartheta_{q+1}^{-3 N} \label{4.25}\tag{3.25}\\
\left\|\rho_{u}^{1 / 2}\right\|_{C_{t, x}} \lesssim \vartheta_{q+1}^{-1}, \quad\left\|\rho_{u}^{1 / 2}\right\|_{C_{t, x}^{N}} \lesssim \vartheta_{q+1}^{-3 N},  \label{4.26}\tag{3.26}\\
\left\|\rho_{u}^{-1}\right\|_{C_{t, x}} \lesssim \vartheta_{q+1}^{-1}, \quad\left\|\rho_{u}^{-1}\right\|_{C_{t, x}^{N}} \lesssim \vartheta_{q+1}^{-3 N}, \label{4.27}\tag{3.27}
\end{gather*}
where the implicit constants are independent of the variable $q$. To ensure temporal compatibility between the perturbations and the concentrated Reynolds stress $\mathring{\widetilde{\mathsf{R}}}_{q}$ discussed in Section \ref{B}, we employ a smoothly temporal cut-off function $f_{u}:[0, T] \rightarrow[0,1]$ that fulfills
\begin{itemize}
  \item $0 \leq f_{u} \leq 1$ and $f \equiv 1$ on $\operatorname{supp}_{t} \mathring{\widetilde{\mathsf{R}}}_{q} ;$
  \item $\operatorname{supp}_t f_{u} \subseteq N_{\vartheta_{q+1} / 2}\left(\operatorname{supp}_t \mathring{\widetilde{\mathsf{R}}}_{q}\right);$
  \item $\|f_{u}\|_{C_{t}^{N}} \lesssim \vartheta_{q+1}^{-N}, \quad 1 \leq N \leq 9$.
\end{itemize}

Now, we provide the specific form of amplitude as follows
\begin{equation*}
a_{(k)}(t, x):=\rho_{u}^{\frac{1}{2}}(t, x) f_{u}(t) \gamma_{(k)}\left(\operatorname{Id}-\frac{\mathring{\widetilde{\mathsf{R}}}_{q}(t, x)}{\rho_{u}(t, x)}\right), \quad k \in \Lambda, \label{4.28}\tag{3.28}
\end{equation*}
where the symbols $\gamma_{(k)}$ and $\Lambda$ are given in the Geometric Lemma \ref{G.1}. Furthermore, we possess the analytic estimates for the amplitudes, the derivation of which follows a similar proof in \cite{YLZ22}.
\begin{Lemma}\label{D.5}(Estimates of amplitudes)
For $1 \leq N \leq 9, k \in \Lambda$, it holds
\begin{align*}
\left\|a_{(k)}\right\|_{L_{t, x}^{2}} & \lesssim \delta_{q+1}^{\frac{1}{2}},  \label{4.29}\tag{3.29}\\
\left\|a_{(k)}\right\|_{C_{t, x}} & \lesssim \vartheta_{q+1}^{-1}, \quad\left\|a_{(k)}\right\|_{C_{t, x}^{N}} \lesssim \vartheta_{q+1}^{-5 N} , \label{4.30}\tag{3.30}
\end{align*}
where the implicit constants do not rely on the variable $q$.
\end{Lemma}
\noindent {\bf Velocity Perturbations.}
After the selection of ``building block'' and the introduction of amplitudes, we now prepare construct the velocity perturbations. Firstly, the principal part of the velocity perturbations, denoted as $w_{q+1}^{(p)}$, is defined as
\begin{equation*}
w_{q+1}^{(p)}:=\sum_{k \in \Lambda} a_{(k)} g_{(k)} \mathbf{W}_{k}. \label{4.31}\tag{3.31}
\end{equation*}

By applying Geometric Lemma \ref{G.1} and utilizing the expression given in \eqref{4.28}, the concentrated Reynolds stress $\mathring{\widetilde{\mathsf{R}}}_q$ can be eliminated by the zero frequency part of $w_{q+1}^{(p)} \otimes w_{q+1}^{(p)}$, i.e.,
\begin{align*}
w_{q+1}^{(p)} \otimes w_{q+1}^{(p)}+\mathring{\widetilde{\mathsf{R}}}_q&=\rho_{u} f_{u}^2 \operatorname{Id}+\sum_{k \in \Lambda} a_{(k)}^2 g_{(k)}^2 \mathbb{P}_{\neq 0}\left(\mathbf{W}_{k} \otimes \mathbf{W}_{k}\right)\\
& +\sum_{k \in \Lambda} a_{(k)}^{2}\left(g_{(k)}^{2}-1\right) \fint_{\mathbb{T}^{2}} \mathbf{W}_{k} \otimes \mathbf{W}_{k} \mathrm{d} x .\label{4.32}\tag{3.32}
\end{align*}
In the given context,  $\mathbb{P}_{\neq 0}$ represents the spatial projection onto the set of nonzero Fourier modes. Since the primary component of the perturbation $w_{q+1}^{(p)}$ lacks divergence-free property, we introduce an incompressibility corrector defined as follows:
\begin{equation*}
w_{q+1}^{(c)}:=\sum_{k \in \Lambda} a_{(k)} g_{(k)}{\bf W}_{k}^{c}+(\lambda r_{\perp}N_{\Lambda})^{-1}\nabla^{\bot} a_{(k)} g_{(k)}{\bf\Psi}_{k}, \label{4.33}\tag{3.33}
\end{equation*}
where $\mathbf{W}_{k}^{c}$ and ${\bf\Psi}_{k} $ are respectively defined by \eqref{4.6} and \eqref{4.7}. subsequently,
\begin{align*}
w_{q+1}^{(p)}+w_{q+1}^{(c)}&=\sum_{k \in \Lambda}a_{(k)} g_{(k)}({\bf W}_{k}+{\bf W}_{k}^{c})
+(\lambda r_{\perp}N_{\Lambda})^{-1}\nabla^{\bot} a_{(k)} g_{(k)}{\bf\Psi}_{k}\\
&=(\lambda r_{\perp}N_{\Lambda})^{-1}\sum_{k \in \Lambda} g_{(k)}[a_{(k)}\nabla^{\bot}{\bf\Psi}_{k}
+\nabla^{\bot} a_{(k)}{\bf\Psi}_{k}]\\
&= (\lambda r_{\perp}N_{\Lambda})^{-1}\sum_{k \in \Lambda}\nabla^{\bot}[a_{(k)} g_{(k)}{\bf\Psi}_{k}],\label{4.34}\tag{3.34}
\end{align*}
and thus
\begin{equation*}
\operatorname{div}\left(w_{q+1}^{(p)}+w_{q+1}^{(c)}\right)=0. \label{4.35}\tag{3.35}
\end{equation*}

To address the high-frequency spatial and temporal errors in \eqref{4.32}, we introduce two additional  temporal correctors. The first temporal corrector $w_{q+1}^{(t)}$ is defined as
\begin{equation*}
w_{q+1}^{(t)}:=-\mu^{-1} (\lambda r_{\perp})\sum_{k \in \Lambda} \mathbb{P}_{H} \mathbb{P}_{\neq 0}\left(a_{(k)}^{2} g_{(k)} \varphi_{\left(k_{1}\right)} ^{2} (\psi_{(k)}^{\prime})^{2} k_{1}\right) \label{4.36}\tag{3.36}
\end{equation*}
to balance the high spatial frequency oscillations:
\begin{align*}
\partial_{t} w_{q+1}^{(t)}&+\sum_{k \in \Lambda} \mathbb{P}_{\neq 0}\left(a_{(k)}^{2} g_{(k)}^{2} \operatorname{div}\left(\mathbf{W}_{k} \otimes \mathbf{W}_{k}\right)\right)=\\
 & \left(\nabla \Delta^{-1} \operatorname{div}\right) \mu^{-1}\lambda r_{\perp} \sum_{k \in \Lambda} \mathbb{P}_{\neq 0} \partial_{t}\left(a_{(k)}^{2} g_{(k)} \varphi_{\left(k_{1}\right)} ^{2} (\psi_{(k)}^{\prime})^{2} k_{1}\right) \\
& -\mu^{-1}\lambda r_{\perp} \sum_{k \in \Lambda} \mathbb{P}_{\neq 0}\left(\partial_{t}\left(a_{(k)}^{2} g_{(k)}\right) \varphi_{\left(k_{1}\right)} ^{2} (\psi_{(k)}^{\prime})^{2} k_{1}\right) , \label{4.37}\tag{3.37}
\end{align*}
where the first term on the right-hand-side can be processed by Helmholtz-Leray projector $\mathbb{P}_{H}$.

The other temporal corrector $w_{q+1}^{(o)}$ is defined by
\begin{equation*}
w_{q+1}^{(o)}:=-\sigma^{-1} \sum_{k \in \Lambda} \mathbb{P}_{H} \mathbb{P}_{\neq 0}\left(h_{(k)} \fint_{\mathbb{T}^{2}} \mathbf{W}_{k} \otimes \mathbf{W}_{k} \mathrm{d} x \nabla\left(a_{(k)}^{2}\right)\right) \label{4.38}\tag{3.38}
\end{equation*}
to balance the high temporal frequency oscillations:
\begin{align*}
& \partial_{t} w_{q+1}^{(o)}+\sum_{k \in \Lambda} \mathbb{P}_{\neq 0}\left(\left(g_{(k)}^{2}-1\right) \fint_{\mathbb{T}^{2}} \mathbf{W}_{k} \otimes \mathbf{W}_{k} \mathrm{d} x \nabla\left(a_{(k)}^{2}\right)\right) \\
= & \left(\nabla \Delta^{-1} \operatorname{div}\right) \sigma^{-1} \sum_{k \in \Lambda} \mathbb{P}_{\neq 0} \partial_{t}\left(h_{(k)} \fint_{\mathbb{T}^{2}} \mathbf{W}_{k} \otimes \mathbf{W}_{k} \mathrm{d} x \nabla\left(a_{(k)}^{2}\right)\right) \\
& -\sigma^{-1} \sum_{k \in \Lambda} \mathbb{P}_{\neq 0}\left(h_{(k)} \fint_{\mathbb{T}^{2}} \mathbf{W}_{k} \otimes \mathbf{W}_{k} \mathrm{d} x \partial_{t} \nabla\left(a_{(k)}^{2}\right)\right) . \label{4.39}\tag{3.39}
\end{align*}
where the right-hand-side above only remain the low frequency part $\partial_{t} \nabla\left(a_{(k)}^{2}\right)$ and the harmless pressure term.

In summary, the velocity perturbation $w_{q+1}$ at level $q+1$ is defined as the sum of four components:
\begin{equation*}
w_{q+1}:=w_{q+1}^{(p)}+w_{q+1}^{(c)}+w_{q+1}^{(t)}+w_{q+1}^{(o)} . \label{4.40}\tag{3.40}
\end{equation*}
These components are constructed in such a way that $w_{q+1}$ possesses both mean-free and divergence-free properties. Subsequently, the velocity field at level $q+1$ is determined by
\begin{equation*}
u_{q+1}:=\widetilde{u}_{q}+w_{q+1}, \label{4.41}\tag{3.41}
\end{equation*}
where $\widetilde{u}_{q}$ is obtained from the gluing stage as detailed in \eqref{3.21}.

Furthermore, the following Lemma \ref{D.6} presents the crucial estimates associated with these velocity perturbations.
\begin{Lemma}\label{D.6} (Estimates of velocity perturbations). Given any $p \in(1, \infty), \gamma \in[1, \infty]$ and integers $0 \leq N \leq 7$, it holds that
\begin{align*}
&\left\|\nabla^{N} w_{q+1}^{(p)}\right\|_{L_{t}^{\gamma} L_{x}^{p}} \lesssim \vartheta_{q+1}^{-1} \lambda^{N} (r_{\perp} r_{\|})^{\frac{1}{p}-\frac{1}{2}} \tau^{\frac{1}{2}-\frac{1}{\gamma}},  \label{4.42}\tag{3.42}\\
&\left\|\nabla^{N} w_{q+1}^{(c)}\right\|_{L_{t}^{\gamma} L_{x}^{p}} \lesssim \vartheta_{q+1}^{-6} \lambda^{N} (r_{\perp} r_{\|})^{\frac{1}{p}-\frac{1}{2}} \tau^{\frac{1}{2}-\frac{1}{\gamma}}  \frac{r_{\perp}}{ r_{\|}},\label{4.43}\tag{3.43}\\
&\left\|\nabla^{N} w_{q+1}^{(t)}\right\|_{L_{t}^{\gamma} L_{x}^{p}} \lesssim \vartheta_{q+1}^{-2} \lambda^{N} \mu^{-1}\lambda r_{\perp} (r_{\perp}r_{\|})^{\frac{1}{p}-1} \tau^{\frac{1}{2}-\frac{1}{\gamma}},  \label{4.44}\tag{3.44}\\
&\left\|\nabla^{N} w_{q+1}^{(o)}\right\|_{L_{t}^{\gamma} L_{x}^{p}} \lesssim \vartheta_{q+1}^{-5N-7} \sigma^{-1}. \label{4.45}\tag{3.45}
\end{align*}
Especially, for integers $1 \leq N \leq 7$,
\begin{align*}
& \left\|w_{q+1}^{(p)}\right\|_{L_{t}^{\infty} H_{x}^{N}}+\left\|w_{q+1}^{(c)}\right\|_{L_{t}^{\infty} H_{x}^{N}}+\left\|w_{q+1}^{(t)}\right\|_{L_{t}^{\infty} H_{x}^{N}}+\left\|w_{q+1}^{(o)}\right\|_{L_{t}^{\infty} H_{x}^{N}} \lesssim \lambda^{N+2},  \label{4.46}\tag{3.46}\\
& \left\|\partial_{t} w_{q+1}^{(p)}\right\|_{L_{t}^{\infty} H_{x}^{N}}+\left\|\partial_{t} w_{q+1}^{(c)}\right\|_{L_{t}^{\infty} H_{x}^{N}}+\left\|\partial_{t} w_{q+1}^{(t)}\right\|_{L_{t}^{\infty} H_{x}^{N}}+\left\|\partial_{t} w_{q+1}^{(o)}\right\|_{L_{t}^{\infty} H_{x}^{N}} \lesssim \lambda^{N+5}, \label{4.47}\tag{3.47}
\end{align*}
where the implicit constants do not rely on $\lambda$.
\end{Lemma}
{\bf Proof.} By \eqref{3.1}, \eqref{4.12}, \eqref{4.17}, \eqref{4.32}, \eqref{4.33} and Lemma \ref{D.5}, for any $p \in(1, \infty)$,
$$
\begin{aligned}
\left\|\nabla^{N} w_{q+1}^{(p)}\right\|_{L_{t}^{\gamma} L_{x}^{p}} & \lesssim \sum_{k \in \Lambda} \sum_{N_{1}+N_{2}=N}\left\|a_{(k)}\right\|_{C_{t, x}^{N_{1}}}\left\|g_{(k)}\right\|_{L_{t}^{\gamma}}\left\|\nabla^{N_{2}} \mathbf{W}_{k}\right\|_{C_{t} L_{x}^{p}} \\
& \lesssim \vartheta_{q+1}^{-1} \lambda^{N} (r_{\perp} r_{\|})^{\frac{1}{p}-\frac{1}{2}} \tau^{\frac{1}{2}-\frac{1}{\gamma}},
\end{aligned}
$$
and
$$
\begin{aligned}
&\left\|\nabla^{N} w_{q+1}^{(c)}\right\|_{L_{t}^{\gamma} L_{x}^{p}} \\
& \lesssim \sum_{k \in \Lambda}\left\|g_{(k)}\right\|_{L_{t}^{\gamma}} \sum_{N_{1}+N_{2}=N}\left(\left\|a_{(k)}\right\|_{C_{t, x}^{N_{1}}}\left\|\nabla^{N_{2}} \mathbf{W}_{k}^{c}\right\|_{C_{t} L_{x}^{p}}+(\lambda r_{\perp})^{-1}\left\|a_{(k)}\right\|_{C_{t, x}^{N_{1}+1}}\left\|\nabla^{N_{2}} \boldsymbol{\Psi}_{k} \right\|_{C_{t} L_{x}^{p}}\right) \\
& \lesssim \sum_{N_{1}+N_{2}=N} \tau^{\frac{1}{2}-\frac{1}{\gamma}}\left(\vartheta_{q+1}^{-5 N_{1}-1} \lambda^{N_{2}} \frac{r_{\perp}}{r_{\|}}(r_{\perp} r_{\|})^{\frac{1}{p}-\frac{1}{2}}+(\lambda r_{\perp})^{-1}\vartheta_{q+1}^{-5 N_{1}-6} \lambda^{N_{2}} r_{\perp} ( r_{\perp}r_{\|})^{\frac{1}{p}-\frac{1}{2}}\right) \\
& \lesssim \tau^{\frac{1}{2}-\frac{1}{\gamma}}\left(\vartheta_{q+1}^{-1} \lambda^{N} \frac{r_{\perp}}{r_{\|}}(r_{\perp} r_{\|})^{\frac{1}{p}-\frac{1}{2}}+\vartheta_{q+1}^{-6} \lambda^{N-1}(r_{\perp}r_{\|})^{\frac{1}{p}-\frac{1}{2}}\right) \\
& \lesssim \vartheta_{q+1}^{-6} \lambda^{N} (r_{\perp}r_{\|})^{\frac{1}{p}-\frac{1}{2}} \tau^{\frac{1}{2}-\frac{1}{\gamma}}\frac{r_{\perp}}{r_{\|}} ,
\end{aligned}
$$
then, the \eqref{4.42}, \eqref{4.43} are verified.

Regarding the temporal correctors, based on \eqref{4.36} and \eqref{4.38}, Lemmas \ref{D.1}, \ref{D.3}, and \ref{D.5}, as well as the boundedness of the operators $\mathbb{P}{\neq 0}$ and $\mathbb{P}_{H}$ in the space $L_{x}^{p}$, we deduce that
$$
\begin{aligned}
&\left\|\nabla^{N} w_{q+1}^{(t)}\right\|_{L_{t}^{\gamma} L_{x}^{p}}\\
& \lesssim \mu^{-1} \lambda r_{\perp}\sum_{k \in \Lambda}\left\|g_{(k)}\right\|_{L_{t}^{\gamma}} \sum_{N_{1}+N_{2}+N_{3}=N}\left\|\nabla^{N_{1}}\left(a_{(k)}^{2}\right)\right\|_{C_{t, x}}\left\|\nabla^{N_{2}}\left(\varphi_{\left(k_{1}\right)}^{2}\right)\right\|_{C_{t} L_{x}^{p}}\left\|\nabla^{N_{3}}\left(\psi_{(k)}^{\prime}\right)^{2}\right\|_{L_{x}^{p}} \\
& \lesssim \mu^{-1} \lambda r_{\perp}\tau^{\frac{1}{2}-\frac{1}{\gamma}} \sum_{N_{1}+N_{2}+N_{3}=N} \vartheta_{q+1}^{-5 N_{1}-2} (\lambda\frac{r_{\perp}}{r_{\|}})^{N_{2}} \lambda^{N_{3}} (r_{\|}r_{\perp})^{\frac{1}{p}-1} \\
& \lesssim \vartheta_{q+1}^{-2}\mu^{-1}\lambda r_{\perp} \lambda^{N} (r_{\|}r_{\perp})^{\frac{1}{p}-1} \tau^{\frac{1}{2}-\frac{1}{\gamma}} ,
\end{aligned}
$$
and
$$
\left\|\nabla^{N} w_{q+1}^{(o)}\right\|_{L_{t}^{\gamma} L_{x}^{p}} \lesssim \sigma^{-1} \sum_{k \in \Lambda}\left\|h_{(k)}\right\|_{C_{t}}\left\|\nabla^{N+1}\left(a_{(k)}^{2}\right)\right\|_{C_{t, x}} \lesssim \vartheta_{q+1}^{-5N-7} \sigma^{-1}
$$
which yields \eqref{4.44} and \eqref{4.45}.

To obtain the $L_{t}^{\infty} H_{x}^{N}$-estimate of velocity perturbations, by using \eqref{4.42}-\eqref{4.45},
\begin{align*}
& \left\|w_{q+1}^{(p)}\right\|_{L_{t}^{\infty} H_{x}^{N}}+\left\|w_{q+1}^{(c)}\right\|_{L_{t}^{\infty} H_{x}^{N}}+\left\|w_{q+1}^{(t)}\right\|_{L_{t}^{\infty} H_{x}^{N}}+\left\|w_{q+1}^{(o)}\right\|_{L_{t}^{\infty} H_{x}^{N}} \\
\lesssim & \vartheta_{q+1}^{-1} \lambda^{N} \tau^{\frac{1}{2}}+\vartheta_{q+1}^{-6} \lambda^{N} \frac{r_{\perp}}{r_{\|}} \tau^{\frac{1}{2}}+\vartheta_{q+1}^{-2} \lambda r_{\perp}\lambda^{N} \mu^{-1} (r_{\perp} r_{\|})^{-\frac{1}{2}} \tau^{\frac{1}{2}}+\vartheta_{q+1}^{-5N-7} \sigma^{-1} \label{4.48}\tag{3.48}\\
\lesssim & \vartheta_{q+1}^{-1} \lambda^{N+2 \alpha-2+8 \varepsilon}+\vartheta_{q+1}^{-6} \lambda^{N+2 \alpha-2}+\vartheta_{q+1}^{-2} \lambda^{N+2 \alpha-2+2 \varepsilon}+\vartheta_{q+1}^{-5N-7} \lambda^{-2 \varepsilon}\\
\lesssim &\lambda^{N+2},
\end{align*}
where the last step is validated by  \eqref{2.3} and \eqref{2.4}, ultimately leading to the confirmation of \eqref{4.46}.

The last step is to prove \eqref{4.47}, by using \eqref{2.3}, \eqref{4.1} and Lemmas \ref{D.1}, \ref{D.3}, and \ref{D.5}, we have
\begin{align*}
\left\|\partial_{t} w_{q+1}^{(p)}\right\|_{L_{t}^{\infty} H_{x}^{N}} & \lesssim \sum_{k \in \Lambda}\left\|a_{(k)}\right\|_{C_{t, x}^{N+1}} \sum_{M_{1}+M_{2}=1}\left\|\partial_{t}^{M_{1}} g_{(k)}\right\|_{L_{t}^{\infty}}\left\|\partial_{t}^{M_{2}} \mathbf{W}_{k}\right\|_{L_{t}^{\infty} H_{x}^{N}} \\
& \lesssim \sum_{M_{1}+M_{2}=1} \vartheta_{q+1}^{-5 N-6} \sigma^{M_{1}} \tau^{M_{1}+\frac{1}{2}} \lambda^{N}\left(\frac{\mu \tau^{\frac{1}{2}}}{r_{\|}}\right)^{M_{2}} \\
& \lesssim \vartheta_{q+1}^{-5 N-6} \lambda^{N}\frac{\mu \tau}{r_{\|}} \label{4.49}\tag{3.49}
\end{align*}
and
\begin{align*}
&\left\|\partial_{t} w_{q+1}^{(c)}\right\|_{L_{t}^{\infty} H_{x}^{N}} \\
& \lesssim \sum_{k \in \Lambda}\left\|a_{(k)}\right\|_{C_{t, x}^{N+1}} \sum_{M_{1}+M_{2}=1}\left\|g_{(k)}\right\|_{C_{t}^{M_{1}}}\left(\left\|\partial_{t}^{M_{2}} \mathbf{W}_{k}^{c}\right\|_{L_{t}^{\infty} H_{x}^{N}}+(\lambda r_{\perp})^{-1}\left\|\partial_{t}^{M_{2}} {\bf\Psi}_{(k)}\right\|_{L_{t}^{\infty} H_{x}^{N}}\right) \\
& \lesssim \sum_{M_{1}+M_{2}=1} \vartheta_{q+1}^{-5 N-6} \sigma^{M_{1}} \tau^{M_{1}+\frac{1}{2}}\left(\frac{ \mu \tau^{\frac{1}{2}}}{r_{\|}}\right)^{M_{2}}\left(\lambda^{N-1}+\frac{r_{\perp}}{r_{\|}} \lambda^{N}\right) \\
& \lesssim \vartheta_{q+1}^{-5 N-6} \lambda^{N} \frac{\mu \tau}{r_{\|}} \frac{r_{\perp}}{r_{\|}} . \label{4.50}\tag{3.50}
\end{align*}

Furthermore, due to the boundedness of $\mathbb{P}_{H} \mathbb{P}_{\neq 0}$ in $H_{x}^{N}$, we can proceed analogously to the previous argumentation and invoke Lemmas \ref{D.1}, \ref{D.3}, and \ref{D.5} once again, then it yields that
\begin{align*}
&\left\|\partial_{t} w_{q+1}^{(t)}\right\|_{L_{t}^{\infty} H_{x}^{N}}\\
& \lesssim \mu^{-1} \lambda r_{\perp}\sum_{k \in \Lambda}\left\|\partial_{t}\left(a_{(k)}^{2} g_{(k)} \varphi_{\left(k_{1}\right)} ^{2} (\psi_{(k)}^{\prime})^{2}\right)\right\|_{L_{t}^{\infty} H_{x}^{N}} \\
& \lesssim \mu^{-1} \lambda r_{\perp}\sum_{k \in \Lambda}\left\|a_{(k)}^{2}\right\|_{C_{t, x}^{N+1}} \sum_{M_{1}+M_{2}=1}\left\|\partial_{t}^{M_{1}} g_{(k)}\right\|_{L_{t}^{\infty}}\\& \sum_{0 \leq N_{1}+N_{2} \leq N}\left\|\partial_{t}^{M_{2}} \nabla^{N_{1}} \varphi_{\left(k_{1}\right)}^{2}\right\|_{L_{t}^{\infty} L_{x}^{2}}\left\|\nabla^{N_{2}} (\psi_{(k)}^{\prime})^{2}\right\|_{L_{t}^{\infty} L_{x}^{2}} \\
& \lesssim \sum_{M_{1}+M_{2}=1} \vartheta_{q+1}^{-5N-7} \mu^{-1}\lambda r_{\perp} \sigma^{M_{1}} \tau^{M_{1}+\frac{1}{2}} \lambda^{N} (r_{\perp} r_{\|})^{-\frac{1}{2}}\left(\frac{ \mu\tau^{\frac{1}{2}}}{r_{\|}}\right)^{M_{2}} \\
& \lesssim \vartheta_{q+1}^{-5N-7} \lambda^{N+1} \frac{r_{\perp}}{r_{\|}}(r_{\perp} r_{\|})^{-\frac{1}{2}} \tau . \label{4.51}\tag{3.51}
\end{align*}
and
\begin{align*}
\left\|\partial_{t} w_{q+1}^{(o)}\right\|_{L_{t}^{\infty} H_{x}^{N}} & \lesssim \sigma^{-1} \sum_{k \in \Lambda}\left\|\partial_{t}\left(h_{(k)} \nabla\left(a_{(k)}^{2}\right)\right)\right\|_{L_{t}^{\infty} H_{x}^{N}} \\
& \lesssim \sigma^{-1} \sum_{k \in \Lambda} \sum_{M_{1}+M_{2}=1}\left\|\partial_{t}^{M_{1}} h_{(k)}\right\|_{C_{t}}\left\|\partial_{t}^{M_{2}} \nabla\left(a_{(k)}^{2}\right)\right\|_{C_{t, x}^{N}} \\
& \lesssim \sigma^{-1} \sum_{M_{1}+M_{2}=1} \sigma^{M_{1}} \tau^{M_{1}} \vartheta_{q+1}^{-7-5\left(M_{2}+N\right)} \lesssim \vartheta_{q+1}^{-5N-7} \tau . \label{4.52}\tag{3.52}
\end{align*}

Hence, considering the fact that $\vartheta_{q+1}^{-5N-7}\ll\lambda^{\varepsilon}$
, we arrive at
$$
\begin{aligned}
&\left\|\partial_{t} w_{q+1}^{(p)}\right\|_{L_{t}^{\infty} H_{x}^{N}}+\left\|\partial_{t} w_{q+1}^{(c)}\right\|_{L_{t}^{\infty} H_{x}^{N}}+\left\|\partial_{t} w_{q+1}^{(t)}\right\|_{L_{t}^{\infty} H_{x}^{N}}+\left\|\partial_{t} w_{q+1}^{(o)}\right\|_{L_{t}^{\infty} H_{x}^{N}} \\
& \lesssim \vartheta_{q+1}^{-5 N-6} \lambda^{N} \frac{\mu \tau}{r_{\|}}+\vartheta_{q+1}^{-5 N-6} \lambda^{N}  \frac{\mu \tau}{r_{\|}} \frac{r_{\perp}}{r_{\|}}+\vartheta_{q+1}^{-5N-7} \lambda^{N+1} \frac{r_{\perp}}{r_{\|}}(r_{\perp} r_{\|})^{-\frac{1}{2}}+\vartheta_{q+1}^{-5N-7} \tau \\
& \lesssim \vartheta_{q+1}^{-5 N-6} \lambda^{N+4 \alpha-2+8 \varepsilon}+\vartheta_{q+1}^{-5N-7} \lambda^{N+4 \alpha-2+2 \varepsilon}+\vartheta_{q+1}^{-5N-7} \lambda^{4\alpha-4+16\varepsilon} \\
& \lesssim \lambda^{N+5}
\end{aligned}
$$
which verifies \eqref{4.47}. Consequently, the proof of Lemma \ref{D.6} is finished.

\subsection{Inductive estimates for velocity perturbations.}\ \ \ \
Based on the preceding estimations, we proceed to authenticate the inductive estimates \eqref{2.6}, \eqref{2.7}, and \eqref{2.11} - \eqref{2.13} pertaining to the velocity perturbations. Initially, to deduce the decay of the $L_{t, x}^{2}$-norms for these perturbations, we observe that $a_{(k)}$ possesses compact support within the domain $\mathbb{T}^{2}\times [0, T]$ that can be treated as a periodic function defined on $\mathbb{T}^{3}$. We invoke the $L^{p}$ decorrelation Lemma \ref{H.2}, substituting $f$ with $a_{(k)}$ and $g$ with $g_{(k)} \varphi_{\left(k_{1}\right)} \psi_{(k)}^{\prime}$ while setting $\sigma=\lambda^{2 \varepsilon}$. Subsequently, we utilize \eqref{2.2}-\eqref{2.4} alongside Lemmas \ref{D.3}, \ref{D.1}, and \ref{D.5} to obtain the desired results:
\begin{gather*}
\left\|w_{q+1}^{(p)}\right\|_{L_{t, x}^{2}} \lesssim \sum_{k \in \Lambda}\left(\left\|a_{(k)}\right\|_{L_{t, x}^{2}}\left\|g_{(k)}\right\|_{L_{t}^{2}}\left\|\varphi_{\left(k_{1}\right)} \psi_{(k)}^{\prime}\right\|_{C_{t} L_{x}^{2}}\right. \\
\left.+\sigma^{-\frac{1}{2}}\left\|a_{(k)}\right\|_{C_{t, x}^{1}}\left\|g_{(k)}\right\|_{L_{t}^{2}}\left\|\varphi_{\left(k_{1}\right)} \psi_{(k)}^{\prime}\right\|_{C_{t} L_{x}^{2}}\right) \\
\lesssim \delta_{q+1}^{\frac{1}{2}}+\vartheta_{q+1}^{-6} \lambda_{q+1}^{-\varepsilon} \lesssim \delta_{q+1}^{\frac{1}{2}} . \label{4.53}\tag{3.53}
\end{gather*}

In light of \eqref{2.3}, by employing \eqref{4.53} and invoking Lemma \ref{D.6}, we proceed to establish an upper bound for the velocity perturbation,
\begin{align*}
\left\|w_{q+1}\right\|_{L_{t, x}^{2}} & \lesssim\left\|w_{q+1}^{(p)}\right\|_{L_{t, x}^{2}}+\left\|w_{q+1}^{(c)}\right\|_{L_{t, x}^{2}}+\left\|w_{q+1}^{(t)}\right\|_{L_{t, x}^{2}}+\left\|w_{q+1}^{(o)}\right\|_{L_{t, x}^{2}} \\
& \lesssim \delta_{q+1}^{\frac{1}{2}}+\vartheta_{q+1}^{-6} r_{\perp} r_{\|}^{-1}+\vartheta_{q+1}^{-2} \mu^{-1} \lambda r_{\perp} (r_{\perp}r_{\|})^{-\frac{1}{2}}+\vartheta_{q+1}^{-7} \sigma^{-1} \lesssim \delta_{q+1}^{\frac{1}{2}}, \label{4.54}\tag{3.54}
\end{align*}
and
\begin{align*}
\left\|w_{q+1}\right\|_{L_{t}^{1} L_{x}^{2}} & \lesssim\left\|w_{q+1}^{(p)}\right\|_{L_{t}^{1} L_{x}^{2}}+\left\|w_{q+1}^{(c)}\right\|_{L_{t}^{1} L_{x}^{2}}+\left\|w_{q+1}^{(t)}\right\|_{L_{t}^{1} L_{x}^{2}}+\left\|w_{q+1}^{(o)}\right\|_{L_{t}^{1} L_{x}^{2}} \label{4.55}\tag{3.55}\\
& \lesssim \vartheta_{q+1}^{-1} \tau^{-\frac{1}{2}}+\vartheta_{q+1}^{-6} \frac{r_{\perp}}{r_{\|}} \tau^{-\frac{1}{2}}+\vartheta_{q+1}^{-2} \mu^{-1}\lambda r_{\perp} (r_{\perp}r_{\|})^{-\frac{1}{2}} \tau^{-\frac{1}{2}}+\vartheta_{q+1}^{-7} \sigma^{-1} \lesssim \lambda_{q+1}^{-\varepsilon}.
\end{align*}

Next, we proceed to verify the iterative estimates for $u_{q+1}$, we utilize \eqref{2.6}, \eqref{3.4}, \eqref{2.6}, \eqref{2.6}, \eqref{2.6} to derive the following results:
\begin{align*}
\left\|u_{q+1}\right\|_{L_{t}^{\infty} H_{x}^{3}} & \lesssim\left\|\widetilde{u}_{q}\right\|_{L_{t}^{\infty} H_{x}^{3}}+\left\|w_{q+1}\right\|_{L_{t}^{\infty} H_{x}^{3}} \\
& \lesssim \lambda_{q}^{5}+\lambda_{q+1}^{5} \lesssim \lambda_{q+1}^{5}, \label{4.56}\tag{3.56}
\end{align*}
and
\begin{align*}
\left\|\partial_{t} u_{q+1}\right\|_{L_{t}^{\infty} H_{x}^{2}} & \lesssim\left\|\partial_{t} \widetilde{u}_{q}\right\|_{L_{t}^{\infty} H_{x}^{2}}+\left\|\partial_{t} w_{q+1}\right\|_{L_{t}^{\infty} H_{x}^{2}} \\
& \lesssim \sup _{i}\left\|\partial_{t}\left(\mathscr{X} _{i} v_{i}\right)\right\|_{L_{t}^{\infty} H_{x}^{2}}+\lambda_{q+1}^{7} \\
& \lesssim \sup _{i}\left(\left\|\partial_{t} \mathscr{X} _{i}\right\|_{C_{t}}\left\|v_{i}\right\|_{L^{\infty}\left(\operatorname{supp} \mathscr{X} _{i} ; H_{x}^{2}\right)}+\left\|\mathscr{X} _{i}\right\|_{C_{t}}\left\|\partial_{t} v_{i}\right\|_{L^{\infty}\left(\operatorname{supp} \mathscr{X} _{i} ; H_{x}^{2}\right)}\right)+\lambda_{q+1}^{7} \\
& \lesssim \vartheta_{q+1}^{-1} \lambda_{q}^{5}+\ell_{q+1} \lambda_{q}^{5}+\lambda_{q+1}^{7} \lesssim \lambda_{q+1}^{7}. \label{4.57}\tag{3.57}
\end{align*}

Furthermore, given the fact $\lambda_{q}^{-3} \ll \delta_{q+2}^{1 / 2}$, we employ \eqref{2.2}, \eqref{2.3}, \eqref{3.5} and \eqref{4.54} to deduce that
\begin{align*}
\left\|u_{q}-u_{q+1}\right\|_{L_{t, x}^{2}} & \leq\left\|u_{q}-\widetilde{u}_{q}\right\|_{L_{t, x}^{2}}+\left\|\widetilde{u}_{q}-u_{q+1}\right\|_{L_{t, x}^{2}} \\
& \lesssim\left\|u_{q}-\widetilde{u}_{q}\right\|_{L_{t}^{\infty} L_{x}^{2}}+\left\|w_{q+1}\right\|_{L_{t, x}^{2}} \\
& \lesssim \lambda_{q}^{-3}+\delta_{q+1}^{\frac{1}{2}} \leq M^{*} \delta_{q+1}^{\frac{1}{2}}, \label{4.58}\tag{3.58}
\end{align*}
for $M^{*}$ large enough and
\begin{align*}
\left\|u_{q}-u_{q+1}\right\|_{L_{t}^{1} L_{x}^{2}} & \lesssim\left\|u_{q}-\widetilde{u}_{q}\right\|_{L_{t}^{\infty} L_{x}^{2}}+\left\|w_{q+1}\right\|_{L_{t}^{1} L_{x}^{2}} \\
& \lesssim \lambda_{q}^{-3}+\lambda_{q+1}^{-\varepsilon} \leq \delta_{q+2}^{\frac{1}{2}}, \label{4.59}\tag{3.59}
\end{align*}
where we set $a$ to be a sufficiently large value, ensuring the validity of the final inequalities presented in \eqref{4.59}.

Considering the iteration estimate \eqref{2.13}, we assert that the Sobolev embedding
\begin{equation*}
H_{x}^{3} \hookrightarrow W_{x}^{s, p} \label{4.60}\tag{3.60}
\end{equation*}
holds true for any triplet $(s, p, \gamma)\in \mathcal{A}_{1}$. Leveraging the results established in \eqref{2.9}, \eqref{3.16}, \eqref{4.60}, we can deduce that
\begin{align*}
\left\|\widetilde{u}_{q}-u_{q}\right\|_{L_{t}^{\gamma} W_{x}^{s, p}} & \lesssim\left\|\sum_{i} \mathscr{X} _{i}\left(v_{i}-u_{q}\right)\right\|_{L_{t}^{\infty} H_{x}^{3}} \\
& \lesssim \sup _{i}\Big(\left\|\mathscr{X} _{i}\left(v_{i}-u_{q}\right)\right\|_{L^{\infty}\left(\left(\operatorname{supp}\left(\mathscr{X} _{i}\right) ; H_{x}^{3}\right)\right.}\\&\quad+\left\|\left(1-\mathscr{X} _{i}\right)\left(v_{i-1}-u_{q}\right)\right\|_{L^{\infty}\left(\left(\operatorname{supp}\left(\mathscr{X} _{i} \mathscr{X} _{i-1}\right) ; H_{x}^{3}\right)\right.}\Big) \\
& \lesssim \sup _{i}\left|t_{i+1}+\vartheta_{q+1}-t_{i}\right|\left\||\nabla| \mathring{\mathsf{R}}_{q}\right\|_{L_{t}^{\infty} H_{x}^{3}} \lesssim \ell_{q+1}^{-1} \lambda_{q}^{10} \lesssim \lambda_{q}^{-2} . \label{4.61}\tag{3.61}
\end{align*}
Then, for any $\alpha \in[1 , \frac{3}{2})$, by \eqref{2.3}, \eqref{4.1} and Lemma \ref{D.6},
\begin{align*}
\left\|u_{q+1}-u_{q}\right\|_{L_{t}^{\gamma} W_{x}^{s, p}} & \lesssim\left\|\widetilde{u}_{q}-u_{q}\right\|_{L_{t}^{\gamma} W_{x}^{s, p}}+\left\|w_{q+1}\right\|_{L_{t}^{\gamma} W_{x}^{s, p}} \\
& \lesssim \lambda_{q}^{-2}+\vartheta_{q+1}^{-1} \lambda_{q+1}^{s} (r_{\perp} r_{\|})^{\frac{1}{p}-\frac{1}{2}} \tau^{\frac{1}{2}-\frac{1}{\gamma}}+\vartheta_{q+1}^{-17} \sigma^{-1} \\
& \lesssim \lambda_{q}^{-2}+\lambda_{q+1}^{s+2 \alpha-1-\frac{2}{p}-\frac{4 \alpha-4}{\gamma}+\varepsilon\left(2+\frac{12}{p}-\frac{16}{\gamma}\right)}+\lambda_{q+1}^{-\varepsilon}. \label{4.62}\tag{3.62}
\end{align*}
Incorporating the considerations from \eqref{2.4}, we have:
\begin{equation*}
s+2 \alpha-1-\frac{2}{p}-\frac{4 \alpha-4}{\gamma}+\varepsilon\left(2+\frac{12}{p}-\frac{16}{\gamma}\right) \leq s+2 \alpha-1-\frac{2}{p}-\frac{4 \alpha-4}{\gamma}+14 \varepsilon<-6 \varepsilon,\label{4.63} \tag{3.63}
\end{equation*}
which yields that
\begin{equation*}
\left\|u_{q+1}-u_{q}\right\|_{L_{t}^{\gamma} W_{x}^{s, p}} \leq \delta_{q+2}^{\frac{1}{2}}. \label{4.64}\tag{3.64}
\end{equation*}

Hence, the iteration estimates \eqref{2.6}, \eqref{2.7}, and \eqref{2.11}- \eqref{2.13} are verified.

{\centering
\section{Reynolds stress for the super-critical regime $\mathcal{A}_{1}$}\label{D}}
The objective of this section is to verify the inductive estimates \eqref{2.8}- \eqref{2.10} for the new Reynolds stress $\mathring{\mathsf{R}}_{q+1}$ within the super-critical regime $\mathcal{A}_{1}$ when $\alpha \in[1,\frac{3}{2})$, whose borderline includes one endpoint case characterized by the triplet $(s, \gamma, p)=(s,\infty, \frac{2}{2\alpha-1+s})$ in view of gLPS condition \eqref{1.5}.\\
\subsection{Decomposition of Reynolds stress.}\ \ \ \
From \eqref{2.1} and \eqref{4.41}, we deduce that the new Reynolds stress $\mathring{\mathsf{R}}_{q+1}$ is given as follows:
\begin{align*}
& \operatorname{div} \mathring{\mathsf{R}}_{q+1}-\nabla(\mathsf{P}_{q+1}-\mathsf{P}_{q})=\underbrace{\partial_{t}\left(w_{q+1}^{(p)}+w_{q+1}^{(c)}\right)+(-\Delta)^{\alpha} w_{q+1}+\operatorname{div}\left(\widetilde{u}_{q} \otimes w_{q+1}+w_{q+1} \otimes \widetilde{u}_{q}\right)}_{\operatorname{div} \mathring{\mathsf{R}}_{\{l i n\}}+\nabla\mathsf{P}_{\{l i n\}}} \\
& +\underbrace{\operatorname{div}\left(w_{q+1}^{(p)} \otimes w_{q+1}^{(p)}+\mathring{\widetilde{\mathsf{R}}}_{q}\right)+\partial_{t} w_{q+1}^{(t)}+\partial_{t} w_{q+1}^{(o)}}_{\operatorname{div} \mathring{\mathsf{R}}_{\{osc\}}+\nabla\mathsf{P}_{\{o s c\}}}\label{5.1}\tag{4.1} \\
& +\underbrace{\operatorname{div}\left(\left(w_{q+1}^{(c)}+w_{q+1}^{(t)}+w_{q+1}^{(o)}\right) \otimes w_{q+1}+w_{q+1}^{(p)} \otimes\left(w_{q+1}^{(c)}+w_{q+1}^{(t)}+w_{q+1}^{(o)}\right)\right)}_{\operatorname{div} \mathring{\mathsf{R}}_{\{cor\}}+\nabla\mathsf{P}_{\{c o r\}}}.
\end{align*}

By employing the inverse divergence operator designated as $\mathcal{R}$, the definition of Reynolds stress at the level $q+1$ becomes:
\begin{equation*}
\mathring{\mathsf{R}}_{q+1}:=\mathring{\mathsf{R}}_{\{l i n\}}+\mathring{\mathsf{R}}_{\{osc\}}+\mathring{\mathsf{R}}_{\{cor\}}, \label{5.2}\tag{4.2}
\end{equation*}
where the linear error
\begin{equation*}
\mathring{\mathsf{R}}_{ \{lin \}}:=\mathcal{R}\left(\partial_{t}\left(w_{q+1}^{(p)}+w_{q+1}^{(c)}\right)\right)+\mathcal{R}(-\Delta)^{\alpha} w_{q+1}+\mathcal{R} \mathbb{P}_{H} \operatorname{div}\left(\widetilde{u}_{q} \mathring{\otimes} w_{q+1}+w_{q+1} \mathring{\otimes} \tilde{u}_{q}\right) , \label{5.3}\tag{4.3}
\end{equation*}
the oscillation error
\begin{align*}
\mathring{\mathsf{R}}_{\{osc\}}:= & \sum_{k \in \Lambda} \mathcal{R} \mathbb{P}_{H} \mathbb{P}_{\neq 0}\left(g_{(k)}^{2} \mathbb{P}_{\neq 0}\left(\mathbf{W}_{k} \otimes \mathbf{W}_{k}\right) \nabla\left(a_{(k)}^{2}\right)\right) \\
& -\mu^{-1} \sum_{k \in \Lambda} \mathcal{R} \mathbb{P}_{H} \mathbb{P}_{\neq 0}\left(\partial_{t}\left(a_{(k)}^{2} g_{(k)}\right) \varphi_{\left(k_{1}\right)}^{2} (\psi_{(k)}^{\prime})^{2} k_{1}\right) \\
& -\sigma^{-1} \sum_{k \in \Lambda} \mathcal{R} \mathbb{P}_{H} \mathbb{P}_{\neq 0}\left(h_{(k)} \fint_{\mathbb{T}^{2}} \mathbf{W}_{k} \otimes \mathbf{W}_{k} \mathrm{d} x \partial_{t} \nabla\left(a_{(k)}^{2}\right)\right) , \label{5.4}\tag{4.4}
\end{align*}
and the corrector error
\begin{equation*}
\mathring{\mathsf{R}}_{\{cor\}}:=\mathcal{R} \mathbb{P}_{H} \operatorname{div}\left(w_{q+1}^{(p)} \mathring{\otimes}\left(w_{q+1}^{(c)}+w_{q+1}^{(t)}+w_{q+1}^{(o)}\right)+\left(w_{q+1}^{(c)}+w_{q+1}^{(t)}+w_{q+1}^{(o)}\right) \mathring{\otimes} w_{q+1}\right) . \label{5.5}\tag{4.5}
\end{equation*}

Furthermore, it is worthy noting that, as exemplified in references \cite{TBV19,YLZ22}, we have the following mathematical expression:
\begin{equation*}
\mathring{\mathsf{R}}_{q+1}=\mathcal{R} \mathbb{P}_{H} \operatorname{div} \mathring{\mathsf{R}}_{q+1} . \label{5.6}\tag{4.6}
\end{equation*}

\subsection{$L_t^{\infty} H_x^N$-estimates of Reynolds stress.}\ \ \ \
Concerning the Reynolds stress estimates expressed in the $L_{t}^{\infty} H_{x}^{N}$-norm, specifically \eqref{2.8} and \eqref{2.9}. By the identity \eqref{5.6} and \eqref{2.1} for $\left(u_{q+1}, \mathring{\mathsf{R}}_{q+1}\right)$, we obtain that
\begin{align*}
\left\|\mathring{\mathsf{R}}_{q+1}\right\|_{L_{t}^{\infty} H_{x}^{N}} & \leq\left\|\mathcal{R} \mathbb{P}_{H}\left(\operatorname{div} \mathring{\mathsf{R}}_{q+1}\right)\right\|_{L_{t}^{\infty} H_{x}^{N}} \\
& \lesssim\left\|\partial_{t} u_{q+1}+\operatorname{div}\left(u_{q+1} \otimes u_{q+1}\right)+(-\Delta)^{\alpha} u_{q+1}\right\|_{L_{t}^{\infty} H_{x}^{N-1}} \\
& \lesssim\left\|\partial_{t} u_{q+1}\right\|_{L_{t}^{\infty} H_{x}^{N-1}}+\left\|u_{q+1} \otimes u_{q+1}\right\|_{L_{t}^{\infty} H_{x}^{N}}+\left\|u_{q+1}\right\|_{L_{t}^{\infty} H_{x}^{N+3}} \\
& \lesssim\left\|\partial_{t} u_{q+1}\right\|_{L_{t}^{\infty} H_{x}^{N-1}}+\left\|u_{q+1}\right\|_{L_{t}^{\infty} H_{x}^{N}}\left\|u_{q+1}\right\|_{L_{t, x}^{\infty}}+\left\|u_{q+1}\right\|_{L_{t}^{\infty} H_{x}^{N+3}} \label{5.7}\tag{4.7}
\end{align*}

Recalling that
$u_{q+1}:=\widetilde{u}_{q}+w_{q+1}$, by \eqref{3.1}, \eqref{3.4}, \eqref{3.10}, \eqref{4.41} and Lemma \ref{D.6}, for any integers $0 \leq \widetilde{N} \leq 4$, we have,
\begin{equation*}
\left\|u_{q+1}\right\|_{L_{t}^{\infty} H_{x}^{\widetilde{N}+3}} \lesssim \lambda_{q+1}^{\widetilde{N}+5}, \quad\left\|\partial_{t} u_{q+1}\right\|_{L_{t}^{\infty} H_{x}^{\widetilde{N}-1}} \lesssim \lambda_{q+1}^{\widetilde{N}+4}, \quad\left\|u_{q+1}\right\|_{L_{t, x}^{\infty}}  \lesssim \lambda_{q+1}^{4},
\label{5.8}\tag{4.8}
\end{equation*}
which yields that
$$
\left\|\mathring{\mathsf{R}}_{q+1}\right\|_{L_{t}^{\infty} H_{x}^{N}} \lesssim \lambda_{q+1}^{N+4}+\lambda_{q+1}^{N+6}+\lambda_{q+1}^{N+5} \lesssim \lambda_{q+1}^{N+6}.
$$
where the implicit constants are independent of $q$. Hence, the verification of the $L_{t}^{\infty} H_{x}^{N}$-estimates \eqref{2.8} and \eqref{2.9} for $\mathring{\mathsf{R}}_{q+1}$ have been established.
\subsection{$L_{t, x}^1$-decay of Reynolds stress.}\ \ \ \
In the following section, we focus on verifying the precise $L_{t, x}^{1}$-decay \eqref{2.10} exhibited by the Reynolds stress $\mathring{\mathsf{R}}_{q+1}$ at level $q+1$. Given the intricate nature of the spatial and temporal oscillations inherent in the ``building blocks'' constructed in Section \ref{C}, it is imperative to exercise caution when estimating the time/space derivatives within the linear error term $\mathring{\mathsf{R}}_{\{lin\}}$ and the high-frequency components of the oscillation error $\mathring{\mathsf{R}}_{\{osc\}}$.

As Calder\'{o}n-Zygmund operators are bounded in the space $L_{x}^{\varrho}$, where $1 < \varrho < 2$, we give a specific choice:
\begin{equation*}
\varrho:=\frac{2-12 \varepsilon}{2-13 \varepsilon} \in(1,2), \label{5.12}\tag{4.9}
\end{equation*}
where $\varepsilon$ is defined in \eqref{2.4}. Note that,
\begin{equation*}
(2-12 \varepsilon)\left(1-\frac{1}{\varrho}\right)=\varepsilon \label{5.13}\tag{4.10}
\end{equation*}
and
\begin{equation*}
(r_{\perp}r_{\|})^{\frac{1}{\varrho}-1}=\lambda^{\varepsilon}, \quad (r_{\perp}r_{\|})^{\frac{1}{\varrho}-\frac{1}{2}}=\lambda^{-1+7\varepsilon}. \label{5.14}\tag{4.11}
\end{equation*}

\noindent {\bf (i) Linear Error.} For the acceleration part of the linear error, taking time derivative of \eqref{4.34} and using identity \eqref{4.9}, we obtain
$$
\begin{aligned}
\partial_{t}\Big(w&_{q+1}^{(p)}+w_{q+1}^{(c)}\Big)  =(\lambda r_{\perp}N_{\Lambda})^{-1} \sum_{k} \nabla^{\perp}\left[\partial_{t}\left(a_{(k)} g_{(k)}\right) \boldsymbol{\Psi}_{k}\right]+(\lambda r_{\perp}N_{\Lambda})^{-1} \sum_{k} \nabla^{\perp}\left[a_{(k)} g_{(k)} \partial_{t} \boldsymbol{\Psi}_{k}\right] \\
& =(\lambda r_{\perp}N_{\Lambda})^{-1} \sum_{k} \nabla^{\perp}\left[\partial_{t}\left(a_{(k)} g_{(k)}\right) \boldsymbol{\Psi}_{k}\right]+(\lambda r_{\perp}N_{\Lambda})^{-2} \mu N_{\Lambda}\sum_{k} \nabla^{\perp}\left[a_{(k)} g_{(k)}^{2}\left(k_{1} \cdot \nabla\right) \boldsymbol{\Psi}_{k}\right].
\end{aligned}
$$
Note that $\mathcal{R}
\nabla^{\perp}$ is a Calder\'{o}n-Zygmund operator on $\mathbb{T}^2$, we can employ Lemmas \ref{D.1}, \ref{D.3}, and \ref{D.5} to estimate the first term:
$$
\begin{aligned}
(\lambda r_{\perp}N_{\Lambda}&)^{-1}\left\|\sum_{k} \mathcal{R} \nabla^{\perp}\left[\partial_{t}\left(a_{(k)} g_{(k)}\right) \boldsymbol{\Psi}_{k}\right]\right\|_{L^{1}\left(0,1 ; L_{x}^{\varrho}\right)}\\
 & \lesssim (\lambda r_{\perp})^{-1} \sum_{k}\left\|a_{(k)}\right\|_{C_{x, t}^{1}}\left\|g_{(k)}\right\|_{W^{1,1}}\left\|\boldsymbol{\Psi}_{k}\right\|_{L_{t}^{\infty} L_{x}^{\varrho}} \\
& \lesssim \vartheta_{q+1}^{-6}(\lambda r_{\perp})^{-1}\sigma \tau^{\frac{1}{2}} r_{\perp}(r_{\perp} r_{\|})^{\frac{1}{\varrho}-\frac{1}{2}}.
\end{aligned}
$$

As for the second term, we first estimate the derivative of $\boldsymbol{\Psi}_k$ in the direction $k_1$ as
$$
\left\|\left(k_{1} \cdot \nabla\right) \boldsymbol{\Psi}_{k}\right\|_{L_{t}^{\infty} L^{\varrho}_{x}} \lesssim \lambda r_{\perp} \frac{r_{\perp}}{r_{\|}}(r_{\perp}r_{\|})^{\frac{1}{\varrho}-\frac{1}{2}} .
$$
where the order of derivative is $\lambda\frac{r_{\perp}}{r_{\|}}$ (rather than $\lambda$ for the full gradient), together with Lemmas \ref{4.3}, \ref{4.5} imply that
$$
\begin{aligned}
(\lambda r_{\perp}N_{\Lambda})^{-2} \mu N_{\Lambda}&\left\|\sum_{k} \mathcal{R} \nabla^{\perp}\left[a_{(k)} g_{(k)}^{2}\left(k_{1} \cdot \nabla\right) \boldsymbol{\Psi}_{k}\right]\right\|_{L^{1}_{t}L_{x}^{\varrho}} \\
& \lesssim (\lambda r_{\perp})^{-2} \mu \sum_{k}\left\|a_{(k)}\right\|_{L_{x, t}^{\infty}}\left\|g_{(k)}^{2}\right\|_{L^{1}}\left\|\left(k_{1} \cdot \nabla\right) \boldsymbol{\Psi}_{k}\right\|_{L_{t}^{\infty} L_{x}^{\varrho}} \\
& \lesssim\vartheta_{q+1}^{-1} (\lambda r_{\|})^{-1} \mu(r_{\perp} r_{\|})^{\frac{1}{\varrho}-\frac{1}{2}} .
\end{aligned}
$$
Due to \eqref{2.4} and \eqref{4.1}, combining both terms we obtain
\begin{align*}
\left\| \mathcal{R}\partial_{t}\left(w_{q+1}^{(p)}+w_{q+1}^{(c)}\right)
\right\|_{L^{1}_{t}L_{x}^{\varrho}} & \lesssim \vartheta_{q+1}^{-6}\left(\tau^{\frac{1}{2}} \sigma\lambda^{-1}(r_{\perp} r_{\|})^{\frac{1}{\varrho}-\frac{1}{2}}+(\lambda r_{\|})^{-1}\mu(r_{\perp} r_{\|})^{\frac{1}{\varrho}-\frac{1}{2}}\right) \\
& \lesssim \vartheta_{q+1}^{-6} \lambda^{-\varepsilon}.\label{5.15}\tag{4.12}
\end{align*}

Considering the viscosity term $(-\Delta)^\alpha w_{q+1}$ and using \eqref{4.40}, we have
\begin{align*}
\left\|\mathcal{R}(-\Delta)^{\alpha} w_{q+1}\right\|_{L_{t}^{1} L_{x}^{\varrho}} \lesssim & \left\|\mathcal{R}(-\Delta)^{\alpha} w_{q+1}^{(p)}\right\|_{L_{t}^{1} L_{x}^{\varrho}}+\left\|\mathcal{R}(-\Delta)^{\alpha} w_{q+1}^{(c)}\right\|_{L_{t}^{1} L_{x}^{\varrho}} \\
& +\left\|\mathcal{R}(-\Delta)^{\alpha} w_{q+1}^{(t)}\right\|_{L_{t}^{1} L_{x}^{\varrho}}+\left\|\mathcal{R}(-\Delta)^{\alpha} w_{q+1}^{(o)}\right\|_{L_{t}^{1} L_{x}^{\varrho}}. \label{5.16}\tag{4.13}
\end{align*}
To estimate the right-hand side mentioned above, we employ the interpolation inequality (see \cite{HB18}), specifically Lemma \ref{D.6}, along with the condition that $3-2\alpha \geq 20 \varepsilon$ to deduce
\begin{align*}
\left\|\mathcal{R}(-\Delta)^{\alpha} w_{q+1}^{(p)}\right\|_{L_{t}^{1} L_{x}^{\varrho}} & \lesssim\left\||\nabla|^{2 \alpha-1} w_{q+1}^{(p)}\right\|_{L_{t}^{1} L_{x}^{\varrho}} \\
& \lesssim\left\|w_{q+1}^{(p)}\right\|_{L_{t}^{1} L_{x}^{\varrho}}^{\frac{4-2 \alpha}{3}}\left\|w_{q+1}^{(p)}\right\|_{L_{t}^{1} W_{x}^{3, \varrho}}^{\frac{2 \alpha-1}{3}} \\
& \lesssim \vartheta_{q+1}^{-1} \lambda^{2 \alpha-1} (r_{\perp} r_{\|})^{\frac{1}{\varrho}-\frac{1}{2}} \tau^{-\frac{1}{2}} \lesssim \vartheta_{q+1}^{-1} \lambda^{-\varepsilon} , \label{5.17}\tag{4.14}
\end{align*}
and
\begin{align*}
& \left\|\mathcal{R}(-\Delta)^{\alpha} w_{q+1}^{(c)}\right\|_{L_{t}^{1} L_{x}^{\varrho}} \lesssim \vartheta_{q+1}^{-6} \lambda^{2 \alpha-1} (r_{\perp} r_{\|})^{\frac{1}{\varrho}-\frac{1}{2}} \frac{r_{\perp}}{r_{\|}} \tau^{-\frac{1}{2}} \lesssim \vartheta_{q+1}^{-6} \lambda^{-9 \varepsilon},  \label{5.18}\tag{4.15}\\
& \left\|\mathcal{R}(-\Delta)^{\alpha} w_{q+1}^{(t)}\right\|_{L_{t}^{1} L_{x}^{\varrho}} \lesssim \vartheta_{q+1}^{-2} \lambda^{2 \alpha-1} \tau^{-\frac{1}{2}}\mu^{-1} (r_{\perp} r_{\|})^{\frac{1}{\varrho}-1} \lambda r_{\perp}\lesssim \vartheta_{q+1}^{-2} \lambda^{-7\varepsilon}, \label{5.19} \tag{4.16}\\
& \left\|\mathcal{R}(-\Delta)^{\alpha} w_{q+1}^{(o)}\right\|_{L_{t}^{1} L_{x}^{\varrho}} \lesssim \vartheta_{q+1}^{-17} \sigma^{-1} \lesssim \vartheta_{q+1}^{-17} \lambda^{-2 \varepsilon}. \label{5.20}\tag{4.17}
\end{align*}

Therefore, by combining \eqref{5.17}-\eqref{5.20} with the observation that $\vartheta_{q+1}^{-17}\ll\lambda^{\varepsilon}$, we derive that
\begin{equation*}
\left\|\mathcal{R}(-\Delta)^{\alpha} w_{q+1}\right\|_{L_{t}^{1} L_{x}^{\varrho}} \lesssim \vartheta_{q+1}^{-1} \lambda^{-\varepsilon}. \label{5.21}\tag{4.18}
\end{equation*}

The nonlinear part of \eqref{5.3} remains unestimated. Utilizing the Sobolev embedding theorem $H_{x}^{3} \hookrightarrow L_{x}^{\infty}$, \eqref{3.4}, along with Lemma \ref{D.6},
\begin{align*}
& \left\|\mathcal{R} \mathbb{P}_{H} \operatorname{div}\left(w_{q+1} \otimes \widetilde{u}_{q}+\widetilde{u}_{q} \otimes w_{q+1}\right)\right\|_{L_{t}^{1} L_{x}^{\varrho}} \\
\lesssim & \left\|w_{q+1} \otimes \widetilde{u}_{q}+\widetilde{u}_{q} \otimes w_{q+1}\right\|_{L_{t}^{1} L_{x}^{\varrho}} \\
\lesssim & \left\|\widetilde{u}_{q}\right\|_{L_{t}^{\infty} H_{x}^{3}}\left\|w_{q+1}\right\|_{L_{t}^{1} L_{x}^{\varrho}} \\
\lesssim & \lambda_{q}^{5}\left(\vartheta_{q+1}^{-1} (r_{\perp} r_{\|})^{\frac{1}{\varrho}-\frac{1}{2}} \tau^{-\frac{1}{2}}+\vartheta_{q+1}^{-2} \tau^{-\frac{1}{2}}\mu^{-1} (r_{\perp}r_{\|})^{\frac{1}{\varrho}-1}\lambda r_{\perp}+\vartheta_{q+1}^{-7} \sigma^{-1}\right) \lesssim \vartheta_{q+1}^{-8} \lambda^{-2 \varepsilon} . \label{5.22}\tag{4.19}
\end{align*}

Based on the analysis presented in \eqref{5.15}, \eqref{5.21}, and \eqref{5.22}, we arrive at
\begin{equation*}
\left\|\mathring{\mathsf{R}}_{\{l i n\}}\right\|_{L_{t}^{1} L_{x}^{\varrho}} \lesssim \vartheta_{q+1}^{-6} \lambda^{-\varepsilon}+\vartheta_{q+1}^{-1} \lambda^{-\varepsilon}+\vartheta_{q+1}^{-8} \lambda^{-2 \varepsilon} \lesssim \vartheta_{q+1}^{-8} \lambda^{-\varepsilon} . \label{5.23}\tag{4.20}
\end{equation*}

\noindent {\bf (ii) Oscillation Error.} Now we are prepared to estimate the oscillation error, we  decompose it into three distinct components as:
$$
\mathring{\mathsf{R}}_{\{osc\}}=\mathring{\mathsf{R}}_{\{o s c .1\}}+\mathring{\mathsf{R}}_{\{o s c .2\}}+\mathring{\mathsf{R}}_{\{o s c .3\}},
$$
where the low-high spatial oscillation error
$$
\mathring{\mathsf{R}}_{\{o s c .1\}}:=\sum_{k \in \Lambda} \mathcal{R} \mathbb{P}_{H} \mathbb{P}_{\neq 0}\left(g_{(k)}^{2} \mathbb{P}_{\neq 0}\left(\mathbf{W}_{k} \otimes \mathbf{W}_{k}\right) \nabla\left(a_{(k)}^{2}\right)\right),
$$
the high temporal oscillation error
$$
\mathring{\mathsf{R}}_{\{o s c .2\}}:=-\mu^{-1}\lambda r_{\perp} \sum_{k \in \Lambda} \mathcal{R} \mathbb{P}_{H} \mathbb{P}_{\neq 0}\left(\partial_{t}\left(a_{(k)}^{2} g_{(k)}\right) \varphi_{\left(k_{1}\right)} ^{2} (\psi_{(k)}^{\prime})^{2} k_{1}\right),
$$
and the low frequency error
$$
\mathring{\mathsf{R}}_{\{o s c .3\}}:=-\sigma^{-1} \sum_{k \in \Lambda} \mathcal{R} \mathbb{P}_{H} \mathbb{P}_{\neq 0}\left(h_{(k)} \fint_{\mathbb{T}^{2}} \mathbf{W}_{k} \otimes \mathbf{W}_{k} \mathrm{d} x \partial_{t} \nabla\left(a_{(k)}^{2}\right)\right).
$$

In terms of the low-high spatial oscillation error $\mathring{\mathsf{R}}_{\text {osc.1 }}$, note that the velocity flows exhibit high oscillations
$$
\mathbb{P}_{\neq 0}\left(\mathbf{W}_{k} \otimes \mathbf{W}_{k}\right)=\mathbb{P}_{\geq\left(\lambda r_{\perp} / 2\right)}\left(\mathbf{W}_{k} \otimes \mathbf{W}_{k}\right).
$$
Hence, we employ Lemmas \ref{D.1} and \ref{D.5}, and utilize Lemma \ref{H.3} specifically with the substitutions $a=
\nabla \left(a_{(k)}^{2}\right)$ and $f=\varphi_{\left(k_{1}\right)}^{2} (\psi_{(k)}^{\prime})^{2}$ to derive
\begin{align*}
\left\|\mathring{\mathsf{R}}_{\{o s c .1\}}\right\|_{L_{t}^{1} L_{x}^{\varrho}} & \lesssim \sum_{k \in \Lambda}\left\|g_{(k)}\right\|_{L_{t}^{2}}^{2}\left\||\nabla|^{-1} \mathbb{P}_{\neq 0}\left(\mathbb{P}_{\geq\left(\lambda r_{\perp} / 2\right)}\left(\mathbf{W}_{k} \otimes \mathbf{W}_{k}\right) \nabla\left(a_{(k)}^{2}\right)\right)\right\|_{C_{t} L_{x}^{\varrho}} \\
& \lesssim \sum_{k \in \Lambda}\left\||\nabla|^{3}\left(a_{(k)}^{2}\right)\right\|_{C_{t, x}} (\lambda r_{\perp})^{-1}\left\|\varphi_{\left(k_{1}\right)} \right\|^{2}_{C_{t} L_{x}^{2\varrho}}\left\|\psi_{(k)}^{\prime}\right\|^{2}_{C_{t} L_{x}^{2\varrho}} \\
& \lesssim \vartheta_{q+1}^{-17} \lambda^{-1} r_{\perp}^{\frac{1}{\varrho}-2} r_{\|}^{\frac{1}{\varrho}-1}. \label{5.24}\tag{4.21}
\end{align*}

Furthermore, we utilize Lemmas \ref{D.3}, \ref{D.1}, \ref{D.5} and introduce the large temporal oscillation parameter $\mu$ to balance the high temporal oscillation error $\mathring{\mathsf{R}}_{\{osc. 2\}}$:
\begin{align*}
\left\|\mathring{\mathsf{R}}_{\{o s c .2\}}\right\|_{L_{t}^{1} L_{x}^{\varrho}} & \lesssim \mu^{-1} \lambda r_{\perp}\sum_{k \in \Lambda}\left\|\mathcal{R}\mathbb{P}_{H} \mathbb{P}_{\neq 0}\left(\partial_{t}\left(a_{(k)}^{2} g_{(k)}\right) \varphi_{\left(k_{1}\right)} ^{2} (\psi_{(k)}^{\prime})^{2} k_{1}\right)\right\|_{L_{t}^{1} L_{x}^{\varrho}} \\
& \lesssim \mu^{-1} \lambda r_{\perp}\sum_{k \in \Lambda}\left(\left\|\partial_{t}\left(a_{(k)}^{2}\right)\right\|_{C_{t, x}}\left\|g_{(k)}\right\|_{L_{t}^{1}}+\left\|a_{(k)}\right\|_{C_{t, x}}^{2}\left\|\partial_{t} g_{(k)}\right\|_{L_{t}^{1}}\right)\\&\quad
\left\|\varphi_{\left(k_{1}\right)} \right\|_{C_{t} L_{x}^{2 \varrho}}^{2}\left\|\psi_{(k)}^{\prime}\right\|_{L_{x}^{2 \varrho}}^{2} \\
& \lesssim\left(\vartheta_{q+1}^{-7}\tau^{-\frac{1}{2}}+\vartheta_{q+1}^{-2} \tau^{\frac{1}{2}} \sigma\right) \mu^{-1}\lambda r_{\perp}^{\frac{1}{\varrho}} r_{\|}^{\frac{1}{\varrho}-1} \\
& \lesssim \vartheta_{q+1}^{-2} \tau^{\frac{1}{2}} \sigma \mu^{-1}\lambda r_{\perp}^{\frac{1}{\varrho}} r_{\|}^{\frac{1}{\varrho}-1} . \label{5.25}\tag{4.22}
\end{align*}

The estimates of the low-frequency error $\mathring{\mathsf{R}}_{\{osc. 3\}}$ can be estimated through the utilization of \eqref{4.18} and
Lemma \ref{D.5}:
\begin{align*}
\left\|\mathring{\mathsf{R}}_{\{o s c .3\}}\right\|_{L_{t}^{1} L_{x}^{\varrho}} &\lesssim \sigma^{-1} \sum_{k \in \Lambda}\left\|h_{(k)} \partial_{t} \nabla\left(a_{(k)}^{2}\right)\right\|_{L_{t}^{1} L_{x}^{\varrho}}\\
& \lesssim \sigma^{-1} \sum_{k \in \Lambda}\left\|h_{(k)}\right\|_{C_{t}}\left(\left\|a_{(k)}\right\|_{C_{t, x}}\left\|a_{(k)}\right\|_{C_{t, x}^{2}}+\left\|a_{(k)}\right\|_{C_{t, x}^{1}}^{2}\right)\\
& \lesssim \vartheta_{q+1}^{-12} \sigma^{-1} . \label{5.26}\tag{4.23}
\end{align*}

Therefore, combining \eqref{5.24}- \eqref{5.26} and utilizing expressions \eqref{4.1}, \eqref{5.14}, along with the constraint $0 < \varepsilon < (3-2\alpha) / 20$, we arrive at
\begin{align*}
\left\|\mathring{\mathsf{R}}_{\{osc\}}\right\|_{L_{t}^{1} L_{x}^{\varrho}} & \lesssim \vartheta_{q+1}^{-17} \lambda^{-1} r_{\perp}^{\frac{1}{\varrho}-3} r_{\|}^{\frac{1}{\varrho}-1}+\vartheta_{q+1}^{-2} \tau^{\frac{1}{2}} \sigma \mu^{-1} \lambda r_{\perp}^{\frac{1}{\varrho}} r_{\|}^{\frac{1}{\varrho}-1}+\vartheta_{q+1}^{-12} \sigma^{-1} \\
& \lesssim \vartheta_{q+1}^{-17} \lambda^{-\varepsilon}+\vartheta_{q+1}^{-2} \lambda^{2 \alpha-3+13 \varepsilon}+\vartheta_{q+1}^{-12} \lambda^{-2 \varepsilon} \\
& \lesssim \vartheta_{q+1}^{-17} \lambda^{-\varepsilon} . \label{5.27}\tag{4.24}
\end{align*}

\noindent {\bf (iii) Corrector Error.} Drawing inspiration from \cite{YLZ22}, we introduce $p_{1}, p_{2} \in(1, \infty)$ satisfying the conditions
$$
\frac{1}{p_{1}}=1-\widetilde{\eta}, \quad \frac{1}{p_{1}}=\frac{1}{p_{2}}+\frac{1}{2},
$$
with $\widetilde{\eta} \leq 5\varepsilon / (2-12 \varepsilon)$. by H\"{o}lder's inequality alongside Lemmas \ref{D.6}, \eqref{4.54}, we deduce
\begin{align*}
\left\|\mathring{\mathsf{R}}_{\{cor\}}\right\|_{L_{t}^{1} L_{x}^{p_{1}}} & \lesssim\left\|w_{q+1}^{(p)} \otimes\left(w_{q+1}^{(c)}+w_{q+1}^{(t)}+w_{q+1}^{(o)}\right)-\left(w_{q+1}^{(c)}+w_{q+1}^{(t)}+w_{q+1}^{(o)}\right) \otimes w_{q+1}\right\|_{L_{t}^{1} L_{x}^{p_{1}}} \\
& \lesssim\left\|w_{q+1}^{(c)}+w_{q+1}^{(t)}+w_{q+1}^{(o)}\right\|_{L_{t}^{2} L_{x}^{p_{2}}}\left(\left\|w_{q+1}^{(p)}\right\|_{L_{t, x}^{2}}+\left\|w_{q+1}\right\|_{L_{t, x}^{2}}\right) \\
& \lesssim \delta_{q+1}^{\frac{1}{2}}\left(\vartheta_{q+1}^{-6} (r_{\perp} r_{\|})^{\frac{1}{p_{2}}-\frac{1}{2}}\frac{r_{\perp}}{r_{\|}}+\vartheta_{q+1}^{-2} \mu^{-1} \lambda r_{\perp}(r_{\perp}r_{\|})^{\frac{1}{p_{2}}-1}+\vartheta_{q+1}^{-7} \sigma^{-1}\right) \\
& \lesssim \vartheta_{q+1}^{-7} \delta_{q+1}^{\frac{1}{2}}\left(\lambda^{-8 \varepsilon-2 \widetilde{\eta}-12 \widetilde{\eta} \varepsilon}+\lambda^{2 \widetilde{\eta}-6\varepsilon-12 \widetilde{\eta} \varepsilon}+\lambda^{-2 \varepsilon}\right) \lesssim \vartheta_{q+1}^{-7} \lambda^{-\varepsilon} \label{5.28}\tag{4.25}
\end{align*}
where the last step is derived from the inequality $2 \widetilde{\eta}-6\varepsilon-12 \widetilde{\eta} \varepsilon \leq-\varepsilon $.

Therefore, combining the estimates \eqref{5.23}, \eqref{5.27}, and \eqref{5.28}, along with the observation that $\vartheta_{q+1}^{-17}\ll\lambda^{\frac{\varepsilon}{2}}$, we arrive at
\begin{align*}
\left\|\mathring{\mathsf{R}}_{q+1}\right\|_{L_{t, x}^{1}} & \leq\left\|\mathring{\mathsf{R}}_{\{l i n\}}\right\|_{L_{t}^{1} L_{x}^{\varrho}}+\left\|\mathring{\mathsf{R}}_{\{osc\}}\right\|_{L_{t}^{1} L_{x}^{\varrho}}+\left\|\mathring{\mathsf{R}}_{\{cor\}}\right\|_{L_{t}^{1} L_{x}^{p_{1}}} \\
& \lesssim \vartheta_{q+1}^{-8} \lambda^{-\varepsilon}+\vartheta_{q+1}^{-17} \lambda^{-\varepsilon}+\vartheta_{q+1}^{-7} \lambda^{-\varepsilon} \\
& \leq \lambda_{q+1}^{-2\varepsilon_{R}} \delta_{q+2}. \label{5.29}\tag{4.26}
\end{align*}

Hence, the $L_{t,x}^{1}$-estimate \eqref{2.10} of new Reynolds stress is finished.

\begin{center}
	\section{The super-critical regime $\mathcal{A}_2$}\label{E}
\end{center}\ \ \ \
In the present section, our primary focus is on the super-critical regime $\mathcal{A}_{2}$ with $\alpha \in[1,\frac{3}{2})$. When considering the other endpoint $(s, \frac{2\alpha}{2\alpha-1+s}, \infty)$, it might seem natural to employ the fundamental components from the previously discussed endpoint case. However, this approach unfortunately results in an absence of feasible parameter values due to the inclusion of the $\phi_{(k)}(t)$ term in accelerating jets, which leads to conflicting restrictions between the estimates of acceleration error \eqref{5.15} and oscillation error \eqref{5.25} associated with the temporal corrector $w_{q+1}^{t}$, this issue is also mentioned in Section 2.2 of \cite{YL22}. Drawing inspiration from the work in \cite{AC22}, we employ concentrated Mikado flows as spatial building blocks detailed in \eqref{6.3} below and the ``building blocks'' are indexed by four distinct parameters, $r_{\perp}, \lambda, \tau$ and $\sigma$ defined as follows:
\begin{equation*}
r_{\perp}:=\lambda_{q+1}^{-2\alpha+2-8 \varepsilon}, \lambda:=\lambda_{q+1}, \tau:=\lambda_{q+1}^{2 \alpha}, \sigma:=\lambda_{q+1}^{2 \varepsilon} \label{6.1}\tag{5.1}
\end{equation*}
where $\varepsilon$ is well-defined in \eqref{2.5}.

Although Mikado flows primarily exhibit 1D intermittency, the appropriate temporal building blocks can significantly achieve nearing 3D spatial intermittency when $\alpha$ approaches $\frac{3}{2}$. \\

\subsection{Spatial-temporal building blocks.}\ \ \ \
The concentrated Mikado flows are precisely defined by
$$
\mathbf{W}_{k}:=\psi_{r_{\perp}}\left(\lambda r_{\perp} N_{\Lambda} k \cdot x\right) k_{1}, \quad k \in \Lambda,
$$
where we maintain the same notations as in Section \ref{C}, specifically $\psi_{r_{\perp}}, r_{\perp}, N_{\Lambda}$, and the orthogonal basis $(k, k_{1})$. Additionally, we still use the temporal building blocks $g_{(k)}$ and $h_{(k)}$ introduced in \eqref{4.15}, albeit with a distinct choice of parameters $r_{\perp}$ and $\tau$ as specified in \eqref{6.1}.

Let $\boldsymbol{\boldsymbol{\Phi}}: \mathbb{R} \rightarrow \mathbb{R}$ be any smooth cut-off function supported on a ball of radius 1, and then normalize $\boldsymbol{\boldsymbol{\Phi}}$ such that $\psi:=\Delta \boldsymbol{\boldsymbol{\Phi}}$ and rescale it to
$$\mathbf{\Phi}_{r_{\perp}}(x)=r_{\perp}^{-1}\mathbf{\Phi}(\frac{x}{r_{\perp}}).$$
By an abuse of notation, $\boldsymbol{\boldsymbol{\Phi}}_{r_{\perp}}$ is also treated as a periodic function defined on $\mathbb{T}$.
Then, setting
\begin{align*}
& \psi_{(k)}(x):=\psi_{r_{\perp}}\left(\lambda r_{\perp} N_{\Lambda} k \cdot x\right), \label{6.2}\tag{5.2}\\
& \boldsymbol{\boldsymbol{\boldsymbol{\boldsymbol{\Phi}}}}_{(k)}(x):=\lambda^{-2}\boldsymbol{\boldsymbol{\boldsymbol{\Phi}}}_{r_{\perp}}\left(\lambda r_{\perp} N_{\Lambda} k \cdot x\right).
\end{align*}

Thus, the Mikado flows can be reformulated as
\begin{equation*}
\mathbf{W}_{k}= \psi_{(k)} k_{1}, \quad k \in \label{6.3}\Lambda. \tag{5.3}
\end{equation*}
Here, $\mathbf{W}_{k}$ is $(\mathbb{T}/\lambda r_{\perp})$-periodic in space.

To define the incompressible corrector of the perturbations, we employ the skew-symmetric potential $\mathbf{W}_{k}^{c}$, which is given by
\begin{equation*}
\mathbf{W}_{k}^{c}:=\nabla\boldsymbol{\boldsymbol{\boldsymbol{\Phi}}}_{(k)}\otimes k_{1}-k_{1}\otimes\nabla\boldsymbol{\boldsymbol{\boldsymbol{\Phi}}}_{(k)} \label{6.4}\tag{5.4}
\end{equation*}
and satisfying $\operatorname{div} \mathbf{W}_{k}^{c}=\mathbf{W}_{k}$.

In Lemma \ref{F.1}, we present the estimations of spatial building blocks, which are derived from the preceding Lemma \ref{D.1}.
\begin{Lemma}\label{F.1}(Estimates of Mikado flows). For any $p \in[1, \infty]$ and $N \in \mathbb{N}$, we have
\begin{equation*}
\left\|\nabla^{N} \psi_{(k)}\right\|_{L_{x}^{p}}+\lambda^{2}\left\|\nabla^{N} \boldsymbol{\boldsymbol{\boldsymbol{\Phi}}}_{(k)}\right\|_{L_{x}^{p}} \lesssim r_{\perp}^{\frac{1}{p}-\frac{1}{2}} \lambda^{N} \label{6.5}\tag{5.5}
\end{equation*}
where the implicit constants do not rely on $r_{\perp}$ and $\lambda$. Moreover, it holds that
\begin{equation*}
\left\|\nabla^{N} \mathbf{W}_{k}\right\|_{C_{t} L_{x}^{p}}+\lambda\left\|\nabla^{N} \mathbf{W}_{k}^{c}\right\|_{C_{t} L_{x}^{p}} \lesssim r_{\perp}^{\frac{1}{p}-\frac{1}{2}} \lambda^{N}, \quad k \in \Lambda .\label{6.6} \tag{5.6}
\end{equation*}
\end{Lemma}
\begin{Remark}
In the 2D case, the supports of concentrated Mikado flows are long strip in periodic domains, and the different flows would interact with each other even though the support of the smooth cut-off function $\psi$ is small enough. Rather than the spatial shifts used in the super-critical $\mathcal{A}_{1}$, we only use the temporal shifts to decouple different interactions which is different from the 2D case in \cite{AC22}. Due to the effect of non-intersecting temporal concentrated function $g_{(k)}(t)$, our ``building blocks" $g_{(k)}\mathbf{W}_{k}$ are guaranteed to be non-overlapping, without the necessary to consider the errors caused by support intersection in \cite{AC22}, which has also been applied in the results of hypo-dissipative compressible NSE \cite{JL69}.
\end{Remark}
\subsection{Estimates of velocity perturbations.}\ \ \ \
We recall the amplitudes of velocity perturbations in accordance with the previously established as follows:
\begin{equation*}
a_{(k)}(t, x):=\rho_{u}^{\frac{1}{2}}(t, x) f_{u}(t) \gamma_{(k)}\left(\operatorname{Id}-\frac{\mathring{\mathsf{R}}_{q}(t, x)}{\rho_{u}(t, x)}\right), \quad k \in \Lambda \label{6.7}\tag{5.7}
\end{equation*}
where $\rho_{u}, f_{u}, \gamma_{(k)}$ are defined as outlined in Section \ref{C.3.3}. It is worthy noting that the amplitudes $a_{(k)}$ adhere to the identical estimates as stated in Lemma \ref{D.5}. Precisely, we have
\begin{Lemma}\label{F.2}For $1 \leq N \leq 9, k \in \Lambda$, we have
\begin{align*}
\left\|a_{(k)}\right\|_{L_{t, x}^{2}} & \lesssim \delta_{q+1}^{\frac{1}{2}},  \label{6.8}\tag{5.8}\\
\left\|a_{(k)}\right\|_{C_{t, x}} & \lesssim \vartheta_{q+1}^{-1}, \quad\left\|a_{(k)}\right\|_{C_{t, x}^{N}} \lesssim \vartheta_{q+1}^{-5N} \label{6.9}\tag{5.9}
\end{align*}
where the implicit constants do not rely on variable $q$.
\end{Lemma}

Subsequently, we introduce the principal part of the velocity perturbations, which is defined as follows:
\begin{equation*}
w_{q+1}^{(p)}:=\sum_{k \in \Lambda} a_{(k)} g_{(k)} \mathbf{W}_{k}. \label{6.10}\tag{5.10}
\end{equation*}

We also introduce the incompressible corrector $w_{q+1}^{(c)}$ which is different from the previous one in \eqref{4.33}, defined by
\begin{equation*}
w_{q+1}^{(c)}:=\sum_{k \in \Lambda} g_{(k)} \mathbf{W}_{k}^{c} \nabla a_{(k)},\label{6.11}\tag{5.11}
\end{equation*}
where $\mathbf{W}_{k}^{c}$ is the potential given by \eqref{6.4}. Note that, due to the fact $\operatorname{div} \mathbf{W}_{k}^{c}=\mathbf{W}_{k}$, we obtain that
\begin{equation*}
w_{q+1}^{(p)}+w_{q+1}^{(c)}=\operatorname{div} \left(\sum_{k \in \Lambda}a_{(k)} g_{(k)} \mathbf{W}_{k}^{c}\right). \label{6.12}\tag{5.12}
\end{equation*}

Since  $a_{(k)} g_{(k)} \mathbf{W}_{k}^c$ exhibits skew-symmetric for all $k \in \Lambda$, one has
$$
\operatorname{div}\left(w_{q+1}^{(p)}+w_{q+1}^{(c)}\right)=0,
$$
which validates the incompressibility of the velocity perturbations.

Unlike accelerating jets, the absence of $\varphi_{r_{\|}}$ eliminates the necessity for introducing the temporal corrector $w_{q+1}^{(t)}$ as previous \eqref{4.36} to counterbalance high spatial frequency oscillation. Instead, we only rely on the temporal corrector
\begin{equation*}
w_{q+1}^{(o)}:=-\sigma^{-1} \sum_{k \in \Lambda} \mathbb{P}_{H} \mathbb{P}_{\neq 0}\left(h_{(k)} \fint_{\mathbb{T}^{2}} \mathbf{W}_{k} \otimes \mathbf{W}_{k} \mathrm{d} x \nabla\left(a_{(k)}^{2}\right)\right) \label{6.13}\tag{5.13}
\end{equation*}
to balance the high temporal frequency oscillation. Combining \eqref{4.14} and \eqref{6.13}, along with Leibniz's rule, it becomes evident that the algebraic identity \eqref{4.39} remains valid.

We are now in a position to define the velocity perturbation $w_{q+1}$ at level $q+1$ through the following identity:
\begin{equation*}
w_{q+1}:=w_{q+1}^{(p)}+w_{q+1}^{(c)}+w_{q+1}^{(o)}. \label{6.14}\tag{5.14}
\end{equation*}
Similarly, the velocity field at level $q+1$ is given by
\begin{equation*}
u_{q+1}:=\widetilde{u}_{q}+w_{q+1} . \label{6.15}\tag{5.15}
\end{equation*}
where $\widetilde{u}_{q}$ presents the velocity field that has been prepared during the gluing stage, as described in Section \ref{B}. Through these constructions, it can be shown that $w_{q+1}$ is both mean-free and divergence-free.

Furthermore, analogous to the estimates provided in Lemma \ref{D.6}, we can derive the following estimates for the velocity perturbations defined above.

\begin{Lemma}\label{F.3} (Estimates of perturbations). Given any $p \in(1, \infty), \gamma \in[1, \infty]$ and integers $0 \leq N \leq 7$, the following estimates hold true:
\begin{align*}
&\left\|\nabla^{N} w_{q+1}^{(p)}\right\|_{L_{t}^{\gamma} L_{x}^{p}} \lesssim \vartheta_{q+1}^{-1} \lambda^{N} r_{\perp}^{\frac{1}{p}-\frac{1}{2}} \tau^{\frac{1}{2}-\frac{1}{\gamma}},  \label{6.16}\tag{5.16}\\
&\left\|\nabla^{N} w_{q+1}^{(c)}\right\|_{L_{t}^{\gamma} L_{x}^{p}} \lesssim \vartheta_{q+1}^{-6} \lambda^{N-1} r_{\perp}^{\frac{1}{p}-\frac{1}{2}} \tau^{\frac{1}{2}-\frac{1}{\gamma}},  \label{6.17}\tag{5.17}\\
&\left\|\nabla^{N} w_{q+1}^{(o)}\right\|_{L_{t}^{\gamma} L_{x}^{p}} \lesssim \vartheta_{q+1}^{-5 N-7} \sigma^{-1}, \label{6.18}\tag{5.18}
\end{align*}
where the implicit constants depend only on $N$, $\gamma$ and $\varrho$. In particular, for integers $1 \leq N \leq 7$, we have
\begin{align*}
& \left\|w_{q+1}^{(p)}\right\|_{L_{t}^{\infty} H_{x}^{N}}+\left\|w_{q+1}^{(c)}\right\|_{L_{t}^{\infty} H_{x}^{N}}+\left\|w_{q+1}^{(o)}\right\|_{L_{t}^{\infty} H_{x}^{N}} \lesssim \lambda^{N+2},  \label{6.19}\tag{5.19}\\
& \left\|\partial_{t} w_{q+1}^{(p)}\right\|_{L_{t}^{\infty} H_{x}^{N}}+\left\|\partial_{t} w_{q+1}^{(c)}\right\|_{L_{t}^{\infty} H_{x}^{N}}+\left\|\partial_{t} w_{q+1}^{(o)}\right\|_{L_{t}^{\infty} H_{x}^{N}} \lesssim \lambda^{N+5},\label{6.20} \tag{5.20}
\end{align*}
where the implicit constants are independent of $\lambda$.
\end{Lemma}
{\bf Proof.} By Lemmas \ref{D.3}, \ref{F.1}, \ref{F.2}, \eqref{6.10} and \eqref{6.11}, we obtain
\begin{align*}
\left\|\nabla^{N} w_{q+1}^{(p)}\right\|_{L_{t}^{\gamma} L_{x}^{p}} & \lesssim \sum_{k \in \Lambda} \sum_{N_{1}+N_{2}=N}\left\|a_{(k)}\right\|_{C_{t, x}^{N_{1}}}\left\|g_{(k)}\right\|_{L_{t}^{\gamma}}\left\|\nabla^{N_{2}} \mathbf{W}_{k}\right\|_{C_{t} L_{x}^{p}} \\
& \lesssim \vartheta_{q+1}^{-1} \lambda^{N} r_{\perp}^{\frac{1}{p}-\frac{1}{2}} \tau^{\frac{1}{2}-\frac{1}{\gamma}} \label{6.21}\tag{5.21}
\end{align*}
and
$$
\begin{aligned}
\left\|\nabla^{N} w_{q+1}^{(c)}\right\|_{L_{t}^{\gamma} L_{x}^{p}} & \lesssim \sum_{k \in \Lambda}\sum_{N_{1}+N_{2}=N}\left\|g_{(k)}\right\|_{L_{t}^{\gamma}} \left\|a_{(k)}\right\|_{C_{t, x}^{N_{1}+1}}\left\|\nabla^{N_{2}} \mathbf{W}_{k}^{c}\right\|_{C_{t} L_{x}^{p}} \\
& \lesssim \vartheta_{q+1}^{-6} \lambda^{N-1} r_{\perp}^{\frac{1}{p}-\frac{1}{2}} \tau^{\frac{1}{2}-\frac{1}{\gamma}},
\end{aligned}
$$
which verifies \eqref{6.16} and \eqref{6.17}.

Regarding the temporal corrector $w_{q+1}^{(o)}$, by using Lemmas \ref{D.3}, \ref{F.2} and \eqref{6.13},
$$
\left\|\nabla^{N} w_{q+1}^{(o)}\right\|_{L_{t}^{\gamma} L_{x}^{p}} \lesssim \sigma^{-1} \sum_{k \in \Lambda}\left\|h_{(k)}\right\|_{C_{t}}\left\|\nabla^{N+1}\left(a_{(k)}^{2}\right)\right\|_{C_{t, x}} \lesssim \vartheta_{q+1}^{-5N-7} \sigma^{-1}.
$$

Turning to the $L_{t}^{\infty} H_{x}^{N}$-estimates of velocity perturbations, using \eqref{2.5}, \eqref{6.1} and \eqref{6.16}-\eqref{6.18}, we get
$$
\begin{aligned}
& \left\|w_{q+1}^{(p)}\right\|_{L_{t}^{\infty} H_{x}^{N}}+\left\|w_{q+1}^{(c)}\right\|_{L_{t}^{\infty} H_{x}^{N}}+\left\|w_{q+1}^{(o)}\right\|_{L_{t}^{\infty} H_{x}^{N}} \\
\lesssim & \vartheta_{q+1}^{-1} \lambda^{N} \tau^{\frac{1}{2}}+\vartheta_{q+1}^{-6} \lambda^{N-1} \tau^{\frac{1}{2}}+\vartheta_{q+1}^{-5N-7} \sigma^{-1} \\
\lesssim & \vartheta_{q+1}^{-1} \lambda^{\frac{\alpha}{2}+N}+\vartheta_{q+1}^{-6} \lambda^{\frac{\alpha}{2}+N-1}+\vartheta_{q+1}^{-5N-7} \lambda^{-2 \varepsilon} \lesssim \lambda^{N+2}
\end{aligned}
$$
which implies \eqref{6.19}.

For \eqref{6.20}, by \eqref{6.1} and Lemmas \ref{D.3}, \ref{F.3}, \ref{F.2}, we obtain
\begin{equation*}
\left\|\partial_{t} w_{q+1}^{(p)}\right\|_{L_{t}^{\infty} H_{x}^{N}} \lesssim \sum_{k \in \Lambda}\left\|a_{(k)}\right\|_{C_{t, x}^{N+1}}\left\|\partial_{t} g_{(k)}\right\|_{L_{t}^{\infty}}\left\|\mathbf{W}_{k}\right\|_{L_{t}^{\infty} H_{x}^{N}} \lesssim \vartheta_{q+1}^{-5N-7} \lambda^{N} \sigma \tau^{\frac{3}{2}} , \label{6.22}\tag{5.22}
\end{equation*}
and
\begin{align*}
\left\|\partial_{t} w_{q+1}^{(c)}\right\|_{L_{t}^{\infty} H_{x}^{N}} & \lesssim \sum_{k \in \Lambda}\left\|a_{(k)}\right\|_{C_{t, x}^{N+2}}\left\|\partial_{t} g_{(k)}\right\|_{C_{t}}\left\|\mathbf{W}_{k}^{c}\right\|_{L_{t}^{\infty} H_{x}^{N}} \\
& \lesssim \vartheta_{q+1}^{-5 N-12} \lambda^{N-1} \sigma \tau^{\frac{3}{2}} . \label{6.23}\tag{5.23}
\end{align*}

As $\mathbb{P}_{H} \mathbb{P}_{\neq 0}$ is bounded in $H_{x}^{N}$, an analogy can be drawn to \eqref{4.52}, leading to
\begin{equation*}
\left\|\partial_{t} w_{q+1}^{(o)}\right\|_{L_{t}^{\infty} H_{x}^{N}} \lesssim \sigma^{-1} \sum_{k \in \Lambda}\left\|\partial_{t}\left(h_{(k)} \nabla\left(a_{(k)}^{2}\right)\right)\right\|_{L_{t}^{\infty} H_{x}^{N}} \lesssim \vartheta_{q+1}^{-5N-7} \tau . \label{6.24}\tag{5.24}
\end{equation*}

Hence, based on the fact $\vartheta_{q+1}^{-5 N-12} \leq \lambda^{N \varepsilon / 2}$ and $0<\varepsilon \leq(3-2\alpha) / 20$, we arrive at the conclusion that
$$
\begin{aligned}
&\left\|\partial_{t} w_{q+1}^{(p)}\right\|_{L_{t}^{\infty} H_{x}^{N}}+\left\|\partial_{t} w_{q+1}^{(c)}\right\|_{L_{t}^{\infty} H_{x}^{N}}+\left\|\partial_{t} w_{q+1}^{(o)}\right\|_{L_{t}^{\infty} H_{x}^{N}} \\
& \lesssim \vartheta_{q+1}^{-5N-7} \lambda^{N} \sigma \tau^{\frac{3}{2}}+\vartheta_{q+1}^{-5 N-12} \lambda^{N-1} \sigma \tau^{\frac{3}{2}}+\vartheta_{q+1}^{-5N-7} \tau \\
& \lesssim \vartheta_{q+1}^{-5N-7} \lambda^{3 \alpha+N+2 \varepsilon}+\vartheta_{q+1}^{-5 N-12} \lambda^{3 \alpha+N-1+2 \varepsilon}+\vartheta_{q+1}^{-5N-7} \lambda^{2 \alpha} \lesssim \lambda^{N+5}.
\end{aligned}
$$

Consequently, the proof of Lemma \ref{F.3} is finished.

{\bf Verification of the inductive estimates for velocity}. In this part, we aim to verify the inductive estimates for velocity. Employing the $L^{p}$ decorrelation Lemma \ref{H.2} with $f=a_{(k)}, g=g_{(k)} \psi_{(k)}$ and $\sigma=\lambda^{2 \varepsilon}$, combining \eqref{2.2} and \eqref{2.3}, along with Lemmas \ref{D.3}, \ref{F.1}, and \ref{F.3}, we obtain
\begin{align*}
\left\|w_{q+1}^{(p)}\right\|_{L_{t, x}^{2}} & \lesssim \sum_{k \in \Lambda}\left(\left\|a_{(k)}\right\|_{L_{t, x}^{2}}\left\|g_{(k)}\right\|_{L_{t}^{2}}\left\|\psi_{(k)}\right\|_{C_{t} L_{x}^{2}}+\sigma^{-\frac{1}{2}}\left\|a_{(k)}\right\|_{C_{t, x}^{1}}\left\|g_{(k)}\right\|_{L_{t}^{2}}\left\|\psi_{(k)}\right\|_{C_{t} L_{x}^{2}}\right) \\
& \lesssim \delta_{q+1}^{\frac{1}{2}}+\vartheta_{q+1}^{-6} \lambda_{q+1}^{-\varepsilon} \lesssim \delta_{q+1}^{\frac{1}{2}} . \label{6.25}\tag{5.25}
\end{align*}

Then, by \eqref{2.3}, \eqref{6.25} and Lemma \ref{F.3} we have
\begin{align*}
\left\|w_{q+1}\right\|_{L_{t, x}^{2}} & \lesssim\left\|w_{q+1}^{(p)}\right\|_{L_{t, x}^{2}}+\left\|w_{q+1}^{(c)}\right\|_{L_{t, x}^{2}}+\left\|w_{q+1}^{(o)}\right\|_{L_{t, x}^{2}} \\
& \lesssim \delta_{q+1}^{\frac{1}{2}}+\vartheta_{q+1}^{-6} \lambda^{-1}+\vartheta_{q+1}^{-7} \sigma^{-1} \lesssim \delta_{q+1}^{\frac{1}{2}} \label{6.26}\tag{5.26}
\end{align*}
and
\begin{align*}
\left\|w_{q+1}\right\|_{L_{t}^{1} L_{x}^{2}} & \lesssim\left\|w_{q+1}^{(p)}\right\|_{L_{t}^{1} L_{x}^{2}}+\left\|w_{q+1}^{(c)}\right\|_{L_{t}^{1} L_{x}^{2}}+\left\|w_{q+1}^{(o)}\right\|_{L_{t}^{1} L_{x}^{2}} \\
& \lesssim \vartheta_{q+1}^{-1} \tau^{-\frac{1}{2}}+\vartheta_{q+1}^{-6} \lambda^{-1} \tau^{-\frac{1}{2}}+\vartheta_{q+1}^{-7} \sigma^{-1} \lesssim \lambda_{q+1}^{-\varepsilon} . \label{6.27}\tag{5.27}
\end{align*}

We are now in a position to validate the iterative estimates for $u_{q+1}$.

Based on \eqref{2.6}, \eqref{3.5}, \eqref{3.12}, \eqref{6.15}, and \eqref{6.19}, it can be inferred that for a sufficiently large $a$, yields that
\begin{align*}
\left\|u_{q+1}\right\|_{L_{t}^{\infty} H_{x}^{3}} & \lesssim\left\|\widetilde{u}_{q}\right\|_{L_{t}^{\infty} H_{x}^{2}}+\left\|w_{q+1}\right\|_{L_{t}^{\infty} H_{x}^{3}} \\
& \lesssim \lambda_{q}^{5}+\lambda_{q+1}^{5} \lesssim \lambda_{q+1}^{5}, \label{5.28} \tag{5.28}\\
\left\|\partial_{t} u_{q+1}\right\|_{L_{t}^{\infty} H_{x}^{3}} & \lesssim\left\|\partial_{t} \widetilde{u}_{q}\right\|_{L_{t}^{\infty} H_{x}^{2}}+\left\|\partial_{t} w_{q+1}\right\|_{L_{t}^{\infty} H_{x}^{2}} \\
& \lesssim \vartheta_{q+1}^{-1} \lambda_{q}^{5}+\ell_{q+1} \lambda_{q}^{5}+\lambda_{q+1}^{5} \lesssim \lambda_{q+1}^{5}. \label{6.29}\tag{5.29}
\end{align*}

Additionally, by \eqref{3.5}, \eqref{6.26} and \eqref{6.27}, we can deduce that
\begin{align*}
\left\|u_{q}-u_{q+1}\right\|_{L_{t, x}^{2}} & \leq\left\|u_{q}-\widetilde{u}_{q}\right\|_{L_{t, x}^{2}}+\left\|\widetilde{u}_{q}-u_{q+1}\right\|_{L_{t, x}^{2}} \\
& \lesssim\left\|u_{q}-\widetilde{u}_{q}\right\|_{L_{t}^{\infty} L_{x}^{2}}+\left\|w_{q+1}\right\|_{L_{t, x}^{2}} \\
& \lesssim \lambda_{q}^{-3}+\delta_{q+1}^{\frac{1}{2}} \leq M^{*} \delta_{q+1}^{\frac{1}{2}}, \label{6.30}\tag{5.30}
\end{align*}
for $M^{*}$ large enough  and
\begin{align*}
\left\|u_{q}-u_{q+1}\right\|_{L_{t}^{1} L_{x}^{2}} & \lesssim\left\|u_{q}-\widetilde{u}_{q}\right\|_{L_{t}^{\infty} L_{x}^{2}}+\left\|w_{q+1}\right\|_{L_{t}^{1} L_{x}^{2}} \\
& \lesssim \lambda_{q}^{-3}+\lambda_{q+1}^{-\varepsilon} \leq \delta_{q+2}^{\frac{1}{2}}. \label{6.31}\tag{5.31}
\end{align*}

Considering the iteration estimate \eqref{2.13}, we claim that the Sobolev embedding
\begin{equation*}
H_{x}^{3} \hookrightarrow W_{x}^{s, p} \label{6.32}\tag{5.32}
\end{equation*}
holds true for any triplet $(s, p, \gamma) \in \mathcal{A}_{2}$. Thus, analogous to \eqref{4.61}, by virtue of \eqref{3.7}, \eqref{3.16} and \eqref{6.32}, we arrive at
$$
\left\|\widetilde{u}_{q}-u_{q}\right\|_{L_{t}^{\gamma} W_{x}^{s, p}} \lesssim\left\|\sum_{i} \mathscr{X} _{i}\left(v_{i}-u_{q}\right)\right\|_{L_{t}^{\infty} H_{x}^{3}} \lesssim \lambda_{q}^{-2}.
$$

Hence, for all $\alpha \in[1,\frac{3}{2})$, we can take advantage of \eqref{2.3}, \eqref{6.1} and \eqref{6.32}, along with Lemma \ref{F.3}, to deduce
\begin{align*}
\left\|u_{q+1}-u_{q}\right\|_{L_{t}^{\gamma} W_{x}^{s, p}} & \lesssim\left\|\widetilde{u}_{q}-u_{q}\right\|_{L_{t}^{\gamma} W_{x}^{s, p}}+\left\|w_{q+1}\right\|_{L_{t}^{\gamma} W_{x}^{s, p}} \\
& \lesssim \lambda_{q}^{-2}+\vartheta_{q+1}^{-1} \lambda_{q+1}^{s} r_{\perp}^{\frac{1}{p}-\frac{1}{2}} \tau^{\frac{1}{2}-\frac{1}{\gamma}}+\vartheta_{q+1}^{-17} \sigma^{-1} \\
& \lesssim \lambda_{q}^{-2}+\lambda_{q+1}^{s+2 \alpha-1-\frac{2 \alpha}{\gamma}-\frac{2 \alpha-2}{p}+\varepsilon\left(5-\frac{8}{p}\right)}+\lambda_{q+1}^{-\varepsilon} . \label{6.35}\tag{5.33}
\end{align*}
Due to \eqref{2.5},
\begin{equation*}
s+2 \alpha-1-\frac{2 \alpha}{\gamma}-\frac{2 \alpha-2}{p}+\varepsilon\left(5-\frac{8}{p}\right) \leq s+2 \alpha-1-\frac{2 \alpha}{\gamma}-\frac{2 \alpha-2}{p}+5 \varepsilon<-15 \varepsilon, \label{6.36}\tag{5.34}
\end{equation*}
we thus obtain
\begin{equation*}
\left\|u_{q+1}-u_{q}\right\|_{L_{t}^{\gamma} W_{x}^{s, p}} \leq \delta_{q+2}^{\frac{1}{2}} . \label{6.37}\tag{5.35}
\end{equation*}

Hence, the iterative estimates \eqref{2.6}, \eqref{2.7}, \eqref{2.11} and \eqref{2.13} have been rigorously verified.

\subsection{Estimates of Reynolds stress.}\ \ \ \
In the following discussion, we focus on the Reynolds stress pertaining to the endpoint case $(s, \frac{2\alpha}{2\alpha-1+s}, \infty)$. Utilizing the framework established in \eqref{2.1} at the $q+1$ level, we deduce that the new Reynolds stress obeys
$$
\begin{aligned}
& \operatorname{div} \mathring{\mathsf{R}}_{q+1}-\nabla(\mathsf{P}_{q+1}-\mathsf{P}_{q})=\underbrace{\partial_{t}\left(w_{q+1}^{(p)}+w_{q+1}^{(c)}\right)+(-\Delta)^{\alpha} w_{q+1}+\operatorname{div}\left(\widetilde{u}_{q} \otimes w_{q+1}+w_{q+1} \otimes \widetilde{u}_{q}\right)}_{\operatorname{div} \mathring{\mathsf{R}}_{\{l i n\}}+\nabla\mathsf{P}_{\{l i n\}}}
\end{aligned}
$$
\begin{align*}
& +\underbrace{\operatorname{div}\left(w_{q+1}^{(p)} \otimes w_{q+1}^{(p)}+\mathring{\widetilde{\mathsf{R}}}_q\right)+\partial_t w_{q+1}^{(o)}}_{\operatorname{div} \mathring{\mathsf{R}}_{\{osc\}}+\nabla\mathsf{P}_{\{o s c\}}}\label{6.38}\tag{5.36}
\end{align*}
$$
\begin{aligned}
& +\underbrace{\operatorname{div}\left(\left(w_{q+1}^{(c)}+w_{q+1}^{(o)}\right) \otimes w_{q+1}+w_{q+1}^{(p)} \otimes\left(w_{q+1}^{(c)}+w_{q+1}^{(o)}\right)\right)}_{\operatorname{div} \mathring{\mathsf{R}}_{\{cor\}}+\nabla\mathsf{P}_{\{c o r\}}}.
\end{aligned}
$$

Subsequently, by employing the inverse divergence operator $\mathcal{R}$, we can select the Reynolds stress at the $q+1$ level defined as
\begin{equation*}
\mathring{\mathsf{R}}_{q+1}:=\mathring{\mathsf{R}}_{\{l i n\}}+\mathring{\mathsf{R}}_{\{osc\}}+\mathring{\mathsf{R}}_{\{cor\}} , \label{6.39}\tag{5.37}
\end{equation*}
where the linear error
\begin{equation*}
\mathring{\mathsf{R}}_{\{l i n\}}:=\mathcal{R}\left(\partial_{t}\left(w_{q+1}^{(p)}+w_{q+1}^{(c)}\right)\right)+\mathcal{R}(-\Delta)^{\alpha} w_{q+1}+\mathcal{R} \mathbb{P}_{H} \operatorname{div}\left(\widetilde{u}_{q} \mathring{\otimes} w_{q+1}+w_{q+1} \mathring{\otimes} \widetilde{u}_{q}\right), \label{6.40}\tag{5.38}
\end{equation*}
the oscillation error
\begin{align*}
\mathring{\mathsf{R}}_{\{osc\}}:= & \sum_{k \in \Lambda} \mathcal{R} \mathbb{P}_{H} \mathbb{P}_{\neq 0}\left(g_{(k)}^{2} \mathbb{P}_{\neq 0}\left(\mathbf{W}_{k} \otimes \mathbf{W}_{k}\right) \nabla\left(a_{(k)}^{2}\right)\right) \\
& -\sigma^{-1} \sum_{k \in \Lambda} \mathcal{R} \mathbb{P}_{H} \mathbb{P}_{\neq 0}\left(h_{(k)} \fint_{\mathbb{T}^{2}} \mathbf{W}_{k} \otimes \mathbf{W}_{k} \mathrm{d} x \partial_{t} \nabla\left(a_{(k)}^{2}\right)\right) , \label{6.41}\tag{5.39}
\end{align*}
and the corrector error
\begin{equation*}
\mathring{\mathsf{R}}_{\{cor\}}:=\mathcal{R} \mathbb{P}_{H} \operatorname{div}\left(w_{q+1}^{(p)} \mathring{\otimes}\left(w_{q+1}^{(c)}+w_{q+1}^{(o)}\right)+\left(w_{q+1}^{(c)}+w_{q+1}^{(o)}\right) \mathring{\otimes} w_{q+1}\right) . \label{6.42}\tag{5.40}
\end{equation*}

In the subsequent section, we proceed to verify the inductive estimates pertaining to the new Reynolds stress  $\mathring{\mathsf{R}}_{q+1}$.

{\bf Verification of the $L_{t}^{\infty} H_{x}^{N}$-estimate for the Reynolds stress}. Utilizing the bounds established in \eqref{6.19} and \eqref{6.20}, it can be observed that the velocity perturbations adhere to the same upper bound as those specified in \eqref{4.46} and \eqref{4.47}. From this similarity, we adopt a similar reasoning as employed in \eqref{5.7}, leading to the derivation of \eqref{2.8} and \eqref{2.9} within the super-critical regime denoted by $\mathcal{A}_{2}$. For brevity, the detailed derivations are omitted here.

{\bf Verification of $L_{t, x}^{1}$-decay of Reynolds stress}. In order to establish the $L_{t, x}^{1}$-decay \eqref{2.10} of the new Reynolds stress $\mathring{\mathsf{R}}_{q+1}$ at level $q+1$, we also select
\begin{equation*}
\varrho:=\frac{2 \alpha-2+8 \varepsilon}{2 \alpha-2+7 \varepsilon} \in(1,2), \label{6.43}\tag{5.41}
\end{equation*}
where $\varepsilon$ is defined in \eqref{2.5}. Then, we have
\begin{equation*}
(2-2\alpha-8 \varepsilon)\left(\frac{1}{\varrho}-\frac{1}{2}\right)=1-\alpha-3 \varepsilon \label{6.44}\tag{5.42}
\end{equation*}
and
\begin{equation*}
r_{\perp}^{\frac{1}{\varrho}-\frac{1}{2}}=\lambda^{1-\alpha-3 \varepsilon} . \label{6.45}\tag{5.43}
\end{equation*}
{\bf (i) Linear Error}. Using Lemmas \ref{D.2}, \ref{F.1}, \ref{F.2}, \eqref{2.5}, \eqref{6.1}, \eqref{6.12} and \eqref{6.45}, we get
\begin{align*}
&\left\|\mathcal{R} \partial_{t}\left(w_{q+1}^{(p)}+w_{q+1}^{(c)}\right)\right\|_{L_{t}^{1} L_{x}^{\varrho}} \\
& \lesssim \sum_{k \in \Lambda}\left\|\mathcal{R} \operatorname{div} \partial_{t}\left(g_{(k)} a_{(k)} \mathbf{W}_{k}^{c}\right)\right\|_{L_{t}^{1} L_{x}^{\varrho}} \\
& \lesssim \sum_{k \in \Lambda}\left(\left\|g_{(k)}\right\|_{L_{t}^{1}}\left\|a_{(k)}\right\|_{C_{t, x}^{1}}+\left\|\partial_{t} g_{(k)}\right\|_{L_{t}^{1}}\left\|a_{(k)}\right\|_{C_{t, x}}\right)\left\|\mathbf{W}_{k}^{c}\right\|_{C_{t} L_{x}^{\varrho}} \\
& \lesssim \vartheta_{q+1}^{-6} (\tau^{-\frac{1}{2}}+ \sigma \tau^{\frac{1}{2}}) \lambda^{-1} r_{\perp}^{\frac{1}{\varrho}-\frac{1}{2}} \lesssim \vartheta_{q+1}^{-6} \lambda^{- \varepsilon} . \label{6.46}\tag{5.44}
\end{align*}

For the viscosity term $(-\Delta)^{\alpha} w_{q+1}$ , by \eqref{6.14},
\begin{align*}
\left\|\mathcal{R}(-\Delta)^{\alpha} w_{q+1}\right\|_{L_{t}^{1} L_{x}^{\varrho}} &\lesssim\left\|\mathcal{R}(-\Delta)^{\alpha} w_{q+1}^{(p)}\right\|_{L_{t}^{1} L_{x}^{\varrho}}\\&+\left\|\mathcal{R}(-\Delta)^{\alpha} w_{q+1}^{(c)}\right\|_{L_{t}^{1} L_{x}^{\varrho}}+\left\|\mathcal{R}(-\Delta)^{\alpha} w_{q+1}^{(o)}\right\|_{L_{t}^{1} L_{x}^{\varrho}}. \label{6.47}\tag{5.45}
\end{align*}
By the interpolation estimate, \eqref{6.1}, \eqref{6.45} along with Lemma \ref{F.3},
\begin{align*}
\left\|\mathcal{R}(-\Delta)^{\alpha} w_{q+1}^{(p)}\right\|_{L_{t}^{1} L_{x}^{\varrho}} & \lesssim\left\||\nabla|^{2 \alpha-1} w_{q+1}^{(p)}\right\|_{L_{t}^{1} L_{x}^{\varrho}} \\
& \lesssim\left\|w_{q+1}^{(p)}\right\|_{L_{t}^{1} L_{x}^{\varrho}}^{\frac{4-2 \alpha}{3}}\left\|w_{q+1}^{(p)}\right\|_{L_{t}^{1} W_{x}^{3, \varrho}}^{\frac{2 \alpha-1}{3}} \\
& \lesssim \vartheta_{q+1}^{-1} \lambda^{2 \alpha-1} r_{\perp}^{\frac{1}{\varrho}-\frac{1}{2}} \tau^{-\frac{1}{2}} \lesssim \vartheta_{q+1}^{-1} \lambda^{-3 \varepsilon} . \label{6.48}\tag{5.46}
\end{align*}
Similarly,
\begin{align*}
& \left\|\mathcal{R}(-\Delta)^{\alpha} w_{q+1}^{(c)}\right\|_{L_{t}^{1} L_{x}^{\varrho}} \lesssim \vartheta_{q+1}^{-6} \lambda^{2 \alpha-2} r_{\perp}^{\frac{1}{\varrho}-\frac{1}{2}} \tau^{-\frac{1}{2}} \lesssim \vartheta_{q+1}^{-6} \lambda^{-1-3 \varepsilon},  \label{6.49}\tag{5.47}\\
& \left\|\mathcal{R}(-\Delta)^{\alpha} w_{q+1}^{(o)}\right\|_{L_{t}^{1} L_{x}^{\varrho}} \lesssim \vartheta_{q+1}^{-17} \sigma^{-1} \lesssim \vartheta_{q+1}^{-17} \lambda^{-2 \varepsilon}. \label{6.50}\tag{5.48}
\end{align*}

After combining estimates \eqref{6.47} to \eqref{6.50}, we arrive at
\begin{equation*}
\left\|\mathcal{R}(-\Delta)^{\alpha} w_{q+1}\right\|_{L_{t}^{1} L_{x}^{\varrho}} \lesssim \vartheta_{q+1}^{-17} \lambda^{-2 \varepsilon}. \label{6.51}\tag{5.49}
\end{equation*}

Furthermore, by using \eqref{3.4}, \eqref{6.45} and Lemma \ref{F.3}, we can establish the following estimate:
\begin{align*}
& \left\|\mathcal{R} \mathbb{P}_{H} \operatorname{div}\left(w_{q+1} \otimes \widetilde{u}_{q}+\widetilde{u}_{q} \otimes w_{q+1}\right)\right\|_{L_{t}^{1} L_{x}^{\varrho}} \\
\lesssim & \left\|w_{q+1} \otimes \widetilde{u}_{q}+\widetilde{u}_{q} \otimes w_{q+1}\right\|_{L_{t}^{1} L_{x}^{\varrho}} \\
\lesssim & \left\|\widetilde{u}_{q}\right\|_{L_{t}^{\infty} H_{x}^{3}}\left\|w_{q+1}\right\|_{L_{t}^{1} L_{x}^{\varrho}} \\
\lesssim & \lambda_{q}^{5}\left(\vartheta_{q+1}^{-1} r_{\perp}^{\frac{1}{\varrho}-\frac{1}{2}} \tau^{-\frac{1}{2}}+\vartheta_{q+1}^{-7} \sigma^{-1}\right) \lesssim \vartheta_{q+1}^{-8} \lambda^{-2 \varepsilon} . \label{6.52}\tag{5.50}
\end{align*}

Consequently, combining the estimates from \eqref{6.46}, \eqref{6.51}, and \eqref{6.52}, we have
\begin{equation*}
\left\|\mathring{\mathsf{R}}_{\{l i n\}}\right\|_{L_{t}^{1} L_{x}^{\varrho}} \lesssim \vartheta_{q+1}^{-6} \lambda^{- \varepsilon}+\vartheta_{q+1}^{-17} \lambda^{-2 \varepsilon}+\vartheta_{q+1}^{-8} \lambda^{-2 \varepsilon} \lesssim \vartheta_{q+1}^{-17} \lambda^{- \varepsilon}. \label{6.53}\tag{5.51}
\end{equation*}
{\bf (ii) Oscillation Error}. Different from the previous endpoint case, we only need the decomposition of the oscillation error into two distinct components:
$$
\mathring{\mathsf{R}}_{\{osc\}}=\mathring{\mathsf{R}}_{\{o s c .1\}}+\mathring{\mathsf{R}}_{\{o s c .2\}} \text {, }
$$
where the low-high spatial oscillation error
$$
\mathring{\mathsf{R}}_{\{o s c .1\}}:=\sum_{k \in \Lambda} \mathcal{R} \mathbb{P}_{H} \mathbb{P}_{\neq 0}\left(g_{(k)}^{2} \mathbb{P}_{\neq 0}\left(\mathbf{W}_{k} \otimes \mathbf{W}_{k}\right) \nabla\left(a_{(k)}^{2}\right)\right)
$$
and the low frequency error
$$
\mathring{\mathsf{R}}_{\{o s c .2\}}:=-\sigma^{-1} \sum_{k \in \Lambda} \mathcal{R}\mathbb{P}_{H} \mathbb{P}_{\neq 0}\left(h_{(k)} \fint_{\mathbb{T}^{2}} \mathbf{W}_{k} \otimes \mathbf{W}_{k} \mathrm{d} x \partial_{t} \nabla\left(a_{(k)}^{2}\right)\right).
$$

Then, employing Lemmas \ref{F.1},  \ref{F.2} and \ref{H.3} with $a=\nabla\left(a_{(k)}^{2}\right)$ and $f=\psi_{(k)}^{2}$, we have
\begin{align*}
\left\|\mathring{\mathsf{R}}_{\text {osc.} 1}\right\|_{L_{t}^{1} L_{x}^{\varrho}} & \lesssim \sum_{k \in \Lambda}\left\|g_{(k)}\right\|_{L_{t}^{2}}^{2}\left\||\nabla|^{-1} \mathbb{P}_{\neq 0}\left(\mathbb{P}_{\geq\left(\lambda r_{\perp} / 2\right)}\left(\mathbf{W}_{k} \otimes \mathbf{W}_{k}\right) \nabla\left(a_{(k)}^{2}\right)\right)\right\|_{C_{t} L_{x}^{\varrho}} \\
& \lesssim \sum_{k \in \Lambda}\left\||\nabla|^{3}\left(a_{(k)}^{2}\right)\right\|_{C_{t, x}} (\lambda r_{\perp})^{-1}\left\|\psi_{(k)}^{2}\right\|_{C_{t} L_{x}^{\varrho}} \\
& \lesssim \vartheta_{q+1}^{-17} \lambda^{-1} r_{\perp}^{\frac{1}{p}-2} . \label{6.54}\tag{5.52}
\end{align*}

Furthermore, similar to the case in \eqref{5.26}, the estimates of the low-frequency component $\mathring{\mathsf{R}}_{\{osc. 2\} }$ can be achieved by utilizing Lemmas \ref{D.3}, \ref{F.2},
\begin{equation*}
\left\|\mathring{\mathsf{R}}_{\{o s c .2\}}\right\|_{L_{t}^{1} L_{x}^{\varrho}} \lesssim \sigma^{-1} \sum_{k \in \Lambda}\left\|h_{(k)}\right\|_{C_{t}}\left(\left\|a_{(k)}\right\|_{C_{t, x}}\left\|a_{(k)}\right\|_{C_{t, x}^{2}}+\left\|a_{(k)}\right\|_{C_{t, x}^{1}}^{2}\right) \lesssim \vartheta_{q+1}^{-12} \sigma^{-1}. \label{6.55}\tag{5.53}
\end{equation*}

Therefore, by combining \eqref{6.54} and \eqref{6.55} with the utilization of \eqref{2.5}, \eqref{6.1} and \eqref{6.45}, we arrive at
\begin{align*}
\left\|\mathring{\mathsf{R}}_{\{osc\}}\right\|_{L_{t}^{1} L_{x}^{\varrho}} & \lesssim \vartheta_{q+1}^{-17} \lambda^{-1} r_{\perp}^{\frac{1}{\varrho}-2}+\vartheta_{q+1}^{-12} \sigma^{-1} \\
& \lesssim \vartheta_{q+1}^{-17} \lambda^{2\alpha-3+9 \varepsilon}+\vartheta_{q+1}^{-12} \lambda^{-2 \varepsilon} \\
& \lesssim \vartheta_{q+1}^{-17} \lambda^{-2 \varepsilon}. \label{6.56}\tag{5.54}
\end{align*}
{\bf (iii) Corrector Error}. Taking advantage of H\"{o}lder's inequality, Lemma \ref{F.3}, and \eqref{6.45} , we derive the following estimates:
\begin{align*}
\left\|\mathring{\mathsf{R}}_{\{cor\}}\right\|_{L_{t}^{1} L_{x}^{\varrho}} & \lesssim\left\|w_{q+1}^{(p)} \otimes\left(w_{q+1}^{(c)}+w_{q+1}^{(o)}\right)-\left(w_{q+1}^{(c)}+w_{q+1}^{(o)}\right) \otimes w_{q+1}\right\|_{L_{t}^{1} L_{x}^{\varrho}} \\
& \lesssim\left\|w_{q+1}^{(c)}+w_{q+1}^{(o)}\right\|_{L_{t}^{2} L_{x}^{\infty}}\left(\left\|w_{q+1}^{(p)}\right\|_{L_{t}^{2} L_{x}^{\varrho}}+\left\|w_{q+1}\right\|_{\left.L_{t}^{2} L_{x}^{\varrho}\right)}\right) \\
& \lesssim\left(\vartheta_{q+1}^{-6} \lambda^{-1} r_{\perp}^{-\frac{1}{2}}+\vartheta_{q+1}^{-7} \sigma^{-1}\right)\left(\vartheta_{q+1}^{-1} r_{\perp}^{\frac{1}{\varrho}-\frac{1}{2}}+\vartheta_{q+1}^{-7} \sigma^{-1}\right) \\
& \lesssim\left(\vartheta_{q+1}^{-6} \lambda^{\alpha-2+4 \varepsilon}+\vartheta_{q+1}^{-7} \lambda^{-2 \varepsilon}\right)\left(\vartheta_{q+1}^{-1} \lambda^{-\alpha+1-3 \varepsilon}+\vartheta_{q+1}^{-7} \lambda^{-2 \varepsilon}\right) \\
& \lesssim \vartheta_{q+1}^{-14} \lambda^{-4 \varepsilon} . \label{6.57}\tag{5.55}
\end{align*}

Therefore, from estimates \eqref{6.53}, \eqref{6.56}, and \eqref{6.57}, we conclude that
\begin{align*}
\left\|\mathring{\mathsf{R}}_{q+1}\right\|_{L_{t, x}^{1}} & \leq\left\|\mathring{\mathsf{R}}_{\{l i n\}}\right\|_{L_{t}^{1} L_{x}^{\varrho}}+\left\|\mathring{\mathsf{R}}_{\{osc\}}\right\|_{L_{t}^{1} L_{x}^{\varrho}}+\left\|\mathring{\mathsf{R}}_{\{cor\}}\right\|_{L_{t}^{1} L_{x}^{\varrho}} \\
& \lesssim \vartheta_{q+1}^{-17} \lambda^{- \varepsilon}+\vartheta_{q+1}^{-17} \lambda^{-2 \varepsilon}+\vartheta_{q+1}^{-14} \lambda^{-4 \varepsilon} \\
& \leq \vartheta_{q+1}^{-17}\lambda^{-\varepsilon} \label{6.58}\tag{5.56}
\end{align*}
which verifies the inductive estimate \eqref{2.10} for the $L_{t, x}^{1}$-norm of the new Reynolds stress $\mathring{\mathsf{R}}_{q+1}$.

In the context of the established setting for \eqref{3.1} and the introduction of the index set $\mathscr{C}$, it is worthy noting that the temporal concentrated process of the new Reynolds stress remains independent of spatial dimension. Consequently, the conclusion pertaining to the well-preparedness of $\left(u_{q+1}, \mathring{\mathsf{R}}_{q+1}\right)$ and the temporal inductive inclusion \eqref{2.14} can be directly referenced from the proof provided in the 3D results (\cite{YLZ22}, Sect.7). Based on this, Theorem 2.3 has been rigorously established.

{\centering
\section{Proof of main results}\label{F}}
After the successful establishment of the iterative framework in Proposition \ref{B.3}, the proofs of Theorem \ref{A.2} and Lemma \ref{A.3} became self-evident. The specific proof is briefly described as follows:\\
{\bf Proof of Theorem \ref{A.2}.}  Taking $u_{0}=\tilde{u}$ and set
\begin{align*}
\mathring{\mathsf{R}}_{0} & :=\mathcal{R}\left(\partial_{t} u_{0}+(-\Delta)^{\alpha} u_{0}\right)+u_{0} \mathring{\otimes} u_{0},  \tag{6.1}\\
\mathsf{P}_{0} & :=-\frac{1}{2}\left|u_{0}\right|^{2} , \tag{6.2}
\end{align*}
then the pair $\left(u_{0}, \mathring{\mathsf{R}}_{0}\right)$ constitutes a well-prepared solution to the equations \eqref{2.1}, with the interval $I_{0}=[0, T]$ and the length scale $\theta_{0}=T$. We introduce the notation $\delta_{1}=\|\mathring{\mathsf{R}}_{0}\|_{L_{t,x}^{1}}$ and select the parameter $a$ to be sufficiently large, ensuring that the conditions \eqref{2.6}-\eqref{2.10} are met at the base level $q=0$. Subsequently, relying on Proposition \ref{B.3}, we can establish the existence of a sequence of solutions $(u_{q}, \mathring{\mathsf{R}}_{q})_{q}$ to \eqref{2.1} which satisfies the inductive estimates \eqref{2.6}-\eqref{2.14} for all non-negative integers $q$.

By using the interpolation, \eqref{2.2}, \eqref{2.6}, \eqref{2.7} and \eqref{2.11}, we infer that for any $\beta^{\prime} \in\left(0, \frac{\beta}{7+\beta}\right)$,
\begin{align*}
\sum_{q \geq 0}\left\|u_{q+1}-u_{q}\right\|_{H_{t, x}^{\beta^{\prime}}} & \leq \sum_{q \geq 0}\left\|u_{q+1}-u_{q}\right\|_{L_{t, x}^{2}}^{1-\beta^{\prime}}\left\|u_{q+1}-u_{q}\right\|_{H_{t, x}^{1}}^{\beta^{\prime}} \\
& \lesssim \sum_{q \geq 0}\left(M^{*}\right)^{1-\beta^{\prime}} \delta_{q+1}^{\frac{1-\beta^{\prime}}{2}} \lambda_{q+1}^{7 \beta^{\prime}} \\
& \lesssim\left(M^{*}\right)^{1-\beta^{\prime}} \delta_{1}^{\frac{1-\beta^{\prime}}{2}} \lambda_{1}^{7 \beta^{\prime}}+\sum_{q \geq 1}\left(M^{*}\right)^{1-\beta^{\prime}} \lambda_{q+1}^{-\beta\left(1-\beta^{\prime}\right)+7 \beta^{\prime}}<\infty. \tag{6.3}
\end{align*}

By virtue of \eqref{4.64} and \eqref{6.37}, we have
\begin{equation*}
\sum_{q \geq 0}\left\|u_{q+1}-u_{q}\right\|_{L_{t}^{\gamma} W_{x}^{s, p}}<\infty. \tag{6.4}
\end{equation*}

Therefore, $\left\{u_{q}\right\}_{q \geq 0}$ is established as a Cauchy sequence within the space $H_{t, x}^{\beta^{\prime}}\cap L_{t}^{\gamma} W_{x}^{s, p}$, and thus there exists $u$ such that $\lim _{q \rightarrow \infty}\left(u_{q}\right)=u$ in $H_{t, x}^{\beta^{\prime}}\cap L_{t}^{\gamma} W_{x}^{s, p}$. Additionally, considering the fact that $\lim _{q \rightarrow \infty} \mathring{\mathsf{R}}_{q}=0$ in $L_{t, x}^{1}$ and $u_{q}(0)=\widetilde{u}_{0}$ for all $q \geq 0$, we conclude that $u\in H_{t, x}^{\beta^{\prime}}\cap L_{t}^{\gamma} W_{x}^{s, p}$ is a weak solution to \eqref{1.1} with the initial datum $\widetilde{u}_{0}$. This verifies the assertions $(i)$ and $(ii)$.

As the similar construction in \cite{YLZ22},
$$
\mathscr{G}=\bigcup_{q \geq 0} \mathcal{I}_{q}^{c} \backslash\{0, T\}, \quad \mathscr{B}:=[0, T] \backslash \mathscr{G}
$$
where $\mathscr{B}$ contains the singular set of time, the conclusion $(iii)$ related to the Hausdorff dimension of singular time sets also holds in two-dimensional space.
Concerning the $\varepsilon_{*}$-neighborhood of temporal supports, combining the temporal inductive inclusion \eqref{2.14}
and taking into account $\sum_{q \geq 0} \delta_{q+2}^{1 / 2} \leq \varepsilon_{*}$ for $a$ large enough, then the temporal support statement $(iv)$ is verified. Finally, for the $\varepsilon_{*}$-close between the solution $u$ and the given vector field $\tilde{u}$ in $L_{t}^{1} L_{x}^{2}\cap L_{t}^{\gamma} W_{x}^{s, p}$, by \eqref{2.12} and \eqref{2.13}, $(v)$ is easy to be obtained.

Therefore, the proof of Theorem \ref{A.2} is complete.\\
{\bf Proof of Corollary \ref{A.3}. }Let $\widetilde{u}$ be a weak solution to \eqref{1.1} with divergence-free initial data $\widetilde{u}(0)=\widetilde{u}_{0} \in L^{2}$.

If $\widetilde{u}$ is not a Leray-Hopf solution, then due to the classical result \cite{JL69}, there exists a unique Leray-Hopf weak solution to \eqref{1.1} when $\alpha\geq1$, the non-uniqueness is naturally established. If it is a Leray-Hopf solution, it is already a well-established practice to construct a non-Leray-Hopf solution $u\neq \tilde{u}$ in $ L_{t}^{\gamma} W_{x}^{s, p}$ but with the same initial data $u_{0}=\tilde{u}_{0}$ by using the gluing technique based on Theorem \ref{A.2} $(v)$, even indicating that equations \eqref{1.1} has an infinite number of solutions. Please refer to \cite{AC22,YLZ22} for the proof of relevant details.
\\{\bf Proof of Corollary \ref{A.4}}.
For any $\alpha \in[1,\frac{3}{2})$, the embedding relationship on $\mathbb{T}^{2}$: $L_{x}^{2} \hookrightarrow L_{x}^{2 /(2 \alpha-1)} \hookrightarrow L_{x}^{p}$ implies that the initial data $u_{0}$ is in $L_{x}^{2}$ and also belonging to $L_{x}^{p}$ in Theorem \ref{A.2}, by taking $s=0$ and $\gamma=\infty$ in the super-critical regime $\mathcal{A}_{1}$, we obtain the non-uniqueness of solutions in the super-critical spaces $L_{x}^{p}$ for any $1 \leq p < \frac{2}{2\alpha - 1}$ which verifies $(i)$.

Considering ($ii$, $iii$), according to the embedding theorems in \cite{HS87}, by choosing small constant $\epsilon>0$ and defining $s_{p}:=2 / p+1-2 \alpha$, $\bar{p}:=\frac{2}{2 \alpha-1+\epsilon}$ such that $s<s_{p}-\epsilon$, $1<\bar{p}<\frac{2}{2 \alpha-1}$, we have
$$
L^{\bar{p}}=F_{\bar{p}, 2}^{0} \hookrightarrow F_{p, q}^{s_{p}-\epsilon} \hookrightarrow B_{p, p \vee q}^{s_{p}-\epsilon} \hookrightarrow B_{p, q}^{s},
$$
and $
L^{\bar{p}} \hookrightarrow F_{p, q}^{s_{p}-\epsilon}\hookrightarrow F_{p, q}^{s}$, then Theorem \ref{A.2} also implies that the non-uniqueness of solutions in the super-critical Besov and Triebel-Lizorkin spaces for any $s<2 / p+1-2 \alpha$ which verifies ($ii, iii$).

The proof is therefore finished.\\

\begin{appendices}
{\centering
\section{Analytic tools in convex integration scheme}\label{G1}}

In this section, we shall exhibit some preliminary results that have been previously utilized in the preceding sections.

\begin{Lemma}\label{G.1} (\cite{ACL22}, Geometric Lemma 3.1). There exists a set $\Lambda \subset \mathbb{S} \cap \mathbb{Q}^{2}$ consisting vectors $k$, and positive smooth functions $\gamma_{(k)}: B_{C_{\mathsf{R}}}(\mathrm{Id}) \rightarrow \mathbb{R}$, such that for $\mathsf{R} \in B_{C_{\mathsf{R}}}(Id)$, it holds that
\begin{equation*}
\mathsf{R}=\sum_{k \in \Lambda} \gamma_{(k)}^{2}(\mathsf{R}) k_{1} \otimes k_{1},\tag{A.1}
\end{equation*}
where $(k, k_{1})$ is an orthonormal basis and $B_{C_{\mathsf{R}}}(\mathrm{Id})$ represents the ball with a radius of $C_{\mathsf{R}}>0$ centered on the identity matrix in the $2 \times 2$ symmetric matrix space.
\end{Lemma}

As stated in the previous work \cite{BCV22}, there exists  $N_{\Lambda}\in \mathbb{N}$ that fulfills
\begin{equation*}
\left\{N_{\Lambda} k, N_{\Lambda} k_{1}\right\} \subseteq N_{\Lambda} \mathbb{S} \cap \mathbb{Z}^{2}. \label{7.2}\tag{A.2}
\end{equation*}
We introduce the notation $M_{*}$ to denote a geometric constant that satisfies the following inequality:
\begin{equation*}
\sum_{k \in \Lambda}\left\|\gamma_{(k)}\right\|_{C^{4}\left(B_{C_{\mathsf{R}}}(\mathrm{Id})\right)} \leq M_{*}. \tag{A.3}
\end{equation*}

Subsequently, we revisit the $L^{p}$ decorrelation lemma, introduced by Lemma 2.4 in \cite{AC21} (see also \cite{BV19}, Lemma 3.7). This lemma serves as a crucial component in obtaining the $L_{t, x}^{2}$ estimates for the perturbations.
\begin{Lemma}\label{H.2}(\cite{AC21}, Lemma 2.4). For any $p \in[1, \infty]$, there exist $\sigma \in \mathbb{N}$ and smooth functions $f, g: \mathbb{T}^{d} \rightarrow \mathbb{R}$ such that,
\begin{equation*}
\left|\|f g(\sigma \cdot)\|_{L^{p}\left(\mathbb{T}^{d}\right)}-\|f\|_{L^{p}\left(\mathbb{T}^{d}\right)}\|g\|_{L^{p}\left(\mathbb{T}^{d}\right)}\right| \lesssim \sigma^{-\frac{1}{p}}\|f\|_{C^{1}\left(\mathbb{T}^{d}\right)}\|g\|_{L^{p}\left(\mathbb{T}^{d}\right)}. \tag{A.4}
\end{equation*}
\end{Lemma}

The following stationary phase lemma, as a primary tool, plays a pivotal role in addressing the errors associated with Reynolds stress.
\begin{Lemma}\label{H.3} (\cite{LTP20}, Lemma 7.4). For any given $1<p<\infty$, $\lambda\in \mathbb{Z}_{+}$, $a \in C^{2}\left(\mathbb{T}^{2}, \mathbb{R}\right)$ and $f\in L^{p}\left(\mathbb{T}^{2}, \mathbb{R}^{2}\right)$, one has
$$
\left\||\nabla|^{-1} \mathbb{P}_{\neq 0}\left(a \mathbb{P}_{\geq \lambda} f\right)\right\|_{L^{p}\left(\mathbb{T}^{2}\right)} \lesssim \lambda^{-1}\left\|\nabla^{2} a\right\|_{L^{\infty}\left(\mathbb{T}^{2}\right)}\|f\|_{L^{p}\left(\mathbb{T}^{2}\right)}
.$$
\end{Lemma}

{\centering
\section{$\mathbf{X}^{\mathbf{s},\mathbf{\gamma},\mathbf{p}}$ weak solutions on the torus}\label{G2}}
In this section, we present that the sub-critical and critical weak solutions $L_{t}^{\gamma} W_{x}^{s, p}$ to \eqref{1.1}, where the triple $(s, \gamma, p)$ satisfy \eqref{1.41}, are in fact Leray-Hopf weak solutions by duality approach in \cite{PL01}. In view of the celebrated result \cite{JL69} when $\alpha\geq1$ in 2D, these Leray-Hopf weak solutions with the same initial data must coincide.

For convenience, we define the mixed Lebesgue space $X^{s, \gamma, p}([0, T];\mathbb{T}^{2})$ satisfying $\frac{2 \alpha}{\gamma}+\frac{2}{p}\leq2 \alpha-1+s$, for some $\gamma, p\in [1,\infty]$, $s \geq 0$, as
\begin{align*}
X^{s, \gamma, p}([0, T];\mathbb{T}^{2})=\left\{\begin{array}{l}
L^{\gamma}\left(0, T ; W_{x}^{s, p}(\mathbb{T}^{2})\right), \quad \gamma\neq\infty,  \tag{B.1}\\
C\left([0, T] ; W_{x}^{s, p}(\mathbb{T}^{2})\right), \quad \gamma=\infty.
\end{array}\right.
\end{align*}

The main result is formulated in Theorem \ref{I.1} below.

\begin{Theorem}\label{I.1} (Uniqueness of weak solutions). Let $u\in X^{s, \gamma, p} ([0, T];\mathbb{T}^{2})$ and $(s, \gamma, p)$ satisfy \eqref{1.41} with $\alpha \geq 1$, $s \geq 0$. For some $1 \leq \gamma \leq \infty$, $1 \leq p \leq\infty$ and $0<\frac{1}{p}-\frac{s}{2} < \frac{1}{2}$, if $u$ is a weak solution to \eqref{1.1} in the sense of Definition \ref{A.1}, then $u$ is the unique Leray-Hopf solution.
\end{Theorem}

As in the context of NSE \cite{AC22,YL22}, we recall the following existence result for the linearized generalized NSE, which can be proved by a standard Galerkin method.

\begin{Lemma}\label{I.2} Let $u \in X^{s, \gamma, p} ([0, T];\mathbb{T}^{2})$ be a weak solution of \eqref{1.1} and $(s, \gamma, p)$ satisfy \eqref{1.41} with $\alpha \geq 1$, $s \geq 0$, $1 \leq \gamma \leq \infty$, $1 \leq p \leq\infty$ and $0 < \frac{1}{p}-\frac{s}{2} <\frac{1}{2}$. For any divergence-free initial data $v_{0} \in L^{2}\left(\mathbb{T}^{2}\right)$, there exist a weak solution $v\in C_{w}\left([0, T] ; L_{x}^{2}\right) \cap L^{2}\left(0, T ; H_{x}^{\alpha}\right)$ to the linearized NSE:
\begin{align*}
\left\{\begin{array}{l}
\partial_{t} v+(-\Delta)^{\alpha} v+u \cdot \nabla v+\nabla P=0, \label{B.31}\tag{B.2}\\
\operatorname{div} v=0,  \\
v|_{t=0}=v_{0},
\end{array}\right.
\end{align*}
which satisfies the energy inequality
$$
\frac{1}{2}\|v(t)\|_{L_{x}^{2}}^{2}+ \int_{t_{0}}^{t}\left\|\Lambda^{\alpha} v(s) \right\|_{L_{x}^{2}}^{2} d s \leq \frac{1}{2}\left\|v_{0}\right\|_{L_{x}^{2}}^{2}
$$
for all $t \in\left[t_{0}, T\right]$, a.e. $t_{0} \in[0, T]$ (including $t_{0}=0$).
\end{Lemma}

Let $v$ be the weak solution given by Lemma \ref{I.2} with initial data $u_{0}$ and setting $w=u-v$, we derive from \eqref{B.31} that
\begin{align*}
\left\{\begin{array}{l}
\partial_{t} w+(-\Delta)^{\alpha} w+u \cdot \nabla w+\nabla q=0,  \tag{B.3}\\
\operatorname{div} w=0, \\
w_{0}=0,
\end{array}\right.
\end{align*}
and its weak formulation
\begin{align*}
\int_0^T \int_{\mathbb{T}^2} w \cdot\left(\partial_t \mathscr{T}-(-\Delta)^{\alpha} \mathscr{T}+u \cdot \nabla \mathscr{T}\right) d x d t=0, \quad \text { for any } \mathscr{T} \in \mathcal{D}_T, \label{B.13}\tag{B.4}
\end{align*}
where $\mathcal{D}_T$ represents the test function class containing the smooth divergence-free functions vanishing for $t \geq T$.

Fix $F \in C_c^{\infty}\left([0, T] \times \mathbb{T}^2\right)$. Let $\Theta:[0, T] \times \mathbb{T}^2 \rightarrow \mathbb{R}^2$ and $\Pi:[0, T] \times \mathbb{T}^2 \rightarrow$ $\mathbb{R}$ satisfy
\begin{align*}
\left\{\begin{array}{l}
-\partial_t \Theta+(-\Delta)^{\alpha} \Theta-u \cdot \nabla \Theta+\nabla \Pi=F, \\
\operatorname{div} \Theta=0, \\
\Theta(T)=0.\label{B.4}\tag{B.5}
\end{array}\right.
\end{align*}
Arguing as in Appendix A of \cite{AC22}, the problem of whether $u \equiv v$ reduces to proving that $\Theta$ satisfies some regularity, then $\Theta$ can be used the test function in the weak formulation \eqref{B.3}. Specifically, we give the following lemma.
\begin{Lemma}\label{I.3} Let $u \in X^{s, \gamma, p} ([0, T];\mathbb{T}^{2})$ be a weak solution of \eqref{1.1} and $(s, \gamma, p)$ satisfy \eqref{1.41} with $\alpha \geq 1$, $s \geq 0$, $1 \leq \gamma \leq \infty$, $1 \leq p \leq\infty$ and $0 < \frac{1}{p}-\frac{s}{2} < \frac{1}{2}$. For any $F \in C_c^{\infty}([0, T] \times$ $\mathbb{T}^2$), the system \eqref{B.4} has a weak solution $\Theta \in L_t^{\infty} L^2 \cap L_t^2 H^\alpha$ such that $\Theta$ can be used as a test function in \eqref{B.3}.
\end{Lemma}
{\bf Proof:} We will prove the weak solution $\Theta$ satisfies the regularity
\begin{align*}
\partial_t \Theta,\ (-\Delta)^{\alpha} \Theta,\ u \cdot \nabla \Theta,\ \nabla \Pi \in L_{t}^{\gamma_{1}}L_{x}^{p_{1}},\label{B.15}\tag{B.6}
\end{align*}
which implies $\Theta$ can be used in \eqref{B.3} since $w \in L_{t}^{\gamma_{2}}L_{x}^{p_{2}}$ with $\frac{\alpha}{\gamma_{2}}+\frac{1}{p_{2}}=1$.

The solution $\Theta$ can be constructed by the Galerkin method, so it is necessary to provide a priori estimates as following:\\
{\bf Step 1:} $L_t^{\infty} L^2 \cap L_t^2 H^\alpha-$estimates of $\Theta$.

By the classical energy method, we have
\begin{align*}
-\frac{1}{2} \frac{d}{d t}\|\Theta\|_{L^{2}}^2+\|\Lambda^{\alpha} \Theta\|_{L^{2}}^2 \leq \int_{\mathbb{T}^2}|F \cdot \Phi| d x, \label{B.5}\tag{B.7}
\end{align*}
which implies the desired regularity estimates.\\
{\bf Step 2:} $L_t^{\infty} H^1 \cap L_t^2 H^{\alpha+1}-$estimates of $\Theta$.

Taking the $L^2-$inner product of \eqref{B.4} with $\Delta \Theta$ to obtain
$$
-\frac{1}{2} \frac{d}{d t}\|\nabla \Theta\|_{L^{2}}^2+\|\Lambda^{\alpha+1} \Theta\|_{L^{2}}^2 \leq \int_{\mathbb{T}^2}|F \cdot \Delta \Theta| d x+\int_{\mathbb{T}^2}|u \cdot \nabla \Theta \cdot \Delta \Theta| d x.
$$

{\bf $\bullet$ Case\ 1:} $\gamma<\infty \text { and } 1\leq p < \infty$.

In view of \eqref{B.19} in Proposition \ref{I.4}, for any $2<q<\infty$,
\begin{align*}
\int_{\mathbb{T}^2}|u \cdot \nabla \Theta\cdot \Delta \Theta| d x &\lesssim\|u\|^{\frac{2}{\theta+\delta}}_{L^{q}}
\|\nabla \Theta\|_{L^{2}}^{2}+\left\|\Lambda^{\alpha+1} \Theta\right\|_{L^{2}}^{2} \\
&\lesssim\|u\|^{\frac{2}{\theta+\delta}}_{W^{s,p}}
\|\nabla \Theta\|_{L^{2}}^{2}+\left\|\Lambda^{\alpha+1} \Theta\right\|_{L^{2}}^{2},
\end{align*}
where we have used Sobolev embedding theorem $W^{s,p}(\mathbb{T}^{2})\hookrightarrow L^{q}(\mathbb{T}^{2})$ with $\frac{1}{q}=\frac{1}{p}-\frac{s}{2}$.

Then H\"{o}lder's and Young's inequalities yield
$$
-\frac{1}{2} \frac{d}{d t}\|\nabla \Theta\|_{L^{2}}^2+\frac{1}{2}\|\Lambda^{\alpha+1} \Theta\|_{L^{2}}^2 \lesssim\|u\|^{\frac{2}{\theta+\delta}}_{W^{s,p}}
\|\nabla \Theta\|_{L^{2}}^2+\|F\|_{L^{2}}^2 .
$$
Thanks to the integrability of $\|u\|^{\frac{2}{\theta+\delta}}_{W^{s,p}}$, where
$$
\frac{2}{\theta+\delta}=\frac{2 \alpha}{2 \alpha-1-2 / p+s} = \gamma,
$$
and Gronwall's inequality immediately implies that $\Theta \in L_t^{\infty} H^1 \cap L_t^2 H^{\alpha+1}$.

{\bf $\bullet$ Case\ 2:} One endpoint case  $\gamma=\frac{2}{2\alpha-1+s} \text { and } p=\infty$.\vspace{2mm}

By \eqref{B.110}, and let $p=q=\infty$, thanks to the integrability of $\|u\|^{\frac{2\alpha}{2\alpha-1}}_{L^{\infty}}$ when $s=0$, we can verify that $\Theta \in L_t^{\infty} H^1 \cap L_t^2 H^{\alpha+1}$.

When $s>0$, due to the fact \eqref{B.111} and $u \in X^{s, \gamma, p} ([0, T];\mathbb{T}^{2})$, and by Holder's inequality
$$
\left|\int_{0}^{T}\|u\|^{\frac{2\alpha s}{(2\alpha-1)(1+s)}}_{L^{2}}\|u\|^{\frac{2\alpha}
{(2\alpha-1)(1+s)}}_{W^{s,\infty}}dt\right|\lesssim
\|u\|_{L_{t}^{\frac{2\alpha}{2\alpha-2}}L^{2}_{x}}
\|u\|_{X^{s,\frac{2\alpha}{2\alpha-1+s},\infty}}
\lesssim C,$$
where $L_{t}^{\frac{2\alpha}{2\alpha-2}}L^{2}_{x}$ is also belonging to the class $X^{s,\gamma,p}$, which implies that $\Theta \in L_t^{\infty} H^1 \cap L_t^2 H^{\alpha+1}$.

{\bf $\bullet$ Case\ 3:} The other endpoint case $\gamma=\infty \text { and } p=\frac{2}{2\alpha-1+s}$.\vspace{2mm}

Since $u \in C([0,T]; W_x^{s, \frac{2}{2 \alpha-1+s}})$, for any $\varepsilon>0$, we can decompose $u=u_{1}+u_{2}$ (see \cite{TT18} for example) such that
$$\left\|u_1\right\|_{X^{s, \infty, \frac{2}{2 \alpha-1+s}}}\lesssim\epsilon\quad \text{and}\quad u_2\in X^{s, \infty,\infty}.$$
Then, H\"{o}lder's inequality, Sobolev's embedding and Young's inequality yield
\begin{align*}
&\left|\int_{0}^{T} \int_{\mathbb{T}^{2}} u \cdot\nabla\Theta \cdot \Delta \Theta d x d t
\right|\\&  \lesssim\left\|u_{1}\right\|_{X^{s, \infty, \frac{2}{2 \alpha-1+s}}}\|\nabla\Theta \cdot \Delta \Theta\|_{L^{1}_{t}L_{x}^{\frac{2}{3-2\alpha}}}+\left\|u_{2}\right\|_{X ^{s, \infty, \infty}}\|\nabla\Theta \cdot \Delta\Theta\|_{L^{1}_{t} L_{x}^{1}} \\
& \lesssim \epsilon\|\nabla\Theta\|_{L^{2}_{t} L_{x}^{\frac{2}{1-\alpha}}}\|\Delta \Theta\|_{L^{2}_{t} L_{x}^{\frac{2}{2-\alpha}}}+\left\|u_{2}\right\|_{X^{s, \infty, \infty}}\|\nabla\Theta\|_{L^{2}_{t}L_{x}^{2}}
\|\Lambda^{\alpha+1}\Theta\|_{L^{2}_{t} L_{x}^{2}} \\
& \lesssim \epsilon\left\|\Lambda^{\alpha+1} \Theta\right\|_{L^{2}_{t} L_{x}^{2}}^{2}+C_{\epsilon}\|\nabla\Theta\|_{L^{2}_{t}L_{x}^{2}}^{2}, \tag{B.8}
\end{align*}
we arrive at the same conclusion.

{\bf Step 3:} Conclusion from maximal regularity of the fractional parabolic equation.

Applying the projection operator $\mathbb{P}$ to \eqref{B.4}, the equation for $\Theta$ can be rewritten as
\begin{align*}
\left\{\begin{array}{l}
-\partial_t \Theta+(-\Delta)^{\alpha} \Theta=\mathbb{P}(u \cdot \nabla \Theta)+\mathbb{P} F, \\
\Theta(T)=0,
\end{array}\right.
\end{align*}
when  $\alpha\geq1$. Reference to the maximal $L_{t}^{\gamma}L_{x}^{p}$ regularity of heat equation is also applicable to fractional parabolic equation, so we only need to show the estimate \eqref{B.5} for $u \cdot \nabla \Theta$ to conclude the proof. By Sobolev embedding theorem $W^{s,p}(\mathbb{T}^{2})\hookrightarrow L^{q}(\mathbb{T}^{2})$ with $\frac{1}{q}=\frac{1}{p}-\frac{s}{2}$ and Sobolev interpolation, we have
\begin{align*}
\|u \cdot \nabla \Theta\|_{L^{\gamma_{1}}\left([0, T]; L^{p_{1}} (\mathbb{T}^2)\right)}&\lesssim \|u\|_{L_{t}^{\gamma}L_{x}^{q}}\|\nabla\Theta\|_{L_{t}^{\gamma_{3}} L_{x}^{q_{3}}}\\
&\lesssim\|u\|_{L_{t}^{\gamma}W_{x}^{s,p}}\|\nabla\Theta\|_{L_{t}^{\infty} L_{x}^{2}}^{1-\xi}\|\nabla\Theta\|_{L_{t}^{\gamma_{3}} H_{x}^{\alpha}}^{\xi},
\end{align*}
where $\frac{2\alpha}{\gamma_{3}}+\frac{2}{q_{3}}=1$ when $\xi=\frac{2}{\gamma_{3}}\in[0,1]$. Recalling that $u \in X^{s,\gamma,p}$ with $\frac{2\alpha}{\gamma}+\frac{2}{p}\leq2\alpha-1+s$, and combing the conclusion of Step 2, we conclude that
$$
u \cdot \nabla \Theta \in L^{\gamma_{1}}\left([0, T]; L^{p_{1}} (\mathbb{T}^2)\right)
$$
with $\frac{\alpha}{\gamma_{1}}+\frac{1}{p_{1}}\leq\alpha$. After the duality discussion, then $\Theta$ can be chosen to the be test function in \eqref{B.3}.

The proof of Lemma \ref{I.3} is complete.
\begin{Proposition}\label{I.4}  Let $2<q\leq\infty,\ s\geq0$. For any smooth vector fields $u,v\in C_0^{\infty}\left(\mathbb{T}^2\right)$, when $2<q<\infty,$
\begin{align*}
\int_{\mathbb{T}^{2}} \left|u \cdot\nabla v \cdot \Delta v\right| d x\lesssim\|u\|^{\frac{2}{\theta+\delta}}_{L^{q}}
\|\nabla v\|_{L^{2}}^{2}+\left\|\Lambda^{\alpha+1} v\right\|_{L^{2}}^{2}.\label{B.19}\tag{B.9}
\end{align*}
In particularly, when $q=\infty, s=0$, we have
\begin{align*}
\int_{\mathbb{T}^{2}} \left|u \cdot\nabla v \cdot \Delta v\right| d x\lesssim\|u\|^{\frac{2\alpha}{2\alpha-1}}_{L^{\infty}}
\|\nabla v\|_{L^{2}}^{2}+\left\|\Lambda^{\alpha+1} v\right\|_{L^{2}}^{2},\label{B.110}\tag{B.10}
\end{align*}
and when $s>0$,
\begin{align*}
\int_{\mathbb{T}^{2}} \left|u \cdot\nabla v \cdot \Delta v\right| d x\lesssim\|u\|^{\frac{2\alpha s}{(2\alpha-1)(1+s)}}_{L^{2}}\|u\|^{\frac{2\alpha}
{(2\alpha-1)(1+s)}}_{W^{s,\infty}}
\|\nabla v\|_{L^{2}}^{2}+\left\|\Lambda^{\alpha+1} v\right\|_{L^{2}}^{2}.\label{B.111}\tag{B.11}
\end{align*}
\end{Proposition}
{\bf Proof:}
When $2<q<\infty$, we apply the H\"{o}lder's inequality, Young inequality, Gagliardo-Nirenberg inequality to derive that
\begin{align*}
& \int_{\mathbb{T}^{2}} \left|u \cdot\nabla v \cdot \Delta v\right| d x \\
& \lesssim\|u\|_{L^q}\|\nabla v\|_{L^b}\|\Delta v\|_{L^c}\\
& \lesssim\|u\|_{L^{q}}\|\nabla v\|_{L^{2}}^{\theta}
\|\Lambda^{\alpha+1} v\|_{L^{2}}^{1-\theta}
\|\nabla v\|_{L^{2}}^{\delta}
\|\Lambda^{\alpha+1}v\|_{L^{2}}^{1-\delta} \\
& \lesssim\|u\|^{\frac{2}{\theta+\delta}}_{L^{q}}\|\nabla v\|_{L^{2}}^{2}+\left\|\Lambda^{\alpha+1} v\right\|_{L^{2}}^{2}. \label{B.9}\tag{B.12}
\end{align*}
The constants $1<b, c<\infty$ and $0 \leqslant \theta, \delta \leqslant 1$ satisfy
\begin{align*}
\left\{\begin{array}{l}
\frac{1}{b}+\frac{1}{c}+\frac{1}{q}=1,\vspace{2mm} \\
\frac{1}{b}-\frac{1}{2}=(1-\theta)\left(\frac{1}{2}-\frac{\alpha+1}{2}\right), \vspace{2mm}\\
\frac{1}{c}-1=(1-\delta)\left(\frac{1}{2}-\frac{\alpha+1}{2}\right) .
\end{array}\right.
\label{B.10}\tag{B.13}
\end{align*}
System \eqref{B.10} has 4 unknowns but 3 equations, so we choose one solution to \eqref{B.10} as
$$
b=\frac{4 q}{q-2}, \quad c=\frac{4q}{3 q-2}, \quad \theta=\delta=1-(\frac{1}{2\alpha}+\frac{1}{\alpha q})<1.
$$

In particularly, when $q=\infty,\ s=0$, it yields that
$$
\begin{aligned}
&\quad\quad \int_{\mathbb{T}^{2}} \left|u \cdot\nabla v \cdot \Delta v\right| d x \\
& \quad \lesssim \left\|u\right\|_{L^{\infty}}\|\nabla v\|_{L^{4}} \|\Delta v\|_{L^{\frac{4}{3}}} \\
& \quad \lesssim \left\|u\right\|_{L^{\infty}} \left\|\Lambda^{\alpha+1} v\right\|_{L^{2}}^{ \frac{1}{2\alpha}}\|\nabla v\|_{L^{2}}^{1-\frac{1}{2\alpha}}\left\|\Lambda^{\alpha+1} v\right\|_{L^{2}}^{ \frac{1}{2\alpha}}\|\nabla v\|_{L^{2}}^{1-\frac{1}{2\alpha}} \\
& \quad \lesssim \left\|u\right\|_{L^{\infty}}\|\nabla v\|_{L^{2}}^{\frac{2\alpha-1}{\alpha}}\left\|\Lambda^{\alpha+1} v\right\|_{L^{2}}^{ \frac{1}{\alpha}}\\
& \quad \lesssim  \left\|u\right\|_{L^{\infty}}^{\frac{2\alpha}{2\alpha-1}}\|\nabla v\|_{L^{2}}^{2}+\left\|\Lambda^{\alpha+1} v\right\|_{L^{2}}^{2},
\end{aligned}
$$
then \eqref{B.110} is verified.

When $q=\infty, \ s>0$, by Gagliardo-Nirenberg inequality, we deduce that
$$\left\|u\right\|_{L^{\infty}}
\leq\|u\|^{\frac{s}{1+s}}_{L^{2}}\|u\|^{\frac{1}
{1+s}}_{W^{s,\infty}},$$
which implies \eqref{B.111}.

{\centering
\section{Regularity and stability estimates in hyper-dissipative case}\label{G3}}
It is well-established that the NSE is global (in time) well-posed in $H^{3}(\mathbb{T}^{2})$, which serves as a scaling sub-critical space. For $\alpha\in[1,\frac{3}{2})$, we have the similar regularity estimates of strong solution, facilitating the temporal gluing of local solutions in Section \ref{B}.
\begin{Proposition}\label{C.2}
$\alpha\in[1,\frac{3}{2})$. Let $v_{0}=\left.v\right|_{t=t_{0}} \in H^{3}\left(\mathbb{T}^{2}\right)$ with zero mean on $\mathbb{T}^{2}$, the Cauchy problem \eqref{1.1} admits a unique strong solution on the interval $\left[t_{0}, t_{*}\right)$ that satisfies
\begin{align*}
& \sup _{t \in\left[t_{0}, t_{*}\right]}\|v(t)\|_{L^{2}}^{2}+2 \int_{t_{0}}^{t_{*}}\|v(t)\|_{\dot{H}^{\alpha}}^{2} d t \leq\left\|v_{0}\right\|_{L^{2}}^{2}, \label{3.8a}\tag{C.1a}\\
& \sup _{t \in\left[t_{0}, t_{*}\right]}\|v(t)\|_{H^{3}} \leq 2\left\|v_{0}\right\|_{H^{3}}. \label{3.8b}\tag{C.1b}
\end{align*}

Moreover, assuming that
\begin{equation*}
0<t_{*}-t_{0} \leq \frac{1}{\left\|v_{0}\right\|_{H^{3}}\left(1+\left\|v_{0}\right\|_{L^{2}}\right)^{\frac{1}{2 \alpha-1}}}, \label{3.9}\tag{C.2}
\end{equation*}
we have
\begin{equation*}
\sup _{t \in\left(t_{0}, t_{*}\right]}\left|t-t_{0}\right|^{\frac{N}{2 \alpha}+M}\left\|\partial_{t}^{M} D^{N} v(t)\right\|_{H^{3}} \lesssim\left\|v_{0}\right\|_{H^{3}} \label{3.10}\tag{C.3}
\end{equation*}
for any $N \geq 0$ and $M \in\{0,1\}$. The implicit constants only depend on $\alpha, N, M$.
\end{Proposition}
{\bf Proof.} The energy inequality gives a global in time control on $\|v(t)\|_{L^{2}}$:
$$
\frac{1}{2} \frac{d}{d t}\|v\|_{L^{2}}^{2} \leq-\|v\|_{\dot{H}^{\alpha}}^{2},
$$
and thanks to the fact $\int_{\mathbb{T}^{2}}(v\cdot\nabla v)\cdot \Delta v dx=0$ , we have
\begin{equation*}
\frac{1}{2} \frac{d}{d t}\|v\|_{\dot{H}^{1}}^{2} \leq-\|v\|_{\dot{H}^{\alpha+1}}^{2}.
\label{3.11}\tag{C.4}
\end{equation*}
From the commutator estimates in \cite{TK88}, the Young inequality, and using the divergence free condition, we obtain:
\begin{equation*}
\frac{1}{2} \frac{d}{d t}\|v\|_{\dot{H}^{3}}^{2}+\|v\|_{\dot{H}^{3+\alpha}}^{2} \leq\frac{1}{2}\|v\|_{\dot{H}^{3}}^{2}\| v\|_{\dot{H}^{2}}^{2}+\frac{1}{2}\|v\|_{\dot{H}^{4}}^{2}.
\label{3.12}\tag{C.5}
\end{equation*}
Combining \eqref{3.11} for $\alpha\in[1,\frac{3}{2})$ in bounded domain, we conclude that \eqref{3.8b} holds true by the Gronwall inequality.
The bound \eqref{3.8b} exhibits a sub-critical property, as an $L_{t}^{\infty} H_{x}^{3}$ apriori estimate alone suffices to establish the uniqueness of the solution. The higher regularity claimed in \eqref{3.10} could potentially be derived from the proof of Proposition 3.1 in \cite{BCV22}, by using  two-dimensional Gagliardo-Nirenberg-Sobolev inequalities. Then, the proof of Proposition \ref{C.2} is finished.

Define  the vector field $\omega_{i}:\left[t_{i}, t_{i+1}\right] \times \mathbb{T}^{2} \rightarrow \mathbb{R}^{2}$ as $\omega_{i}:=u_{q}-v_{i}$ and the scalar field $p_{i}:\left[t_{i}, t_{i+1}\right] \times \mathbb{T}^{2} \rightarrow \mathbb{R}$ such that satisfy the following equations:
\[
\left\{\begin{array}{l}
\partial_{t} \omega_{i}+(-\Delta)^{\alpha} \omega_{i}+\operatorname{div}\left(v_{i} \otimes \omega_{i}+\omega_{i} \otimes u_{q}\right)+\nabla p_{i}=\operatorname{div} \mathring{\mathsf{R}}_{q}  \label{3.14}\tag{C.6}\\
\operatorname{div} \omega_{i}=0, \\
\left.\omega_{i}\right|_{t=t_{i}}=0.
\end{array}\right.
\]

We recall the stability estimates for the solutions $\omega_{i}$ to \eqref{3.14} in \cite{YL22}, the proof therein is applicable to the 2D case as well, due to the embedding $H^{3}(\mathbb{T}^{2})\hookrightarrow W^{1,\infty}(\mathbb{T}^{2})\hookrightarrow L^{\infty}(\mathbb{T}^{2})$. Therefore, we omit the detailed proof for brevity.

\begin{Lemma}\label{C.3} (\cite{YL22}, Lemma 3.3). For $\alpha \in[1,\frac{3}{2}), 1<\varrho \leq 2$. Let $\left(u_{q}, \mathring{\mathsf{R}}_{q}\right)$ be the well-prepared solution to \eqref{2.1}, the solutions $\omega_{i}$ to \eqref{3.14} hold that
\begin{align*}
& \left\|\omega_{i}\right\|_{L^{\infty}\left(\left[t_{i}, t_{i+1}+\vartheta_{q+1}\right] ; L_{x}^{\varrho}\right)} \leq C \int_{t_{i}}^{t_{i+1}+\vartheta_{q+1}}\left\||\nabla| \mathring{\mathsf{R}}_{q}(s)\right\|_{L_{x}^{\varrho}} d s, \label{3.15} \tag{C.7}\\
& \left\|\omega_{i}\right\|_{L^{\infty}\left(\left[t_{i}, t_{i+1}+\vartheta_{q+1}\right] ; H_{x}^{3}\right)} \leq C \int_{t_{i}}^{t_{i+1}+\vartheta_{q+1}}\left\||\nabla| \mathring{\mathsf{R}}_{q}(s)\right\|_{H_{x}^{3}} d s,  \label{3.16}\tag{C.8}\\
& \left\|\mathcal{R} \omega_{i}\right\|_{L^{\infty}\left(\left[t_{i}, t_{i+1}+\vartheta_{q+1}\right] ; L_{x}^{\varrho}\right)} \leq C \int_{t_{i}}^{t_{i+1}+\vartheta_{q+1}}\left\|\mathring{\mathsf{R}}_{q}(s)\right\|_{L_{x}^{\varrho}} d s, \label{3.17}\tag{C.9}
\end{align*}
where $\mathcal{R}$ is the inverse divergence operator defined as in \cite{CD13}, and the constants $C$ only depend on $\alpha$ and $\varrho$.
\end{Lemma}

\end{appendices}
\section*{Acknowledgment.}
Lili Du is supported by National Natural Science Foundation of China Grant 12125102. Xinliang Li is supported by National Natural Science Foundation of
China Grant 12401297 and the project funded by China Postdoctoral Science Foundation 2024T170577, 2023M742401.
\section*{Conflicts of interest}\ \ \ \
The authors declare that they have no conflict of interest.
\section*{Data Availability}\ \ \ \
Data sharing is not applicable to this article as no data were created or analyzed in this study.

\begin{center}

\end{center}
\end{document}